\documentclass[a4paper,10pt]{article}

\usepackage[title]{appendix}
\usepackage{setspace}
\usepackage[bottom]{footmisc}
\usepackage{color}
\usepackage{subfig}
\usepackage{svg}
\usepackage{hyperref}
\usepackage{multirow}
\usepackage{cite}
\usepackage{makecell}
\usepackage{amsmath} 
\usepackage{mathtools}
\usepackage{cleveref}
\usepackage{graphicx}
\usepackage{epsfig}
\usepackage{epstopdf}
\usepackage{float}
\usepackage{titling}
\usepackage[margin=1.0in]{geometry}
\usepackage{etoolbox}

\usepackage{standalone}     

\usepackage{tikz}
\usepackage{tkz-base}       
\usepackage{tkz-euclide}

\usepackage{xcolor}
\usepackage{color, colortbl}
\usepackage{enumitem}
\usepackage{empheq}
\usepackage{siunitx}  
\usepackage{booktabs} 

\newcommand{\volume}{{\ooalign{\hfil$V$\hfil\cr\kern0.08em--\hfil\cr}}}
\newcommand{\bcb}[1]{\fbox{\color{blue} #1}}
\newcommand{\bcr}[1]{\fbox{\color{red} #1}}

\usepackage{algorithm}
\usepackage{algpseudocode}

\title{Application of troubled-cells to finite volume methods - an optimality study using a novel monotonicity parameter}
\makeatletter
\renewcommand\@date{{%
  \vspace{-\baselineskip}%
  \large\centering
  \begin{tabular}{@{}c@{}}
    R. Shivananda Rao\textsuperscript{1} \\
    \normalsize nandu.rsr@gmail.com
  \end{tabular}%
  \quad and\quad
  \begin{tabular}{@{}c@{}}
    M Ramakrishna\textsuperscript{2} \\
    \normalsize krishna@ae.iitm.ac.in
  \end{tabular}

  \bigskip

  \textsuperscript{1}Dept. of Aerospace Engg., IIT Madras.\par
  \textsuperscript{2}Professor, Dept. of Aerospace Engg., IIT Madras.

  \bigskip

  \today
}}
\makeatother
\begin{document}

\maketitle

\begin{abstract}
We adapt a troubled-cell indicator from discontinuous Galerkin (DG) methods to finite volume methods (FVM) with MUSCL reconstruction and using a novel monotonicity parameter show there is a trade-off between convergence and quality of the solution. Employing two dimensional compressible Euler equations for flows with oblique shocks, this trade-off is studied by varying the number of troubled-cells systematically. An oblique shock is characterized primarily by the upstream Mach number, the shock angle $\beta$, and the deflection angle $\theta$. We study these factors and their combinations and find that the degree of the shock misalignment with the grid determines the optimal number of troubled-cells. On each side of the shock, the optimal set consists of three troubled-cells for aligned shocks, and the troubled-cells identified by tracing the shock and four lines parallel to it, separated by the grid spacing, for nonaligned shocks. We show that the adapted troubled-cell indicator identifies a set of cells that is close to and contains the optimal set of cells for a threshold constant $K = 0.05$, and consequently, produces a solution close to that obtained by limiting everywhere, but with improved convergence.

\end{abstract}

\section{Introduction}
\label{sec:Intro}
In this paper, we show that employing a troubled-cell indicator along with a standard MUSCL based finite volume solver enhances convergence while avoiding non-physical oscillations. We study the effect of the number of troubled-cells showing there is a trade-off between convergence and quality. We employ a novel monotonicity parameter to compare candidate solutions to identify and recommend an optimal number of troubled-cells.

The Finite Volume Method (FVM) is a widely employed numerical technique for solving hyperbolic conservation laws arising in various scientific and engineering disciplines, particularly compressible fluid dynamics \cite{LeVeque2002} due to its robustness and adaptability across diverse applications. The Monotone Upstream-centered Schemes for Conservation Laws (MUSCL) \cite{VANLEER1979} reconstruction is widely used within the finite volume framework to attain high-order accuracy by reconstructing cell interface values from neighboring cell-average data. However, like other high-order numerical schemes, MUSCL face challenges in accurately capturing shocks and discontinuities without introducing non-physical oscillations, which can potentially lead to numerical instability.

Limiting techniques are often incorporated into the MUSCL reconstruction to prevent such oscillations. These limiters control spurious oscillations in the vicinity of discontinuities, thereby preserving monotonicity while retaining high-order accuracy of the numerical scheme in regions of smooth solution \cite{VANLEER1979,VANLEER1974}. In conventional implementations of the MUSCL reconstruction, the limiter function is typically applied to all computational cells. However, this function remains inactive in most of the cells, particularly, in regions of smooth solution and it often identifies regions near smooth extrema as requiring limiting \cite{BISWAS1994}. This approach guarantees monotonicity but it can lead to unnecessary computational effort and may adversely affect the convergence of the solution to a steady state \cite{WAN2022}. In the subsequent sections, we refer to this as the ``\textbf{limiting everywhere approach}''.

Krivodonova et al. \cite{KRIVODONOVA2004} proposed a discontinuity detector for discontinuous Galerkin (DG) methods. By applying limiters exclusively in the computational cells suggested by the discontinuity detector, they demonstrated a significant improvement over the conventional approach of applying limiters in all computational cells. Qiu and Shu \cite{SHU2005} proposed a weighted essentially non-oscillatory (WENO) finite volume methodology as limiters for DG methods. Their approach involves detecting ``troubled-cells" that potentially contain discontinuities and may require limiting. Within these troubled cells, the original high-order DG solution polynomials are replaced with WENO reconstructed polynomials to retain the original high-order accuracy of the DG method. In the subsequent sections, we refer to the method of applying the limiter function only in troubled-cells as the ``\textbf{limiting restricted region approach}''.

We use this limiting restricted region approach for an FVM solver employing MUSCL reconstruction. To detect such regions near discontinuities, where limiting is needed, we adapt the troubled-cell indicator developed by Fu and Shu \cite{FU2017} for DG methods to our FVM framework. To validate the effectiveness of this limiting approach, we employ various two-dimensional test cases solving compressible Euler equations involving discontinuities, particularly oblique shocks. Limiting restricted region approach shows better convergence without compromising the solution accuracy compared to the limiting everywhere approach, when there are enough number of troubled-cells near the shocks. As the number of troubled-cells near the shocks decreases, the solution exhibits overshoots and undershoots though the convergence is much better.This suggests that there exists an optimal number of troubled-cells that yields a solution without any oscillations while ensuring improved convergence.

In the present paper, we seek the optimal number of troubled-cells required in the neighbourhood of the shocks to obtain a solution without much oscillations or closer to the solution of the limiting everywhere approach with enhanced convergence. For this study, we employ a series of two-dimensional test cases with known solution for given flow conditions, involving oblique shocks both aligned and not aligned with grid lines, at various shock angles. We refer to the shocks which are aligned with grid lines as aligned shocks and the shocks which are not aligned with grid lines as non-aligned shocks.

We pre-label the cells in the neighbourhood of the shock as troubled-cells since the exact location of the shock is known. For the aligned shocks, we simply the label the cells immediately before and after the shock location. For the non-aligned shocks, we use a simple line tracing algorithm to label the troubled cells in the neighbourhood of the shock.
Starting with an initial minimal set of troubled cells, we progressively expand this region of troubled-cells. We then evaluate the quality of the solution obtained by limiting only in the troubled-cells of each set.

To quantify the quality of the solution, particularly in the neighbourhood of the shock, we introduce a novel parameter called the monotonicity parameter $\mu$. It is defined as the difference between the total variation of the error in density and the $L_{\infty}$ norm of the error in density. This parameter serves as a quantitative measure of how closely a solution adheres to monotonicity.

The outline of the rest of the paper is as follows. In section \ref{sec:Methodology}, we briefly review the two-dimensional compressible Euler equations governing the inviscid flow and the FVM with MUSCL reconstruction scheme employed for their solution. In section \ref{sec:Validation}, we validate our solver using a standard test case to ensure accuracy of the solver. In section \ref{sec:Motivation}, we introduce the troubled-cell indicator for finite volume methods and present the results of the limiting restricted region approach, which served as a motivation for the study to find the optimal number of troubled-cells. We also introduce a novel monotonicity parameter to quantify the quality of the solution. In section \ref{sec:TC_methods}, we introduce the methods to label the troubled-cells and the various configurations of number of troubled-cells studied in this investigation. In section \ref{sec:Results}, we present the results of various two-dimensional test cases involving oblique shocks for the study of optimal number of troubled-cells in the neighbourhood of the shock. In section \ref{sec:Comparison}, we compare the results obtained using the limiting restricted region approach for the troubled-cells identified by the indicator and the respective optimal configuration. We conclude in section \ref{sec:Conclusion}.

%


\section{Numerical Methodology}
\label{sec:Methodology}
In this section, we provide a concise overview of the governing equations and the base solver employed to solve them. We describe the various two-dimensional test cases employed in this study. We present the methodologies employed to assess the quality of the solutions obtained.

The base solver utilizes a finite volume approach with MUSCL reconstruction and the Koren limiter \cite{HEMKER1988}. No other scheme is used to improve the quality of the solution or enhance convergence to steady-state so as to isolate the effect of the limiting strategy. As a consequence, there are instances where the base solver stalls and convergence to steady-state is not achieved.

\subsection*{Governing Equations and Base Solver:}
The governing equations for compressible inviscid flows can be written in integral form over a control volume with volume $\mathcal{V}$ and surface area $S$ as
\begin{equation}\label{eq:2d_euler}
\frac{d}{dt}\int_{\mathcal{V}}\textbf{Q} \,d\mathcal{V} + \int_{S} \textbf{H}(\textbf{Q}) \cdot \,\hat{\textbf{n}}\, dS= 0.
\end{equation}

In two dimensions, the vector of conservative variables $\textbf{Q} = [\rho, \rho u, \rho v, \rho E]^T$, the convective flux vector $\textbf{H} = (\textbf{F}, \textbf{G})$, $\textbf{F} = [\rho u, \rho u^2 + p, \rho uv, \rho uH]^T$, $\textbf{G} = [\rho v, \rho uv, \rho v^2 + p, \rho vH]^T$, pressure $p = (\gamma - 1)\rho\left[E - \frac{1}{2}(u^2 + v^2)\right]$, enthalpy $H = E + p/\rho$ and $\gamma = 1.4$.

The region of interest in the flow field is discretized and replaced by an array of non-overlapping cells (or volumes). For a given cell, using the finite volume method, equation (\ref{eq:2d_euler}) can be discretized in space to obtain the semi-discretized form
\begin{equation}\label{eq:semidiscrete}
 \Omega\,\frac{d\bar{\textbf{Q}}}{dt} + \textbf{R}(\bar{\textbf{Q}})= 0
\end{equation}
where, $\bar{\textbf{Q}}$ is the cell average of $\textbf{Q}$ and $\Omega$ is the volume of the cell. The residue $\textbf{R}(\bar{\textbf{Q}})$ is given as
\begin{equation}\label{eq:residue}
\textbf{R}(\bar{\textbf{Q}}) = \sum_{f} \tilde{\textbf{H}}_f(\bar{\textbf{Q}}) \cdot \,\hat{\textbf{n}}_f \,s_f
\end{equation}
where, $f$ is an index over the faces of the cell, $\tilde{\textbf{H}}_f$ is the numerical inviscid flux vector at the face, $\hat{\textbf{n}}_f$ is the unit normal vector at the face and $s_f$ is the area of the cell face.

In this investigation, we use the Advection Upstream Splitting Method (AUSM) family scheme AUSM+ \cite{LIOU1996} to evaluate the numerical inviscid flux $\tilde{\textbf{H}}_f$. This requires the state to the left and right of the face where $\tilde{\textbf{H}}_f$ is being evaluated.

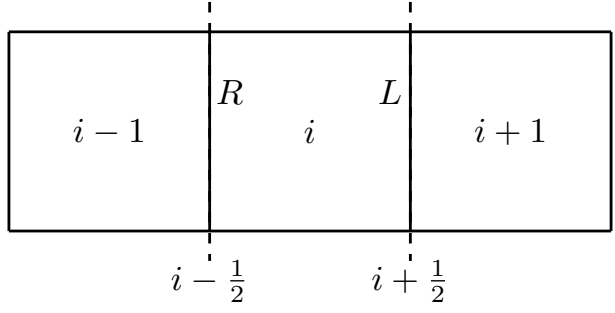
\begin{figure}
\centering
\fontsize{10}{12}\selectfont

\begin{tikzpicture}
 \centering

\draw[thick] (0,0) -- (6,0);   
\draw[thick] (0,-2) -- (6,-2); 

\draw[thick](0,0) -- (0,-2);
\draw[thick](2,0) -- (2,-2);
\draw[thick](4,0) -- (4,-2);
\draw[thick](6,0) -- (6,-2);

\draw (1, -1.0) node[black] {$i-1$};
\draw (3, -1.0) node[black] {$i$};
\draw (5, -1.0) node[black] {$i+1$};

\draw[dashed, thick](2,0.3) -- (2,-2.3);
\draw[dashed, thick](4,0.3) -- (4,-2.3);

\draw (2, -2.5) node[black]{$i-\frac{1}{2}$};
\draw (4, -2.5) node[black]{$i+\frac{1}{2}$};

\draw (2.2, -0.6) node[black] {$R$};
\draw (3.8, -0.6) node[black] {$L$};

%

\end{tikzpicture}
\caption{Stencil of a cell $i$ for the MUSCL reconstruction scheme to obtain the left $(L)$ state at the face $i+1/2$ and the right $(R)$ state at the face $i-1/2$.}
\label{fig:MUSCL_setup}
\end{figure}

We use the higher order Monotone Upwind Schemes for Scalar Conservation Laws (MUSCL) reconstruction strategy \cite{VANLEER1979} to reconstruct the left and the right states at the interface. In the MUSCL reconstruction strategy, the left $(L)$ and the right $(R)$ states at the faces of a cell $i$ can be computed using the cell averages of cells $i-1$, $i$, and $i+1$. These three cells together constitute the stencil $(i-1,\, i,\, i+1)$ for cell $i$ and are shown in Figure \ref{fig:MUSCL_setup}.

The left and right states at the faces of a cell $i$ using the MUSCL reconstruction scheme in terms of primitive variables are determined as follows.
\begin{equation}\label{eq:Muscl_left_wol}
 \textbf{W}^L_{i+\frac{1}{2}} = \textbf{W}_i + \frac{1}{4}\left[(1 - k)\Delta^-\textbf{W}_i + (1 + k)\Delta^+\textbf{W}_i\right]
\end{equation}
\begin{equation}\label{eq:Muscl_right_wol}
 \textbf{W}^R_{i-\frac{1}{2}} = \textbf{W}_{i} - \frac{1}{4}\left[(1 + k)\Delta^-\textbf{W}_{i} + (1 - k)\Delta^+\textbf{W}_{i}\right]
\end{equation}
where, $\textbf{W}$ is the vector of primitive variables, i.e., $\textbf{W} = (\rho, u, v, p)^T$. $\textbf{W}^L_{i+\frac{1}{2}}$ is the desired left state at the face $i+\frac{1}{2}$ and $\textbf{W}^R_{i-\frac{1}{2}}$ is the desired right state at the face $i-\frac{1}{2}$. The forward difference ($\Delta^+$) and the backward difference ($\Delta^-$) operators are defined as $\Delta^+\textbf{W}_i = \textbf{W}_{i+1} - \textbf{W}_i$, $\Delta^-\textbf{W}_i = \textbf{W}_{i} - \textbf{W}_{i-1}$. The parameter $k$ controls the upwind biasing and is taken here as 1/3 to achieve a quadratic reconstruction for smooth solutions.

Traditionally, a limiter is used in the reconstruction to prevent spurious oscillations due to discontinuities in the solution. The left and right states in the MUSCL reconstruction scheme with a slope limiter, $\Phi$, are given by
\begin{equation}\label{eq:Muscl_left_wl}
 \textbf{W}^L_{i+\frac{1}{2}} = \textbf{W}_i + \frac{1}{4}\left[(1 - k) \, \Phi\left(r\right)\, \Delta^-\textbf{W}_i + (1 + k)\, \Phi\left(\frac{1}{r}\right) \,\Delta^+\textbf{W}_i\right]
\end{equation}
\begin{equation}\label{eq:Muscl_right_wl}
 \textbf{W}^R_{i-\frac{1}{2}} = \textbf{W}_{i} - \frac{1}{4}\left[(1 + k) \, \Phi\left(r\right) \, \Delta^-\textbf{W}_{i} + (1 - k) \, \Phi\left(\frac{1}{r}\right) \, \Delta^+\textbf{W}_{i}\right]
\end{equation}
where, $r = \displaystyle \frac{\Delta^+w}{\Delta^-w}$ is the ratio of forward and backward differences of respective primitive variable ($w$). For the slope limiter with symmetry property $\left( \text{i.e.,} \, \Phi(r) = \Phi\left(\frac{1}{r}\right)\right)$, the limited MUSCL reconstruction equations (\ref{eq:Muscl_left_wl}) and (\ref{eq:Muscl_right_wl}) becomes
\begin{equation}\label{eq:Muscl_left_wl_simplified}
 \textbf{W}^L_{i+\frac{1}{2}} = \textbf{W}_i + \frac{1}{2}\, \Psi^L(r)\, \Delta^-\textbf{W}_i
\end{equation}
\begin{equation}\label{eq:Muscl_right_wl_simplified}
 \textbf{W}^R_{i-\frac{1}{2}} = \textbf{W}_{i} - \frac{1}{2}\, \Psi^R(r)\, \Delta^-\textbf{W}_{i}
\end{equation}
with the limiter function $\Psi$ defined as
\begin{equation}
 \Psi^L(r) = \frac{1}{2}\left[(1 - k) + (1 + k) r \right]\Phi(r)
\end{equation}
\begin{equation}
 \Psi^R(r) = \frac{1}{2}\left[(1 + k) + (1 - k) r\right]\Phi(r)
\end{equation}

In the present investigation, we used a slope limiter $\Phi$ given by
\begin{equation}\label{eq:vanAlbada}
 \Phi(r) = \displaystyle \frac{3r}{2r^2 - r + 2}
\end{equation}
with $k = 1/3$. For this slope limiter, the limiter function $\Psi(r)$ corresponds to the limiter of Hemker and Koren \cite{HEMKER1988}.

The semi-discretized form given in equation (\ref{eq:semidiscrete}) can be solved by integrating it in time to obtain the evolution of dependent variables. For this purpose, an implicit matrix-free Lower-Upper Symmetric Gauss-Seidel (LU-SGS) \cite{SHAROV1997} method is utilized with a Courant-Friedrichs-Lewy (CFL) number set to 1 for steady-state problems. For unsteady problems, we use a Total-Variation-Diminishing Runge-Kutta Method (TVD RK3) \cite{SHU1988}, with a CFL number of 0.3.

\subsection*{Test Cases:}
\noindent Test Case 1: Isentropic vortex. This test problem was suggested by \cite{SHU1998} and is primarily used to determine the accuracy of higher-order methods. The problem involves introducing a vortex in the free stream flow, initially centered at $(x_c, y_c)$, using perturbations in velocity and temperature such that the entropy is constant throughout the domain. The perturbations in velocity and temperature are given by
\begin{empheq}{align}
\delta u &= - \frac{\beta}{2\pi} e^{1 - r^2} (y - y_c)  \\
\delta v &= \frac{\beta}{2\pi} e^{1 - r^2} (x - x_c)  \\
\delta T &= -\frac{\gamma - 1}{16\gamma\pi^2} \beta^2 e^{2(1 - r^2)}
\end{empheq}
where, $\beta$ represents the vortex strength, and $r = ((x - x_c)^2 + (y - y_c)^2)^{1/2}$ is the distance from the vortex center to the point $(x, y)$. Using the isentropic relations between density, pressure and temperature, the initial conditions for the flow are given by
\begin{empheq}{align}
\rho &= \left[T_{\infty} - \frac{\gamma - 1}{16\gamma\pi^2} \beta^2 e^{2(1 - r^2)}\right]^{\frac{1}{\gamma - 1}} \label{eq:iv_initial_start}\\
 p &= \rho^{\gamma} \\
u &= u_{\infty} - \frac{\beta}{2\pi} e^{1 - r^2} (y - y_c) \\
v &= v_{\infty} + \frac{\beta}{2\pi} e^{1 - r^2} (x - x_c) \label{eq:iv_initial_end}
\end{empheq}

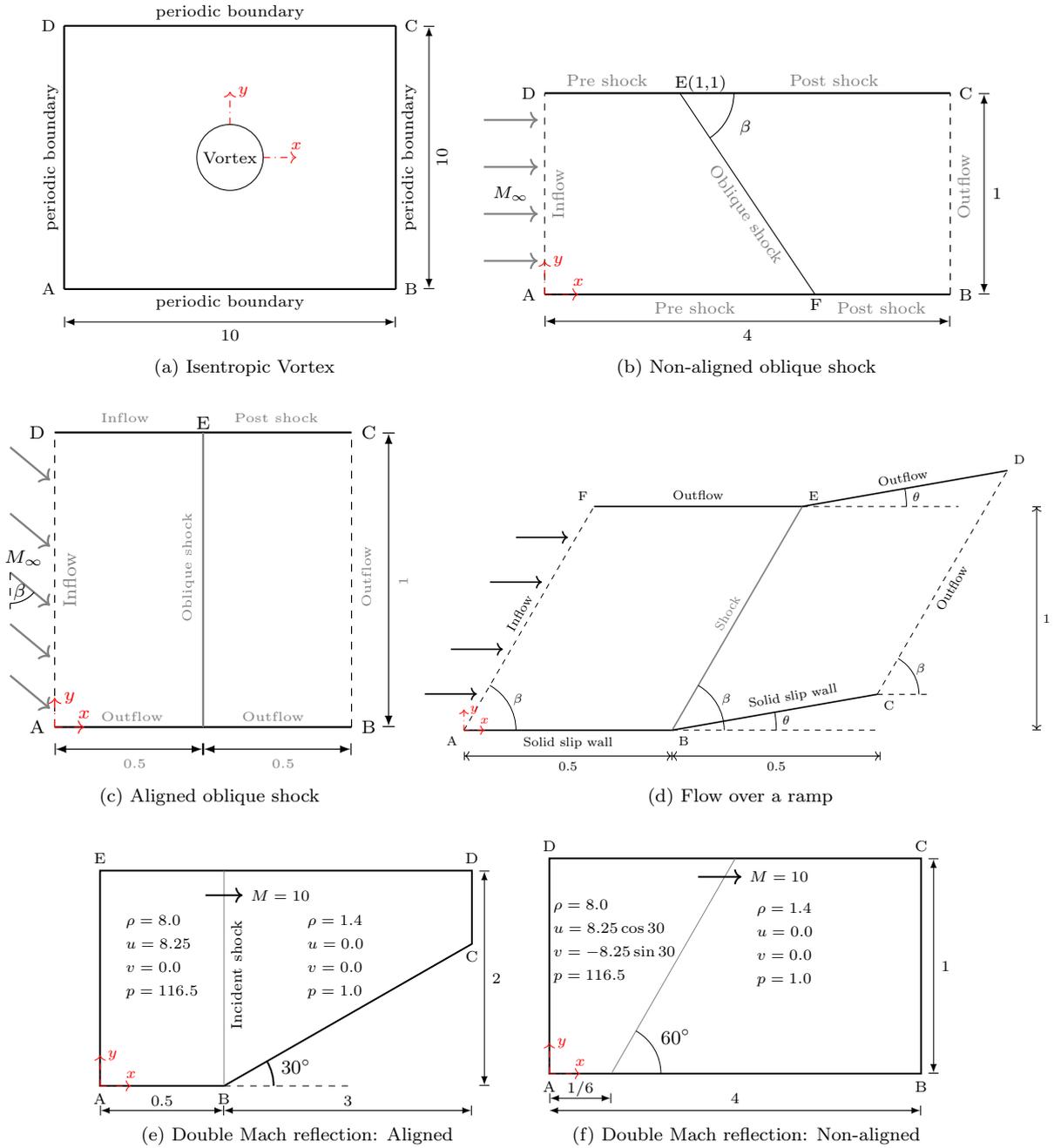
\begin{figure}
\centering
\subfloat[Isentropic Vortex]{

\begin{tikzpicture}[font=\scriptsize]
\centering

\draw[thick] (0,0) node[left]{D} -- (5,0) node[right]{C};
\draw[thick] (0,-4) node[left]{A} -- (5,-4) node[right]{B};
\draw[thick] (5,0) -- (5,-4);
\draw[thick] (0,0) -- (0,-4);
%



\draw (2.5, -2) circle (0.5) node[]{Vortex};

\draw (-0.2, -2)node[rotate=90]{periodic boundary};
\draw (5.2, -2)node[rotate=90]{periodic boundary};
\draw (2.5, -4.2)node[rotate=0]{periodic boundary};
\draw (2.5, 0.2)node[rotate=0]{periodic boundary};




\draw[dashdotted,red,->] (3,-2) -- (3.5,-2)node[above]{$x$};
\draw[dashdotted,red, ->] (2.5,-1.5) -- (2.5,-1.0)node[right]{$y$};

\draw[>=latex, |<->|] (5.5,-4) -- ++(90:4cm);
\draw (5.7, -2)node[rotate=90]{10};

\draw[>=latex, |<->|] (0.0,-4.5) -- ++(0:5cm);
\draw (2.5, -4.7)node[rotate=0]{10};

\end{tikzpicture}\label{fig:IV_setup}}  \hspace{0.3cm}
\subfloat[Non-aligned oblique shock]{

 \begin{tikzpicture}[font=\scriptsize]

 \centering

\draw[thick] (0,0)node[left]{D} -- (6,0)node[right]{C};   
  
\draw[thick] (0,-3) node[left]{A}-- (6,-3)node[right]{B}; 

\draw[dashed](0,-3) -- (0,0); 

\draw[dashed](6,-3) -- (6,0); 

\draw (2,0) -- (4, -3);		

\coordinate (A) at (2cm,0cm);
\coordinate  (C) at (6.0cm,0.0cm);
\coordinate (B) at (4.0cm,-3.0cm);

\tkzMarkAngle[size=0.8cm](B,A,C)
\tkzLabelAngle[pos = 1.1, font=\scriptsize](B,A,C){$\beta$}


\node [] at (-0.5,-1.5) {$M_{\infty}$};

\draw[dashdotted,red,->] (0,-3) -- (0.5,-3)node[above]{$x$};
\draw[dashdotted,red, ->] (0,-3) -- (0,-2.5)node[right]{$y$};

\draw (3.5, 0.2)node[gray,right]{Post shock};
\draw (0.2, 0.2)node[gray,right]{Pre shock};

\draw (4.2, -3.2)node[gray,right]{Post shock};
\draw (1.5, -3.2)node[gray,right]{Pre shock};

\draw (2.4, -1.0)node[gray,right,rotate=-55]{Oblique shock};

\draw[>=latex, |<->|] (6.5,0) -- (6.5,-3);
\draw[>=latex, |<->|] (0,-3.4) -- (6,-3.4);

\draw (6.5, -1.5)node[black, right]{1};
\draw (3.0, -3.4)node[black, below]{4};

\draw (1.8,0.15)node[black, right]{E(1,1)};
\draw (3.8,-3.15)node[black, right]{F};
\draw (6.2,-2.0)node[gray,right, rotate=90]{Outflow};
\draw (0.2,-2.0)node[gray,right, rotate=90]{Inflow};

\draw[thick,gray, ->] (-0.9,-0.4) -- ++(0:0.8cm);
\draw[thick,gray, ->] (-0.9,-1.1) -- ++(0:0.8cm);
\draw[thick,gray, ->] (-0.9,-1.8) -- ++(0:0.8cm);
\draw[thick,gray, ->] (-0.9,-2.5) -- ++(0:0.8cm);

\end{tikzpicture}\label{fig:NonAlignOS_setup}}\\
\subfloat[Aligned oblique shock]{

 \begin{tikzpicture}[font=\scriptsize]

 \centering

\draw[thick] (0,0)node[left]{A} -- ++(0:4cm)node[right]{B};       
\draw[thick] (0,4)node[left, rotate=0]{D} -- ++(0:4cm)node[right, rotate=0]{C};      

\draw(2, 3.92) node[above]{E};

\draw[dashed] (0,0) -- ++(90:4cm);         
\draw[dashed] (4,0) -- ++(90:4cm);      
\draw[thick,gray] (2,0) -- ++(90:4cm);    

\draw (0.5, 4.2)node[gray,right, rotate=0]{\tiny Inflow};
\draw (2.3, 4.2)node[gray,right,rotate=0]{\tiny Post shock};

\draw[>=latex, |<->|] (4.5,0.0) -- ++(90:4cm);
\draw (4.7, 1.8)node[gray,right, rotate=90]{\tiny 1};

\draw[>=latex, |<->|] (0.0,-0.3) -- ++(0:2cm);
\draw (0.8, -0.5)node[gray,right, rotate=0]{\tiny 0.5};

\draw[>=latex, |<->|] (2,-0.3) -- ++(0:2cm);
\draw (2.8, -0.5)node[gray,right, rotate=0]{\tiny 0.5};

\draw[thick,gray, ->] (-0.6,3.8) -- ++(-40:0.7cm);

\draw[thick,gray, ->] (-0.6,2.9) -- ++(-40:0.7cm);
\draw[thick,gray, ->] (-0.6,2.1) -- ++(-40:0.7cm);

\draw (-0.8, 2.3)node[black,right]{$M_{\infty}$};


\draw[dashdotted,red,->] (0,0) -- (0.4,0)node[above]{$x$};
\draw[dashdotted,red, ->] (0,0) -- (0,0.4)node[right]{$y$};

\draw[thick,gray, ->] (-0.6,1.4) -- ++(-40:0.7cm);
\draw[thick,gray, ->] (-0.6,0.7) -- ++(-40:0.7cm);

 \coordinate (A) at (-0.6,2.1);
\coordinate  (C) at (-0.6,1.6);
\coordinate (B) at (0.0,1.6);
\tkzMarkAngle[size=0.42cm](C,A,B)
\tkzLabelAngle[pos = 0.3, font=\scriptsize](C,A,B){$\beta$}
\draw[dashed](-0.6,2.1)--(-0.6,1.6);

\draw (1.8,1.3)node[gray,right, rotate=90]{\tiny Oblique shock};

\draw (4.2,1.5)node[gray,right, rotate=90]{\tiny Outflow};
\draw (0.5,0.15)node[gray,right, rotate=0]{\tiny Outflow};
\draw (2.4,0.15)node[gray,right, rotate=0]{\tiny Outflow};

\draw (0.2,1.5)node[gray,right, rotate=90]{Inflow};

\end{tikzpicture} \label{fig:AlignOS_setup}}
\subfloat[Flow over a ramp]{

 \begin{tikzpicture}[font=\scriptsize]

 \centering

\coordinate (A) at (0, 0);
\coordinate (B) at (4, 0);
\coordinate (F) at ($(A) + (60:5)$); 
\coordinate (E) at ($(B) + (60:5)$);


\draw[thick] (A) -- (B) node[midway, below, sloped] {Solid slip wall};
\draw[thick] (F) -- (E) node[midway, above, sloped] {Outflow};
\draw[dashed] (A) -- (F) node[midway, above, sloped] {Inflow};
\draw[thick, gray] (B) -- (E) node[midway, above, sloped] {Shock};

\node at (A) [below left] {A};
\node at (B) [below right] {B};
\node at (F) [above left] {F};
\node at (E) [above right] {E};

\coordinate (C) at ($(B) + (10:4)$);
\coordinate (D) at ($(E) + (10:4)$);

\draw[thick] (B) -- (C) node[pos = 0.6, above, sloped] {Solid slip wall};
\draw[thick] (E) -- (D) node[midway, above, sloped] {Outflow};
\draw[dashed] (C) -- (D) node[midway, below, sloped] {Outflow};

\node at (C) [below right] {C};
\node at (D) [above right] {D};

\coordinate (G) at (8, 0);
\draw[dashed] (B) -- (G);

\coordinate (H) at ($(E) + (0:3)$);
\draw[dashed] (E) -- (H);

\coordinate (I) at ($(C) + (0:1)$);
\draw[dashed] (C) -- (I);

\draw[|<->|] ([yshift=-0.5cm]A) -- ([yshift=-0.5cm]B) node[midway, below] {0.5};
\draw[<->|] ([yshift=-0.5cm]B) -- ([yshift=-0.5cm]G) node[midway, below] {0.5};
\draw[|<->|] ([xshift=7cm]B) -- ([xshift=7.0cm]B |- E) node[midway, right] {1};

\pic [draw, angle radius=1cm, angle eccentricity=1.2, "$\beta$"] {angle = G--B--E};
\pic [draw, angle radius=1cm, angle eccentricity=1.2, "$\beta$"] {angle = B--A--F};
\pic [draw, angle radius=0.8cm, angle eccentricity=1.2, "$\beta$"] {angle = I--C--D};

\pic [draw, angle radius=2cm, angle eccentricity=1.1, "$\theta$"] {angle = G--B--C};

\pic [draw, angle radius=2cm, angle eccentricity=1.1, "$\theta$"] {angle = H--E--D};

\draw[->, thick] ($($(A)!0.4!(F)$) + (0, 2)$) -- ++(1, 0);
\draw[->, thick] ($($(A)!0.2!(F)$) + (0, 2)$) -- ++(1, 0);
\draw[->, thick] ($($(A)!-0.1!(F)$) + (0, 2)$) -- ++(1, 0);
\draw[->, thick] ($($(A)!-0.3!(F)$) + (0, 2)$) -- ++(1, 0);

\draw[dashdotted,red,->] (0,0) -- (0.4,0)node[above]{$x$};
\draw[dashdotted,red, ->] (0,0) -- (0,0.4)node[right]{$y$};
\end{tikzpicture}\label{fig:AlignedRamp_setup}}\\
\subfloat[Double Mach reflection: Aligned]{

 \begin{tikzpicture}[font=\scriptsize]

 \centering

\coordinate (A) at (0,0); 
\coordinate (B) at (2,0); 
\coordinate (C) at (6,{4*tan(30)}); 
\coordinate (D) at (6,3.5); 
\coordinate (E) at (0,3.5); 

\draw[thick] (A) -- (B) -- (C) -- (D) -- (E) -- cycle;

\coordinate (F) at (6,0);
\draw[dashed] (B) -- (4,0);

\draw[thick] pic["\small $30^\circ$", draw=, angle radius=8mm, angle eccentricity=1.5] {angle=F--B--C};

\draw[gray] (2,0) -- (2,3.5);
\node[right, rotate=90] at (2.2, 0.8) {Incident shock};

\node[align=left] at (1,2.1) {$\rho=8.0$ \\[0.1cm] $u=8.25$ \\[0.1cm] $v=0.0$ \\[0.1cm] $p=116.5$};

\node[align=left] at (3.8,2.1) {$\rho=1.4$ \\[0.1cm] $u=0.0$ \\[0.1cm] $v=0.0$ \\[0.1cm] $p=1.0$};

\draw[dashdotted,red,->] (0,0) -- (0.5,0)node[above]{$x$};
\draw[dashdotted,red, ->] (0,0) -- (0,0.5)node[right]{$y$};

\draw[>=latex, |<->|] (0,-0.4) -- (2.0,-0.4);
\draw (1.0, -0.25)node[black]{0.5};

\draw[>=latex, <->|] (2.0,-0.4) -- (6.0,-0.4);
\draw (4.0, -0.25)node[black]{3};

\draw[>=latex, |<->|] (6.2,0) -- (6.2, 3.5);
\draw (6.2, 1.75)node[black, right]{2};

\draw[thick, ->](1.7, 3.1) -- (2.3, 3.1)node[right]{$M=10$};

\node[below] at (A) {A};
\node[below] at (B) {B};
\node[below] at (C) {C};
\node[above] at (D) {D};
\node[above] at (E) {E};

\end{tikzpicture}\label{fig:DMR_ramp_setup}}  \hspace{0.3cm}
\subfloat[Double Mach reflection: Non-aligned]{

 \begin{tikzpicture}[font=\scriptsize]

 \centering

\coordinate (A) at (0,0); 
\coordinate (B) at (6,0); 
\coordinate (C) at (6,3.5); 
\coordinate (D) at (0,3.5); 

\coordinate (E) at (1,0);
\coordinate (F) at (3,{2.02*tan(60)}); 

\draw[thick] (A) -- (B) -- (C) -- (D) -- cycle;
\draw[gray] (E) -- (F);

\draw[] pic["\small $60^\circ$", draw=, angle radius=8mm, angle eccentricity=1.5] {angle=B--E--F};


\node[align=left] at (1.05,2.15) {$\rho=8.0$ \\[0.1cm] $u=8.25 \cos30$ \\[0.1cm] $v=-8.25 \sin30$ \\[0.1cm] $p=116.5$};

\node[align=left] at (3.8,2.1) {$\rho=1.4$ \\[0.1cm] $u=0.0$ \\[0.1cm] $v=0.0$ \\[0.1cm] $p=1.0$};

\draw[dashdotted,red,->] (0,0) -- (0.5,0)node[above]{$x$};
\draw[dashdotted,red, ->] (0,0) -- (0,0.5)node[right]{$y$};

\draw[>=latex, |<->|] (0,-0.4) -- (1.0,-0.4);
\draw (0.5, -0.25)node[black]{$1/6$};

\draw[>=latex, <->|] (0.0,-0.6) -- (6.0,-0.6);
\draw (3.0, -0.4)node[black]{4};

\draw[>=latex, |<->|] (6.2,0) -- (6.2, 3.5);
\draw (6.2, 1.75)node[black, right]{1};

\draw[thick, ->](2.4, 3.2) -- (3.1, 3.2)node[right]{$M=10$};

\node[below] at (A) {A};
\node[below] at (B) {B};
\node[above] at (C) {C};
\node[above] at (D) {D};

\end{tikzpicture}\label{fig:DMR_setup}}
\caption{Schematic for various test cases, illustrating the varying geometric and flow conditions, and boundary conditions. (not drawn to scale).}
\label{fig:Steady_setup}
\end{figure}

\noindent Test Case 2: Aligned oblique shock. We refer to shocks that coincide with the grid lines as aligned shocks. We use this test case to examine the scenario where the shock aligns with the grid while the flow in both the pre-shock and post-shock regions remains misaligned with the grid. The computational domain for the test case is $[0, 1] \times [0, 1]$. Figure (\ref{fig:AlignOS_setup}) illustrates the computational setup of the test case along with the boundary conditions. \\ \\

Test Case 3: Non-aligned oblique shock. We refer to shocks that do not coincide with the grid lines as non-aligned shocks. This test case is employed to investigate two scenarios: the first, where the flow in the pre-shock region aligns with the grid, while the shock and the flow in the post-region remain misaligned with the grid; and the second, where all three -- the shock and the flow in the pre- and post-shock regions -- are misaligned with the grid. The computational domain for this test case is $[0, 4] \times [0, 1]$. Boundary conditions are applied as shown in Figure (\ref{fig:NonAlignOS_setup}). \\ \\

Test Case 4: Flow over a ramp with an aligned shock. This test case investigates the scenario where all three -- the shock and the flow in the pre- and post-shock regions -- are aligned with the grid. To achieve this, we use a computational domain as illustrated in Figure (\ref{fig:AlignedRamp_setup}) for a flow over a ramp test case. The computational domain comprises two parallelograms, ABEF and BCDE.

The parallelogram ABEF is designed such that the sides AB and EF are parallel to the x-axis, while sides BE and AF are inclined at an angle $\beta$ (shock angle) to the positive x-axis. Similarly, the parallelogram BCDE is constructed with sides BC and DE inclined at an angle $\theta$ (turn angle) to the positive x-axis, and side CD parallel to side BE. Various combinations of shock and turn angles result in different computational domains. This construction aligns the shock but results in a skewed grid. However, the focus of this paper is not on the skewness of the grid and the effects of skewness are not studied in this investigation. \\


\noindent Test Case 5: Double Mach Reflection. This problem is originally suggested by Woodward and Collela \cite{WOODWARD1984}. In this study, we investigate the test case under two configurations: aligned and non-aligned. The aligned setup refers to the case where the shock is initially aligned with the grid, while the non-aligned setup corresponds to the case where the shock is not initially aligned with the grid.

Aligned: The computational domain for this problem is $[0, 3.5] \times [0, 2]$ and is divided into a grid of $840 \times 480$ cells. Initially a right-moving Mach 10 normal shock is positioned at $(x,y) = (0.5,0)$. Along the bottom boundary, we assign the values of initial post-shock flow for the short region AB, and for the rest of the boundary BC, we use the solid slip wall boundary conditions. For the top boundary (i.e. $y = 2$), we set the values to describe the exact motion of the initial Mach 10 normal shock. The computational setup and initial conditions for this problem are given in Figure (\ref{fig:DMR_ramp_setup}). We compute the solution up to $t = 0.2$.

Non-aligned: The computational domain for this problem is $[0, 4] \times [0, 1]$ and is divided into a grid of $960 \times 240$ cells. Initially a right-moving Mach 10 normal shock is positioned at $(x,y) = (1/6,0)$ and makes a \ang{60} angle with positive x-axis. Along the bottom boundary (i.e. $y = 0$), we assign the values of initial post-shock flow for the short region from $x = 0$ to $x = \frac{1}{6}$, and for the rest of the boundary, we use the solid slip wall boundary conditions. For the top boundary (i.e. $y = 1$), we set the values to describe the exact motion of the initial Mach 10 normal shock. The computational setup and initial conditions for this problem are given in Figure (\ref{fig:DMR_setup}). We compute the solution up to $t = 0.2$. \\ \\

To ensure grid independence, we adopt three grid resolutions for test cases 2, 3, and 4, maintaining the same cell sizes across both the aligned and non-aligned shock configurations:
\begin{itemize}
 \item Test cases 2 and 4 (aligned shock):  Coarse grid ($100 \times 100$), Medium grid ($200 \times 200$), and Fine grid ($400 \times 400$).
 \item Test case 3 (non-aligned shock): Coarse grid ($400 \times 100$), Medium grid ($800 \times 200$), and Fine grid ($1600 \times 400$).
\end{itemize}

Our simulations indicate that the results exhibit similar behavior regardless of the grid size employed. Consequently, we have chosen to present the results obtained using the medium grid size for test cases 2, 3, and 4.

\subsection*{Methods of evaluation:}
We evaluate the performance of the solver using the following: \\ \\
Test of convergence to steady-state: To determine if a steady-state solution has been achieved, a residual norm, denoted by $RN$, is calculated at each iteration. This norm is defined as:
\begin{equation}\label{eq:RN}
 RN = \sqrt{\sum_{i=1}^{N} \sum_{j=1}^{4} \textbf{R}_{ij}^2(\textbf{Q}) \, \Omega_i}
\end{equation}
where, $N$ represents the number of cells in the domain, $\Omega_i$ is the volume of the $i$-th cell and $\textbf{R}_{ij}(\textbf{Q})$ is the $j$-th component of residue for that cell and is given in equation (\ref{eq:residue}). The convergence criterion is set as $10^{-14}$. The simulation is stopped when the convergence criterion is met or at the end of 15,000 iterations.\\ \\
Line plots: Density profiles of the numerical solution along a line crossing the shock are visualized to assess the accuracy of the solver in capturing the shock. \\ \\
Error norms: To evaluate and quantify the accuracy of the solution, we calculate the $L_2$ and $L_{\infty}$ norms of the error in density along the line on which density is plotted. These are defined as follows:
\begin{equation}\label{eq:L_2}
 L_2(e) = \sqrt{\frac{1}{N} \sum_{i=1}^{N} |e_i|^2}
\end{equation}
\begin{equation}\label{eq:L_inf}
 L_{\infty}(e) = \max_i |e_i|
\end{equation}
Here, the error, $e$, is defined as
\begin{equation}\label{eq:error}
 e = \rho_{\text{exact}} - \rho_{\text{numerical}}
\end{equation}
where, $\rho_{\text{exact}}$ and $\rho_{\text{numerical}}$ represent the densities of the exact and numerical solutions, respectively. \\ \\
Order of accuracy: The numerical order of accuracy, $n$, is estimated using the following expression:
\begin{equation}\label{eq:order}
 n =  \frac{\ln \left(\frac{\text{norm}_{h_1}}{\text{norm}_{h_2}}\right)}{\ln \left(\frac{h_1}{h_2}\right)}
\end{equation}
where, $\text{norm}_{h_1}$ and  $\text{norm}_{h_2}$ denote the error norms corresponding to the solutions on grids with spacings $h_1$ and $h_2$, respectively. The norm in the equation (\ref{eq:order}) can be either of the two error norms mentioned earlier.


\section{Validation of the solver}
\label{sec:Validation}
We solve the two-dimensional Euler equations (\ref{eq:2d_euler}) for the isentropic vortex test case. The computational domain for the test case is $[-5, -5] \times [5, 5]$, as shown in Figure (\ref{fig:IV_setup}). The domain is discretized using grids of four different resolutions: 100 $\times$ 100, 200 $\times$ 200, 400 $\times$ 400, and 800 $\times$ 800 cells. A vortex of strength $\beta = 5$ is introduced at $(x_c, y_c) = (0, 0)$ into a free stream flow characterized by the following conditions: $\rho_{\infty} = 1.0$, $u_{\infty} = 1.0$, $v_{\infty} = 0.0$, and $p_{\infty} = 1.0$. We initialise the computational domain using the conditions for the flow given in equations (\ref{eq:iv_initial_start}) to (\ref{eq:iv_initial_end}). Periodic boundary conditions are applied in both x and y directions. The solution is computed up to $t = 20$.

\begin{figure}
\centering
\includegraphics[width=0.5\textwidth]{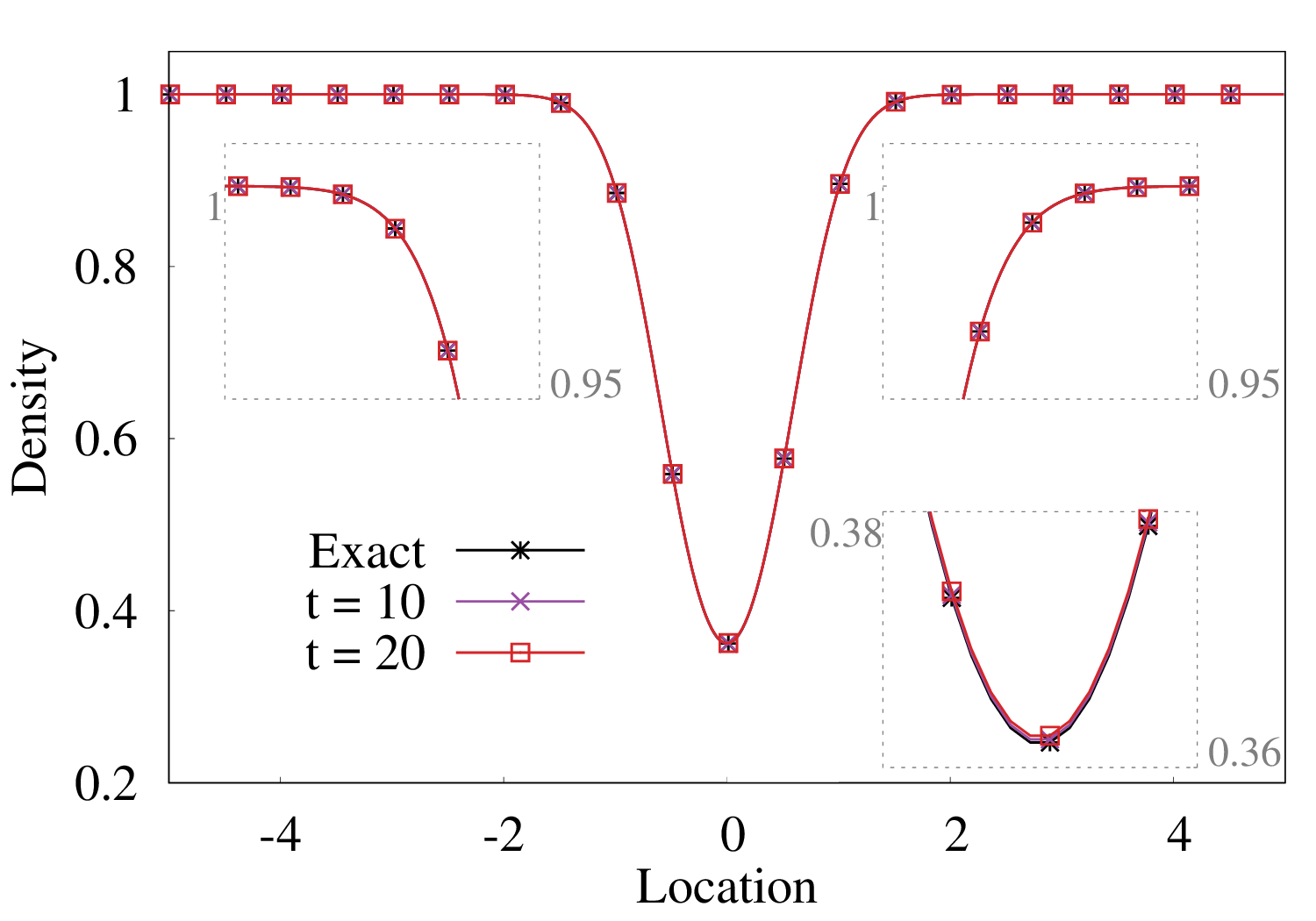}\label{fig:density_iv}
\caption{Isentropic vortex. Density profiles along the line $y = 0$ at time $t = 10$ and $t = 20$.}
\label{fig:Density_profile_iv}
\end{figure}

Figure (\ref{fig:Density_profile_iv}) presents the density profiles along the $y = 0$ line for the solutions at time $t = 10$ and $t = 20$. The results indicate that the computed solution closely matches the exact solution. To estimate the order of accuracy of the scheme, we calculate the $L_2$ and $L_{\infty}$ norms of the error in density for the solution at $t = 20$ over the entire computational domain. As expected for the MUSCL reconstruction with $k = 1/3$, third-order accuracy is achieved. Table (\ref{tab:IV_order}) summarizes the numerical orders of accuracy of the scheme.

\begin{table}
\centering
\caption{Numerical orders of accuracy for the isentropic vortex problem at $t = 20$.}
\begin{tabular}{c c c c c}
    \toprule
  Number of Cells & $L_2$ & order & $L_{\infty}$ & order \\
    \midrule
    100 $\times$ 100 & 2.42E-02 & - & 3.23E-02 & - \\
200 $\times$ 200 & 3.20E-03 & 2.92 & 4.52E-03 & 2.84 \\
400 $\times$ 400 & 4.22E-04 & 2.92 & 5.38E-04 & 3.07 \\
800 $\times$ 800 & 6.47E-05 & 2.71 & 6.34E-05 & 3.09 \\
\bottomrule
\end{tabular}
\label{tab:IV_order}
\end{table}


\section{Motivation}
\label{sec:Motivation}


In this section, we introduce a troubled-cell indicator adapted from DG method. We present sample results using the limiting restricted region approach, which served as motivation for the optimal study outlined in this paper. We use these results to explain the necessity of the novel monotonicity parameter.

\subsection*{Troubled-cell indicator:}
\begin{figure}
\centering

\begin{tikzpicture}
 \centering

\draw[thick] (0,0) -- (6,0);   
\draw[thick] (0,-2) -- (6,-2); 

\draw[thick] (2,-4) -- (2,2);
\draw[thick] (4,-4) -- (4,2);

\draw[thick](0,0) -- (0,-2);
\draw[thick](6,0) -- (6,-2);
\draw[thick](2,-4) -- (4,-4);
\draw[thick](2,2) -- (4,2);

\draw (3, -1) node[black] {$0$};

\draw (1, -1) node[black] {$1$};
\draw (3, -3) node[black] {$2$};
\draw (5, -1) node[black] {$3$};
\draw (3, 1) node[black] {$4$};

\end{tikzpicture}
\caption{Stencil, $S = \{C_0, C_1, C_2, C_3, C_4\}$, used (in the troubled cell indicator) to determine whether the cell $C_0$ is a troubled cell.}
\label{fig:tci_stencil}
\end{figure}
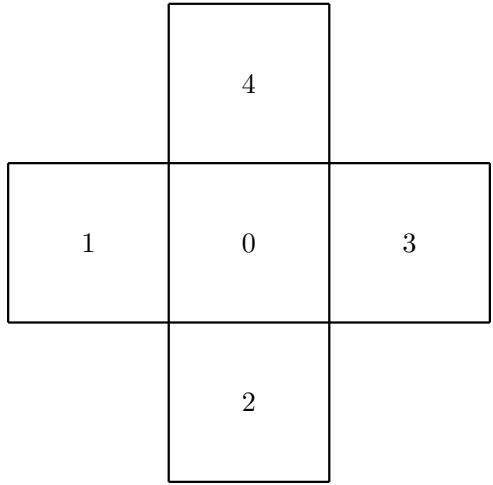
We label the target cell as $C_0$ and denote the stencil to calculate the indicator value as $S = \{C_0, C_1, C_2, C_3, C_4\}$ as shown in Figure (\ref{fig:tci_stencil}). We define the following quantity for the target cell $C_0$,
\begin{equation}\label{eq:TCI}
I_{C_0} = \displaystyle\frac{\sum_{j=1}^4|\bar{\rho}_{C_0} - \bar{\rho}_{C_j}|}{\text{max}_{j\epsilon\{0,1,2,3,4\}}\{\bar{\rho}_{C_j}\}}
\end{equation}
where, $\bar{\rho}_{C_j}$ is the cell average of density of the cell in the stencil. The cell $C_0$ is considered as a troubled cell if $I_{C_0} \geq K$ for a threshold constant $K$.

\begin{figure}
\centering
\subfloat[$K = 0.02$]{\includegraphics[width=0.5\textwidth, height=0.5\textwidth]{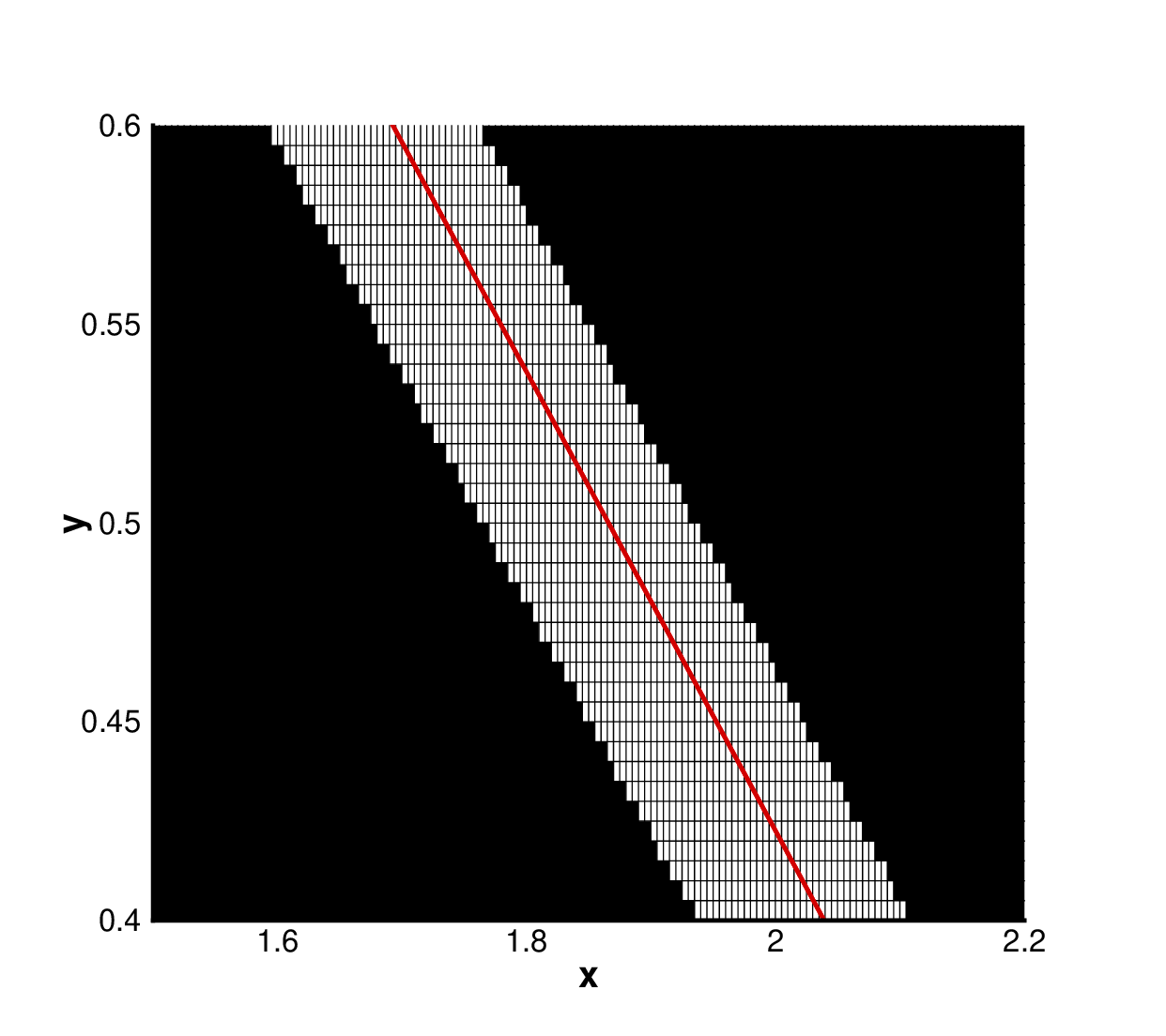}\label{fig:NA30_K002_TC}}
\subfloat[$K = 0.1$]{\includegraphics[width=0.5\textwidth, height=0.5\textwidth]{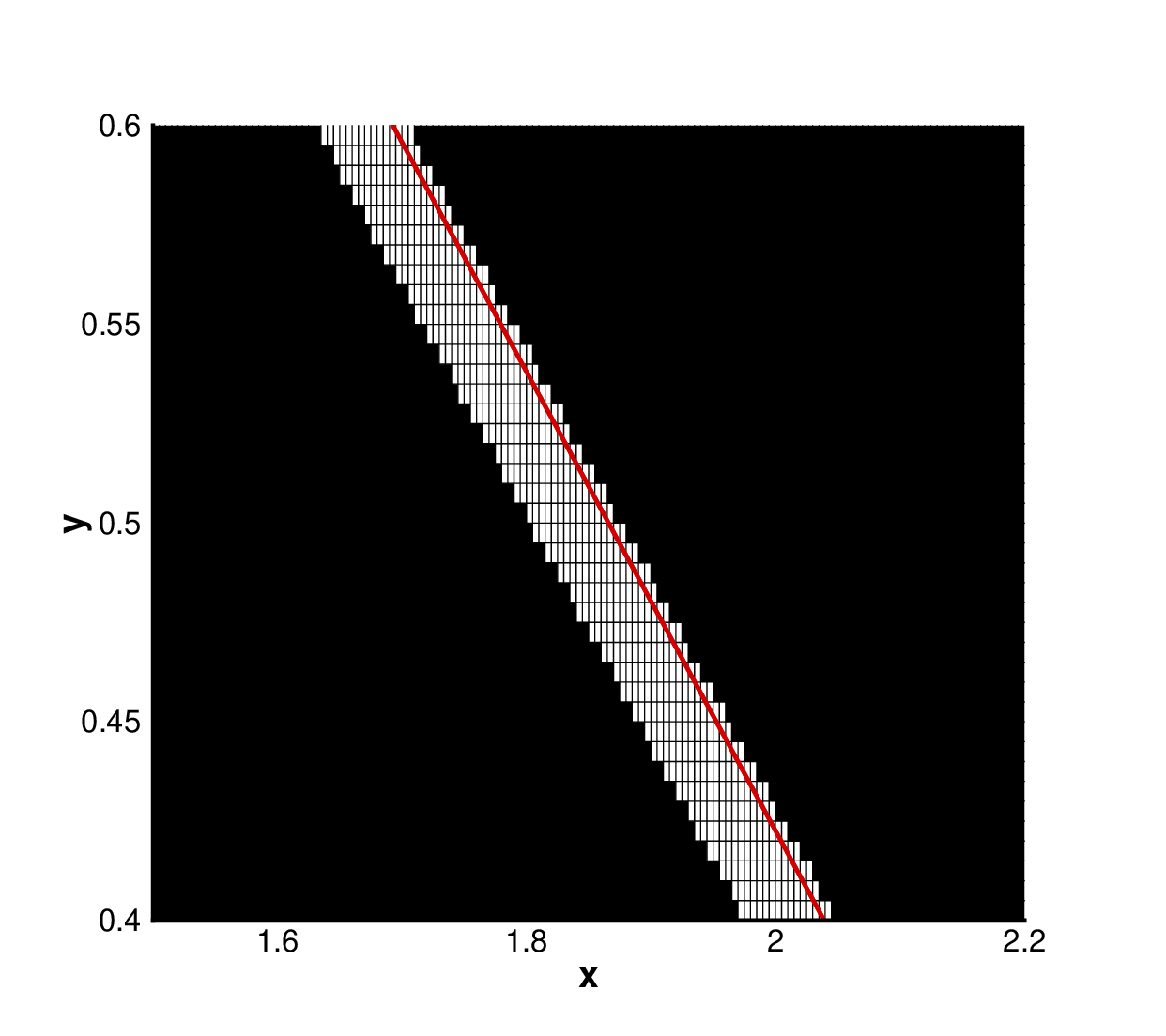}\label{fig:NA30_K01_TC}}
\caption{NonAligned oblique shock. Zoomed-in view of troubled-cells identified by the indicator for two threshold constants. Red line represents the exact shock.}
\label{fig:NonAlignOS_TC}
\end{figure}

We solve the two-dimensional Euler equations (\ref{eq:2d_euler}) for the non-aligned oblique shock test case shown in Figure (\ref{fig:NonAlignOS_setup}). We compute this problem for a Mach number of 3 with shock angle of \ang{30}. We run the simulation using Lax-Friedrich flux scheme to obtain the first-order solution. Troubled-cells are identified using the trouble-cell indicator defined in Equation (\ref{eq:TCI}), applied to the first-order converged solution with threshold constants of $K = 0.02$ and $K = 0.1$. Zoomed-in view of the troubled-cells identified for these threshold constants is illustrated in Figure ({\ref{fig:NonAlignOS_TC}}).


\begin{figure}
\centering
\subfloat[]{\includegraphics[width=0.5\textwidth]{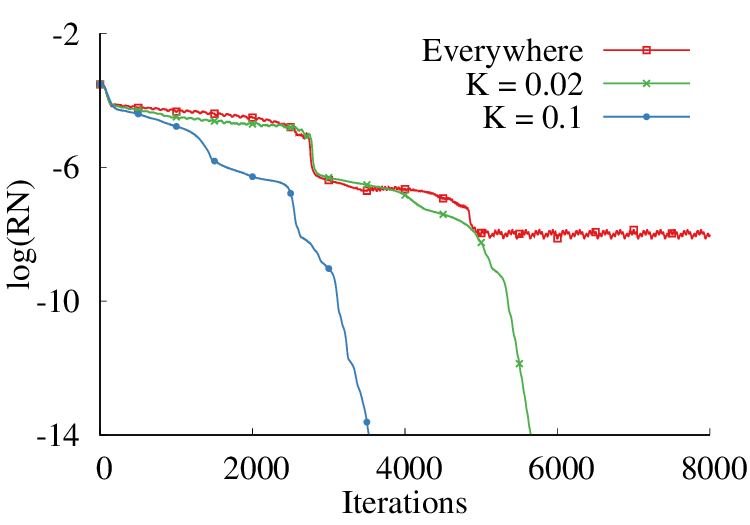}\label{fig:NA30_RN}}
\subfloat[]{\includegraphics[width=0.5\textwidth]{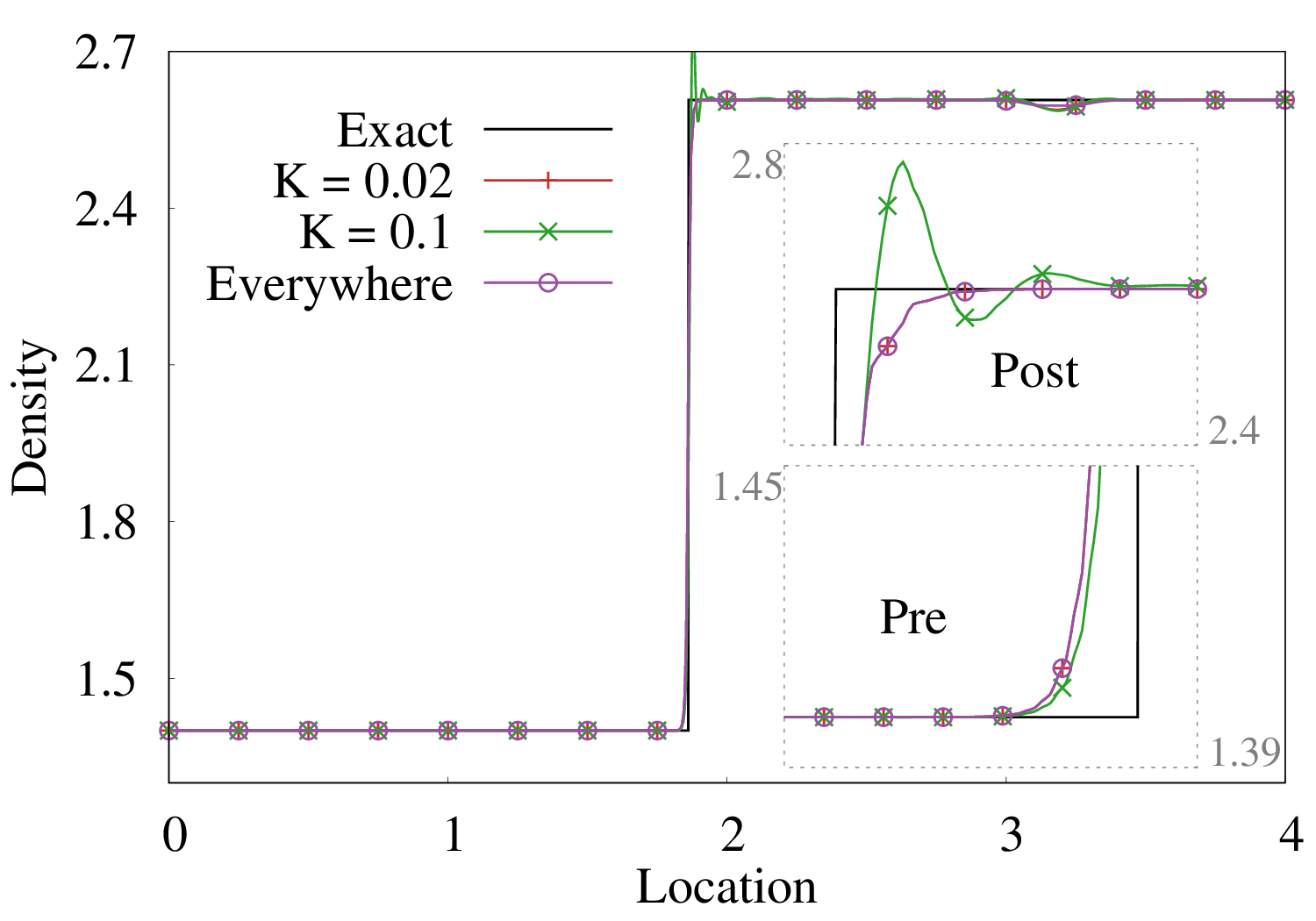}\label{fig:NA30_density}}
\caption{NonAligned oblique shock. (a) The convergence history of the residual norm as a
function of number of iterations. (b) Density profiles along the line $y = 0.5$.}
\label{fig:NonAlignOS_Results}
\end{figure}

We initialise the domain with the first-order converged solution along with the information of troubled-cells. We run the simulation for high-order solution with limiting only in troubled-cells. Figure (\ref{fig:NonAlignOS_Results}) presents the convergence history of the residual norm as a function of the number of iterations and the density profiles along the line $y = 0.5$ for both the solutions of limiting only in troubled-cells for threshold constants $K = 0.02$ and $K = 0.1$ and the solutions of limiting everywhere approach.

The results indicate that the convergence of the limiting restricted region approach is significantly better than that of the limiting everywhere approach. Furthermore, for the limiting restricted region approach, as the threshold constant increases, resulting in fewer identified troubled-cells, the convergence of the residual norm reaches machine zero in fewer iterations compared to a smaller threshold constant.

The density profiles of the limiting restricted region approach, particularly for $K = 0.02$, closely match those of the limiting everywhere approach. However, overshoots are evident in the solution for $K = 0.1$ due to the limited number of troubled cells identified in the post-shock region (Figure (\ref{fig:NA30_K01_TC})).

Table (\ref{tab:Solution_Compare_norms}) presents the $L_2$ and $L_{\infty}$ norms of the error in density for the solutions shown in Figure (\ref{fig:NonAlignOS_Results}). To focus on the vicinity of the shock, the computations of the $L_2$ and $L_{\infty}$ norms are confined to a restricted domain: 20 cells in pre-shock and 20 cells in post-shock from the shock position. The $L_2$ and $L_{\infty}$ norms for these solutions are quite comparable, in fact slightly better for the solution of the limiting restricted region approach for $K = 0.1$. This can lead to the conclusion that the solution to the limiting restricted region for $K = 0.1$ is ``better" than the limiting everywhere despite the presence of significant oscillations (overshoots and undershoots) in the solution.

\begin{table}
\centering
\caption{$L_2$ and $L_{\infty}$ norms of the density error for the solutions of both limiting approaches.}
\begin{tabular}{c c c}
\toprule
\makecell{Limiting} & $L_2$ & $L_{\infty}$ \\
\midrule
\makecell{$K = 0.02$} & 0.114926 & 0.695896\\
\midrule
\makecell{$K = 0.1$} & 0.110430 & 0.652466\\
\midrule
\makecell{Everywhere} & 0.114903 & 0.695483\\
\bottomrule
\end{tabular}%
\label{tab:Solution_Compare_norms}
\end{table}

\subsection*{Monotonicity Parameter:}
If the limiter function is not applied in a sufficient number of troubled-cells in the vicinity of the shock, then there will be overshoots and undershoots in the solution. These oscillations can be captured by the total variation since it is a measure of how oscillatory a solution is. We compute the total variation of the error, $e$ (Equation (\ref{eq:error})), in density. This is defined as follows:
\begin{equation}\label{eq:TV}
TV(e) = \sum_{i=1}^{N-1} |e_{i+1} - e_i|
\end{equation}

\textbf{Monotonic solution property:} If the solution monotonically increases (or decreases) from the exact solution in the pre-shock region, then the total variation of the error and $L_{\infty}$ norm values are equal in that region. Similarly, if the solution monotonically increases (or decreases) to the exact solution in post-shock region, then the total variation of the error and $L_{\infty}$ norm values are equal in that region.

We use this property of monotonic solutions to quantify the quality of the solution. For that, we compute the total variation of the error in density and $L_{\infty}$ norm of the error in density separately for both pre- and post-shock regions. 20 cells in the pre-shock region from the shock position is the domain for the pre-shock region and 20 cells in the post-shock region from the shock position is the domain for the post-shock region.

We define a parameter, $\mu$, called the monotonicity parameter, as
\begin{equation}
 \mu = TV - L_{\infty}
\end{equation}
where $TV$ represents the total variation of the error in density, and $L_{\infty}$ denotes the $L_{\infty}$ norm of the error in density. If the solution is monotonically increasing (or decreasing) from (or to) the exact solution, then $\mu$ is zero. Hence, a solution with a monotonicity parameter value close to zero is considered a better solution.

Table (\ref{tab:Solution_Compare_TV}) presents the $L_{\infty}$ norm and the total variation of the density error, along with the monotonicity parameter values for solutions shown in Figure (\ref{fig:NA30_density}). For the solution obtained using limiting everywhere approach, the $\mu$ value is exactly zero in the pre-shock region, whereas it is not exactly zero but remains close to zero in the post-shock region. This indicates that the solution is monotonically increasing from the exact solution in the pre-shock region and the minimal presence of oscillations in the solution in the post-shock region. Similar results can be observed for limiting restricted region approach for $K = 0.02$. However, for $K = 0.1$, the $\mu$ value in post-shock region is significantly higher than that of limiting everywhere approach. This indicates the presence of substantial oscillations in the solution.

\begin{table}
\centering
\caption{$L_{\infty}$ norm and the total variation of the density error and the monotonicity parameter for both limiting approaches.}
\begin{tabular}{c c c c c}
\toprule
\makecell{Limiting} & \makecell{Region} & $TV$ & $L_{\infty}$ & $\mu$ \\
\midrule
\multirow{2}{*}{$K = 0.02$} & Pre-shock & 0.695896 & 0.695896 & 0\\ \cline{2-5}
& Post-shock & 0.451132 & 0.451019 & 1.13e-04\\
\midrule
\multirow{2}{*}{$K = 0.1$} & Pre-shock & 0.470879 & 0.470877 & 2.00e-06 \\ \cline{2-5}
& Post-shock & 1.116430 & 0.652466 & 4.64e-01\\
\midrule
\multirow{2}{*}{Everywhere} & Pre-shock & 0.695483 & 0.695483 & 0 \\ \cline{2-5}
& Post-shock & 0.451545 & 0.451382 & 1.63e-04\\
\bottomrule
\end{tabular}%
\label{tab:Solution_Compare_TV}
\end{table}

%

Overall, the results demonstrate that the limiting restricted region approach produces a solution which is similar to that of limiting everywhere approach with enhanced convergence, when there are sufficient number of troubled-cells. A reduction in the number of troubled-cells further enhances convergence, achieving a steady state solution in fewer number of iterations. However, the solution exhibits overshoots and undershoots. This suggests the existence of an optimal number of troubled-cells that yields a solution without any oscillations while ensuring improved convergence.

In the subsequent sections, we vary the number of troubled-cells in the neighbourhood of the shock for the limiting restricted region approach. We select the solution whose monotonicity parameter value is close to that of the limiting everywhere approach. For that, we calculate the $L_{\infty}$ norm of the error in density separately for the pre-shock and post-shock regions, each confined to 20 cells from the shock location along the line where the density is plotted. Similarly, we evaluate the total variation of the error in density for the pre-shock and post-shock regions. Instead of presenting these values separately, we add the pre- and post-shock values and present an overall value for the entire domain of 40 cells. We compute the monotonicity parameter using these overall $L_{\infty}$ and total variation values.

\section{Methods to Label Troubled cells}
\label{sec:TC_methods}
To determine the optimal number of troubled-cells, we employ test cases with known shock locations for the given flow conditions. In this paper, we investigate a single oblique shock. We systematically label troubled-cells in the vicinity of the shock. Starting with a minimal set of cells, we progressively expand this region.
\subsection*{Aligned Shocks}
For shocks aligned with the grid, labelling troubled-cells is straightforward. We simply select the cells immediately before and after the known shock location.

Table (\ref{tab:TC_CaseNotation_align}) details the various configurations of troubled cells considered for the aligned shock case. It presents the notation used to identify each configuration, along with the number of troubled cells before and after the shock. For instance, configuration `13' means one cell before the shock and three cells after the shock are labelled as troubled cells and maintained throughout the simulation.

\begin{table}
\centering
\caption{The various configurations of number of troubled cells before and after the shock considered for the aligned shock case, along with the notation used to identify each configuration.}
\begin{tabular}{c r@{\hspace{0\tabcolsep}}l}
    \toprule
 \multirow{2}{*}{\makecell{Configuration \\ Notation}}  & \multicolumn{2}{c}{\multirow{2}{*}{\centering \makecell{Troubled Cells \\ (\textcolor{blue}{Before} $|$ \textcolor{red}{After}) the shock}}} \\
 & & \\
    \midrule
`11'&\bcb{1} & \bcr{1} \\
`12'&\bcb{1} & \bcr{1}\bcr{2} \\
`13'&\bcb{1} & \bcr{1}\bcr{2}\bcr{3} \\
`21'&\bcb{2}\bcb{1} & \bcr{1} \\
`22'&\bcb{2}\bcb{1} & \bcr{1}\bcr{2} \\
`23'&\bcb{2}\bcb{1} & \bcr{1}\bcr{2}\bcr{3}\\
`31'&\bcb{3}\bcb{2}\bcb{1} & \bcr{1}\\
`32'&\bcb{3}\bcb{2}\bcb{1} & \bcr{1}\bcr{2}\\
`33'&\bcb{3}\bcb{2}\bcb{1} & \bcr{1}\bcr{2}\bcr{3}\\
\bottomrule
\end{tabular}
\label{tab:TC_CaseNotation_align}
\end{table}

\subsection*{Non-Aligned Shocks}
Labeling troubled cells in the case of non-aligned shocks is less straightforward than for aligned shocks, as the shock passes through cell interiors. A common approach involves calculating the distance from the cell center to the shock and labeling the cell as troubled if this distance falls within a predefined threshold. However, this method offers limited flexibility in selecting troubled cells, particularly when aiming to include a wider region on either side of the shock.

To address this limitation, we developed a line-tracing algorithm inspired by Bresenham's line algorithm \cite{Bresenham1998}. While Bresenham's algorithm is typically used for drawing lines between integer-valued points in computer graphics applications, we adapt it to our CFD application, where shock endpoints may not always coincide with cell corners. The proposed algorithm, outlined in Algorithm (\ref{alg:LTA}), determines the coordinates on the exact shock line corresponding to a given cell center and labels the cell as troubled if these coordinates fall within the cell's boundaries.


To get a broader region of potentially troubled cells, we trace cells along the  lines parallel to the exact shock line in both pre-shock and post-shock regions and cells along the actual shock itself. By varying the number of parallel lines to be traced on each side of the shock, we can customize the total number of troubled cells captured in the pre-shock and post-shock regions.

\begin{algorithm}
\caption{Line Tracing Algorithm}\label{alg:LTA}
\begin{algorithmic}
\Require End points of the exact shock in the computational domain \\ \Comment{points $(x_1, y_1)$, $(x_2, y_2)$}.
\State \textbf{Get:} Equation of a line representing the exact shock
\State \hspace{1.2cm}$y = y_1 + m\,(x - x_1)$ \Comment{Slope $m = \displaystyle\frac{y_2 - y_1}{x_2 - x_1}$}
\For {each cell in the domain}
    \State \textbf{Get:} $(x, y)$ i.e., coordinates of the cell corner points
    \State \textbf{Find:} $x_{\text{min}}$, $x_{\text{max}}$, $y_{\text{min}}$ and $y_{\text{max}}$ of the cell corner points
    \State \textbf{Get:} ($x_{\text{c}}, y_{\text{c}}$) i.e., coordinates of the cell center
    \State \textbf{Find:} $y$ coordinate on the shock line for given $x_{\text{c}}$
    \State \hspace{1.1cm} \Comment{$y_{\text{shock}} \gets y_1 + m\,(x_{\text{c}} - x_1)$}
    \State \textbf{Find:} $x$ coordinate on the shock line for given $y_{\text{c}}$
    \State \hspace{1.1cm} \Comment{$x_{\text{shock}} \gets x_1 + (y_{\text{c}} - y_1)/m$}
    \If{($x_{\text{min}} \leq x_{\text{shock}} \leq x_{\text{max}}$) or ($y_{\text{min}} \leq y_{\text{shock}} \leq y_{\text{max}}$)}
        \State Cell is troubled
    \Else
        \State Cell is not troubled
    \EndIf
\EndFor
\end{algorithmic}
\end{algorithm}

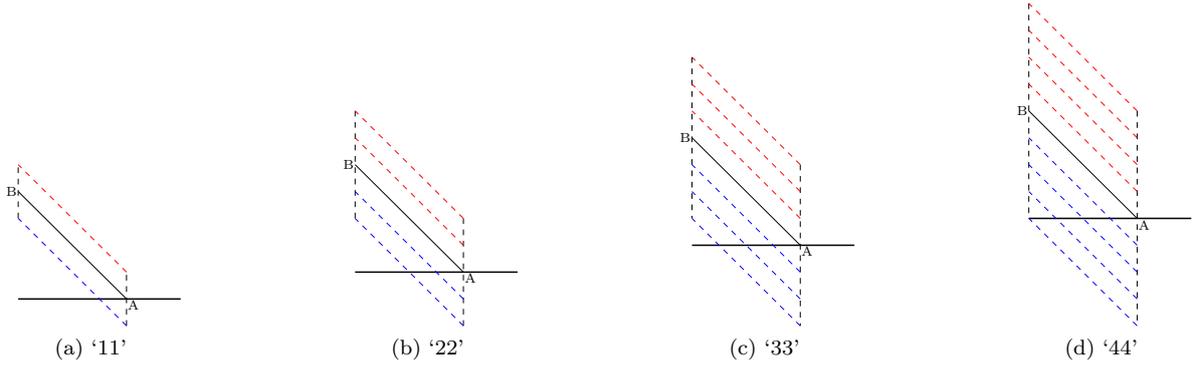
\begin{figure}
\centering
\subfloat[`11']{

 \begin{tikzpicture}[font=\scriptsize]

 \centering
  
\draw[thick] (0,-2) -- (3,-2); 

\draw (0, 0) -- (2, -2);     

\draw[dashed](0,-0.5) -- (0,0.5);
\draw[dashed](2,-2.5) -- (2,-1.5);

\draw [blue, dashed] (2,-2.5) -- (0, -0.5);
\draw [red, dashed] (2,-1.5) -- (0, 0.5);

\draw (1.9, -2.12)node[right]{A};
\draw (0.1, 0)node[left]{B};

\end{tikzpicture}
}
\hfill
\subfloat[`22']{

 \begin{tikzpicture}[font=\scriptsize]

 \centering

\draw[thick] (0,-2) -- (3,-2); 

\draw (0, 0) -- (2, -2);     

\draw[dashed](2,-3) -- (2,-1);
\draw[dashed](0,-1) -- (0,1);

\draw [red, dashed] (2,-1) -- (0, 1);
\draw [red, dashed] (2,-1.5) -- (0, 0.5);
\draw [blue, dashed] (2,-2.5) -- (0, -0.5);
\draw [blue, dashed] (2,-3) -- (0, -1.0);

\draw (1.9, -2.12)node[right]{A};
\draw (0.1, 0)node[left]{B};

\end{tikzpicture}
}
\hfill
\subfloat[`33']{

 \begin{tikzpicture}[font=\scriptsize]

 \centering

\draw[thick] (0,-2) -- (3,-2); 

\draw (0, 0) -- (2, -2);     

\draw[dashed](2,-0.5) -- (2,-3.5);
\draw[dashed](0,1.5) -- (0,-1.5);

\draw [red, dashed] (2,-0.5) -- (0, 1.5);
\draw [red, dashed] (2,-1) -- (0, 1);
\draw [red, dashed] (2,-1.5) -- (0, 0.5);
\draw [blue, dashed] (2,-2.5) -- (0, -0.5);
\draw [blue, dashed] (2,-3) -- (0, -1.0);
\draw [blue, dashed] (2,-3.5) -- (0, -1.5);

\draw (1.9, -2.12)node[right]{A};
\draw (0.1, 0)node[left]{B};

\end{tikzpicture}
}
\hfill
\subfloat[`44']{

 \begin{tikzpicture}[font=\scriptsize]

 \centering

\draw[thick] (0,-2) -- (3,-2); 

\draw (0, 0) -- (2, -2);     

\draw[dashed](2,0) -- (2,-4);
\draw[dashed](0,2) -- (0,-2);

\draw [red, dashed] (2,0) -- (0, 2);
\draw [red, dashed] (2,-0.5) -- (0, 1.5);
\draw [red, dashed] (2,-1) -- (0, 1);
\draw [red, dashed] (2,-1.5) -- (0, 0.5);
\draw [blue, dashed] (2,-2.5) -- (0, -0.5);
\draw [blue, dashed] (2,-3) -- (0, -1.0);
\draw [blue, dashed] (2,-3.5) -- (0, -1.5);
\draw [blue, dashed] (2,-4) -- (0, -2);

\draw (1.9, -2.12)node[right]{A};
\draw (0.1, 0)node[left]{B};

\end{tikzpicture}
}
\caption{A few configurations of parallel lines before and after the shock considered for the non-aligned shock case, along with the notation used to identify each configuration. The other configurations considered in this study are `12', `13', `14', `22', `23', `24', `32', `33', `34', `42', and `43'. Here the line AB represents actual shock, blue dashed lines represent lines parallel to shock in pre-shock region and red dashed lines represent lines parallel to shock in post-shock region.}
\label{fig:Method3_cases}
\end{figure}

Figure (\ref{fig:Method3_cases}) illustrates various configurations of parallel lines considered in this study. The parallel lines are drawn with an offset distance $h$ (cell size) in the vertical direction as shown in Figure (\ref{fig:Method3_cases}). The notation used in the figure helps identify each configuration. For example, configuration `11' signifies one parallel line in the pre-shock and one parallel line in the post-shock regions to be traced, along with the actual shock line itself.

\begin{figure}
\centering
\subfloat[Shock line]{\includegraphics[width=0.48\textwidth]{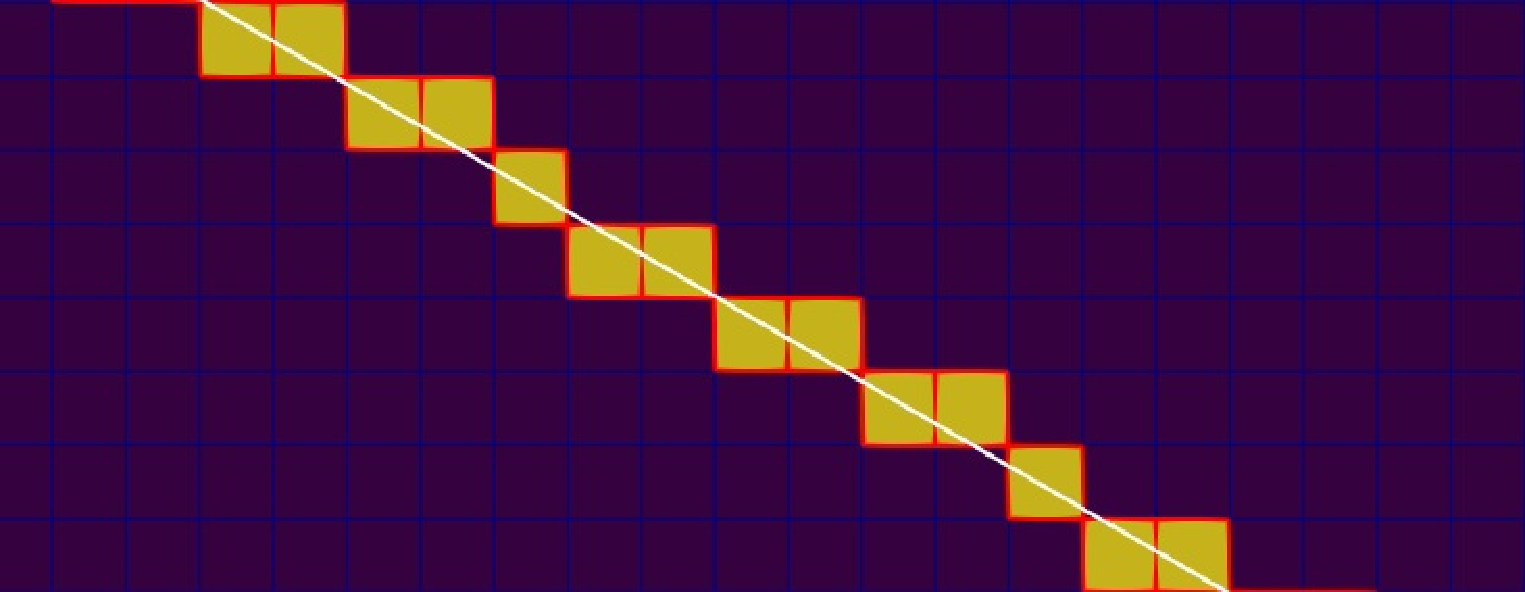}\label{fig:shock_line}}
\hspace{1em}
\subfloat[Parallel line in the pre-shock region]{\includegraphics[width=0.48\textwidth]{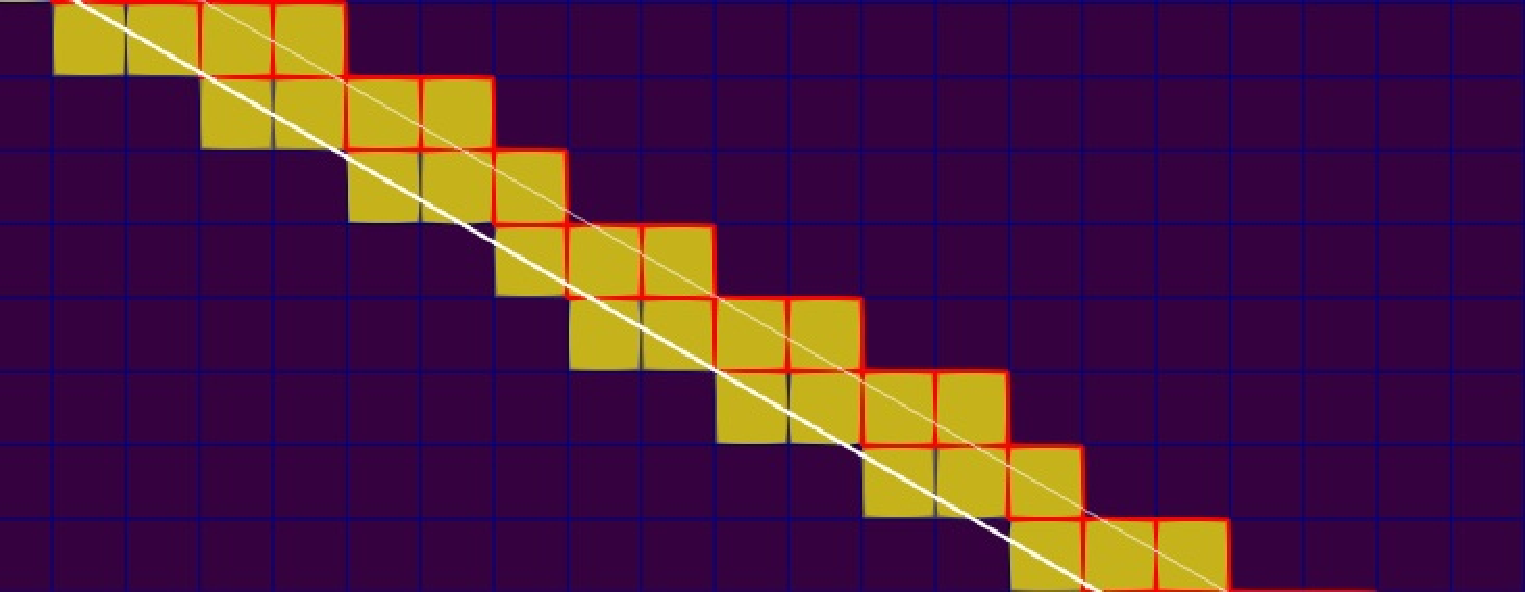}\label{fig:pre_line}}\\
\subfloat[Parallel line in the post-shock region]{\includegraphics[width=0.48\textwidth]{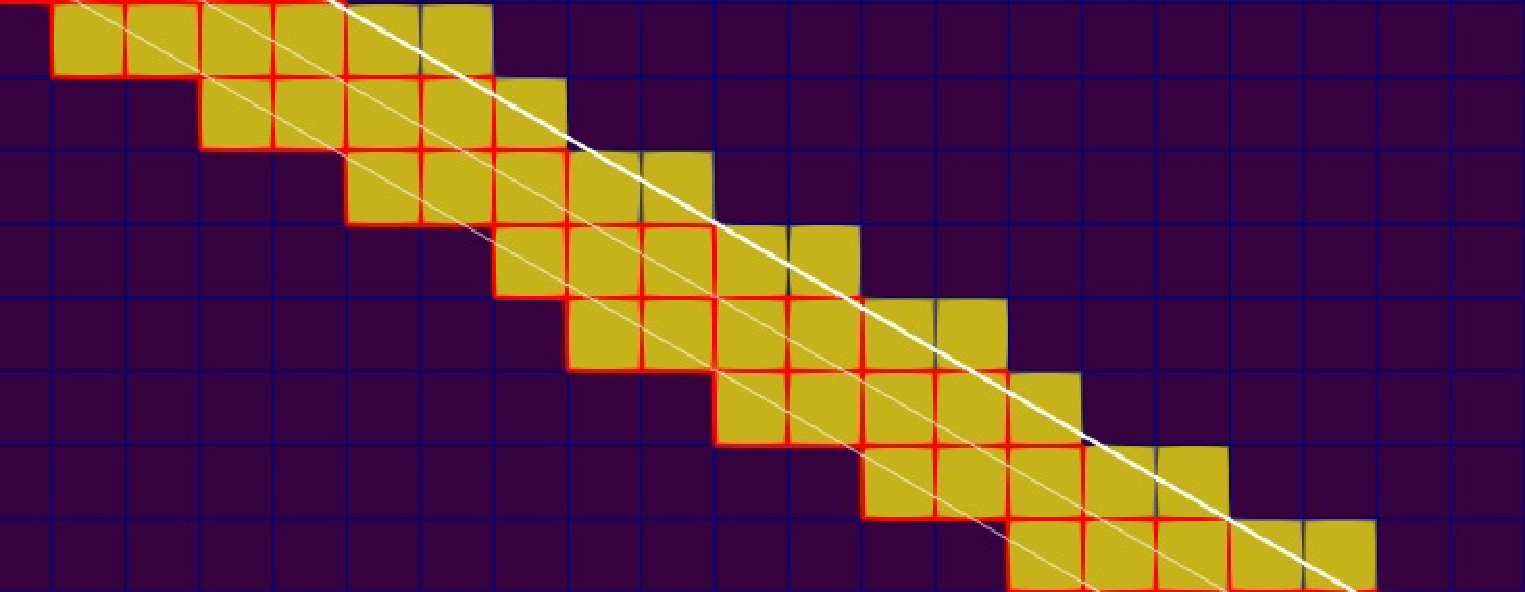}\label{fig:post_line}}
\caption{Zoomed-in view of troubled-cells labelled using line tracing algorithm for the configuration `11'. The highlighted white lines in (a), (b), and (c) are the shock line, parallel line to the shock line in the pre- and post-shock region, respectively. (a): Troubled-cells labelled by tracing the shock line. (b) and (c): cells filled with yellow color except the ones with red color border are the troubled-cells labelled by tracing parallel line to the shock line in the pre- and post-shock region, respectively. Overall, all the cells filled  with yellow color in (c) are the troubled-cells for the configuration `11'.}
\label{fig:NonAlign_TC_11}
\end{figure}

Figure (\ref{fig:NonAlign_TC_11}) presents a zoomed-in view of troubled-cells labelled using the line tracing algorithm provided in Algorithm (\ref{alg:LTA}) for the configuration `11'. Figure (\ref{fig:shock_line}) presents the troubled-cells by tracing just the shock line. Figure (\ref{fig:pre_line}) presents the troubled-cells by tracing both the shock line and parallel line in the pre-shock region. Figure (\ref{fig:post_line}) presents the troubled-cells by tracing all the three lines: the shock line and two parallel lines, one in the pre-shock region and one in the post-shock region. For configuration `23', a total of six lines need to be traced: two parallel lines in the pre-shock region, three parallel lines in the post-shock region, and the actual shock line. Similarly, for other configurations, we traced the parallel lines and the shock line to label the troubled-cells around the shock.

The same configuration notation is used for both aligned and non-aligned shocks for simplicity. For aligned shocks, the notation `XY' signifies that there are `X' number of troubled-cells in the pre-shock region and `Y' number of troubled-cells in the post-shock region. In the case of non-aligned shocks, the notation `XY' implies tracing `X' number of parallel lines in the pre-shock region and `Y' number of parallel lines in the post-shock region in addition to the actual shock line.

\section{Optimality study}
\label{sec:Results}
We solve the two-dimensional Euler equations (\ref{eq:2d_euler}) for both aligned (Figure \ref{fig:AlignOS_setup}) and nonaligned shocks (Figure \ref{fig:NonAlignOS_setup}). The computational domain is initialised with the exact solution corresponding to a upstream Mach number of $3$ and the shock angle of interest. The base solver is employed to solve these equations with the limiting restricted region approach. We investigate various configurations of number of troubled-cells in the pre- and post-shock regions, as detailed in Table (\ref{tab:TC_CaseNotation_align}) for the aligned shock case, and various configurations mentioned in Figure (\ref{fig:Method3_cases}) for the nonaligned shock case.

\begin{figure}
\centering
\subfloat[\ang{30}]{\includegraphics[width=0.5\linewidth]{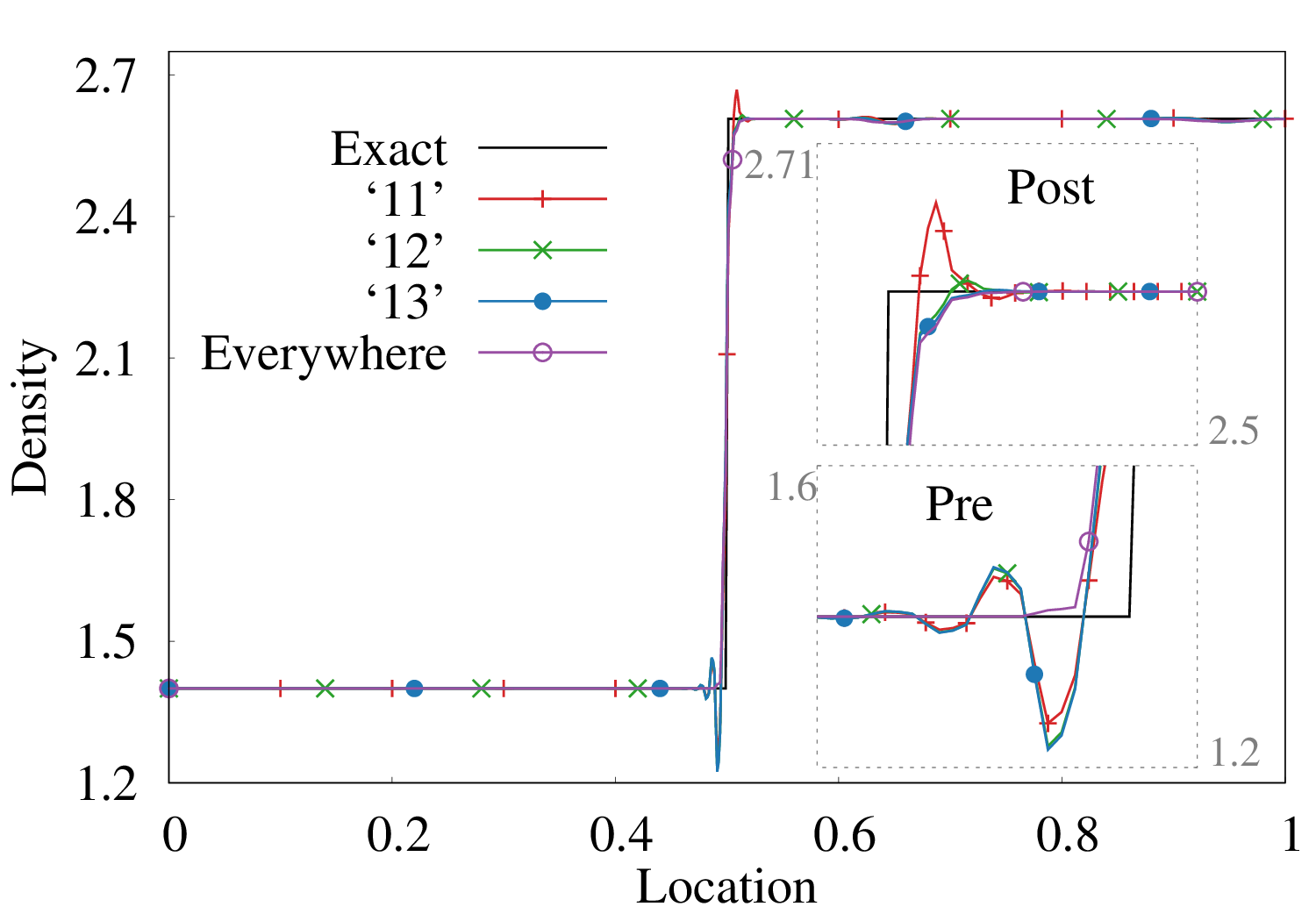}\label{fig:A30_density_pre1}}
\subfloat[\ang{30}]{\includegraphics[width=0.5\linewidth]{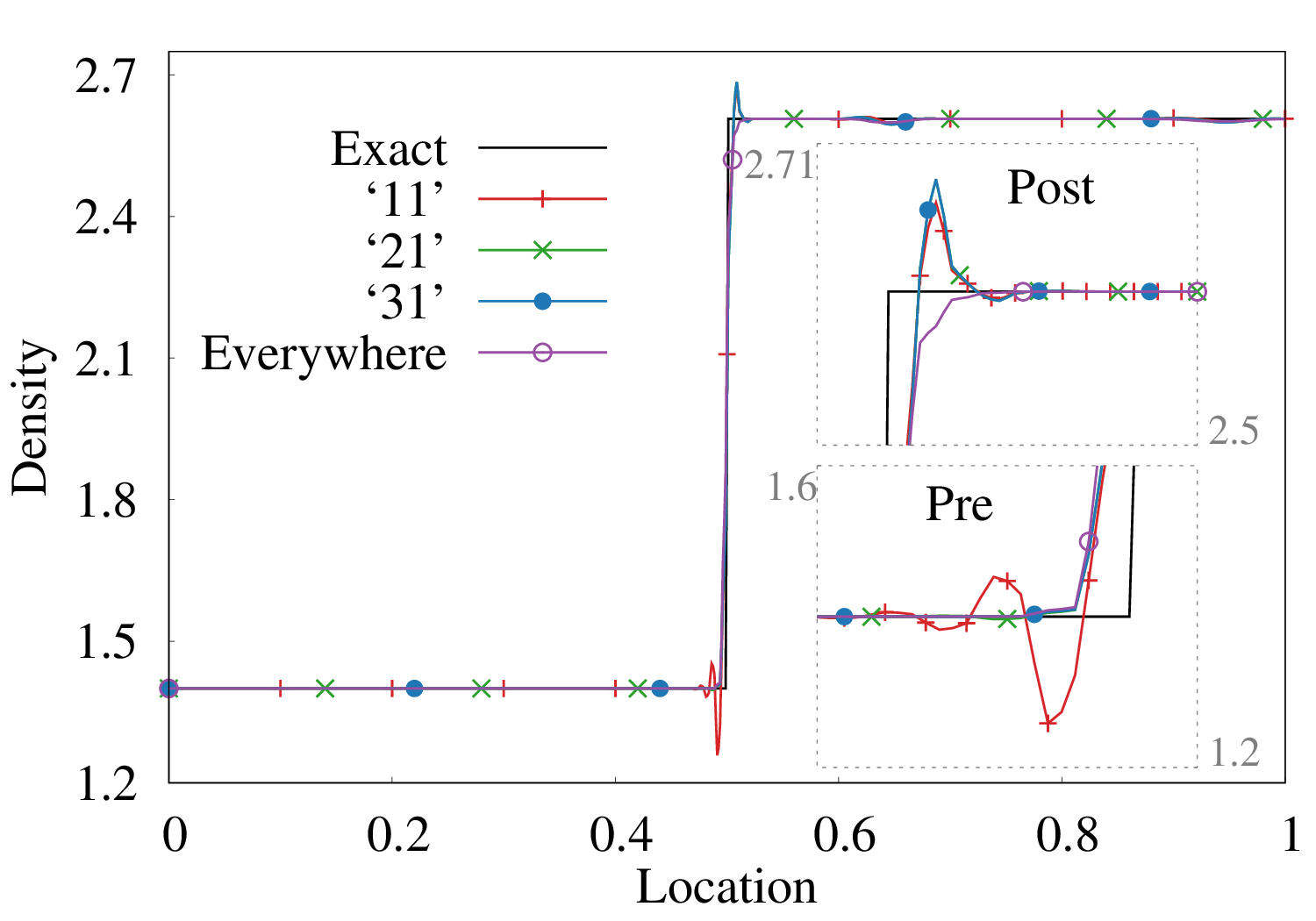}\label{fig:A30_density_post1}}\\
\subfloat[\ang{40}]{\includegraphics[width=0.5\linewidth]{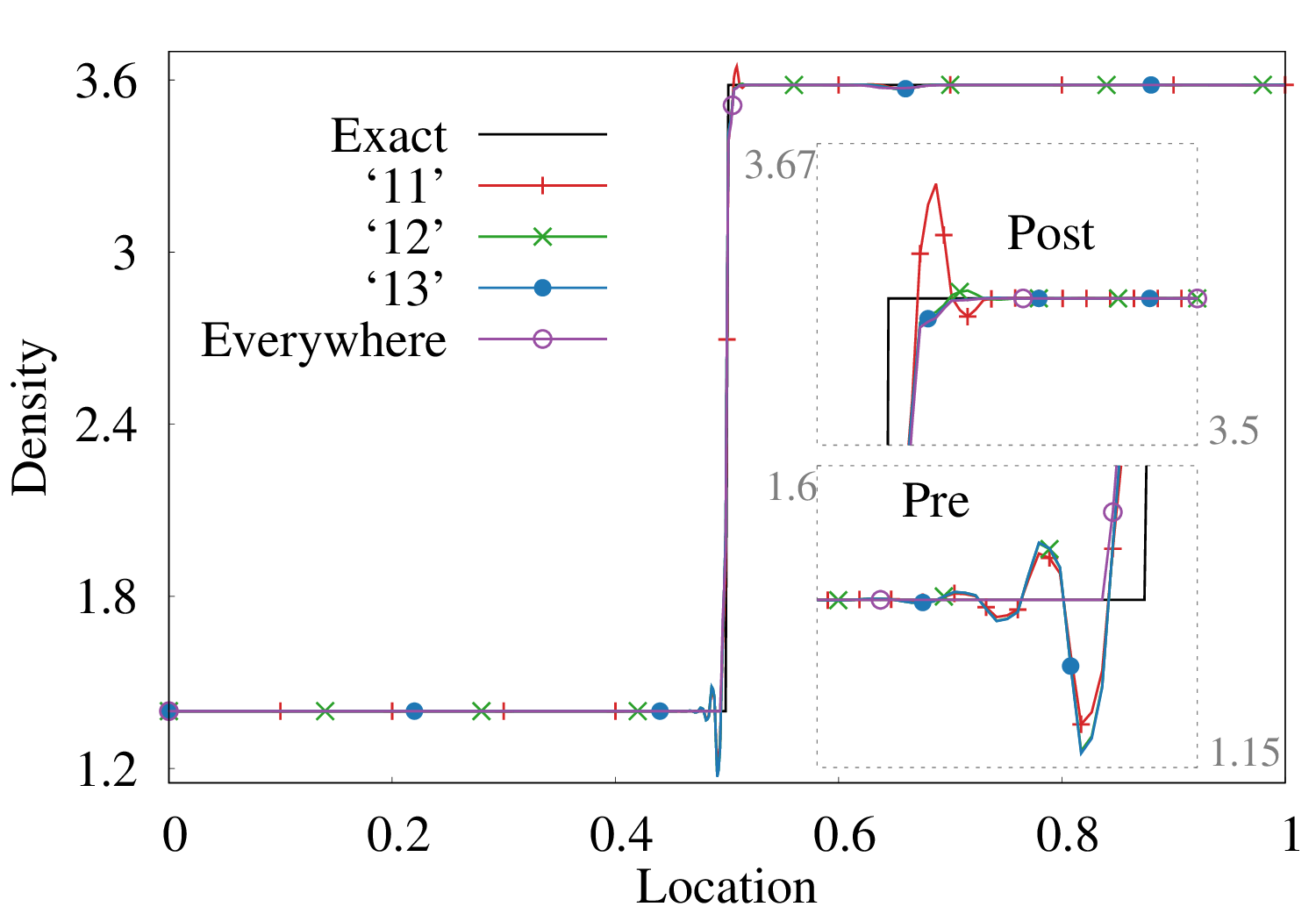}\label{fig:A40_density_pre1}}
\subfloat[\ang{40}]{\includegraphics[width=0.5\linewidth]{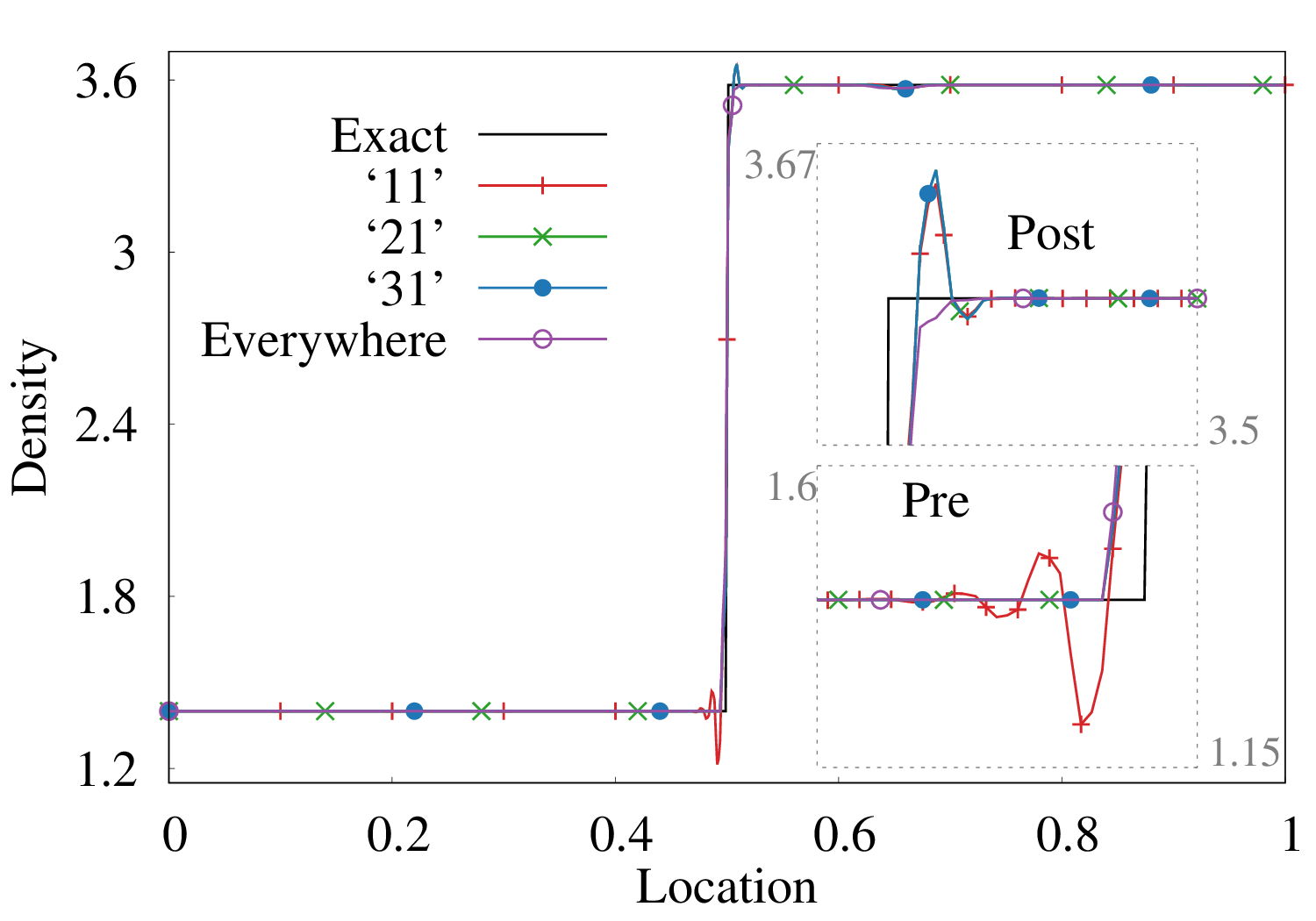}\label{fig:A40_density_post1}}\\
\subfloat[\ang{50}]{\includegraphics[width=0.5\linewidth]{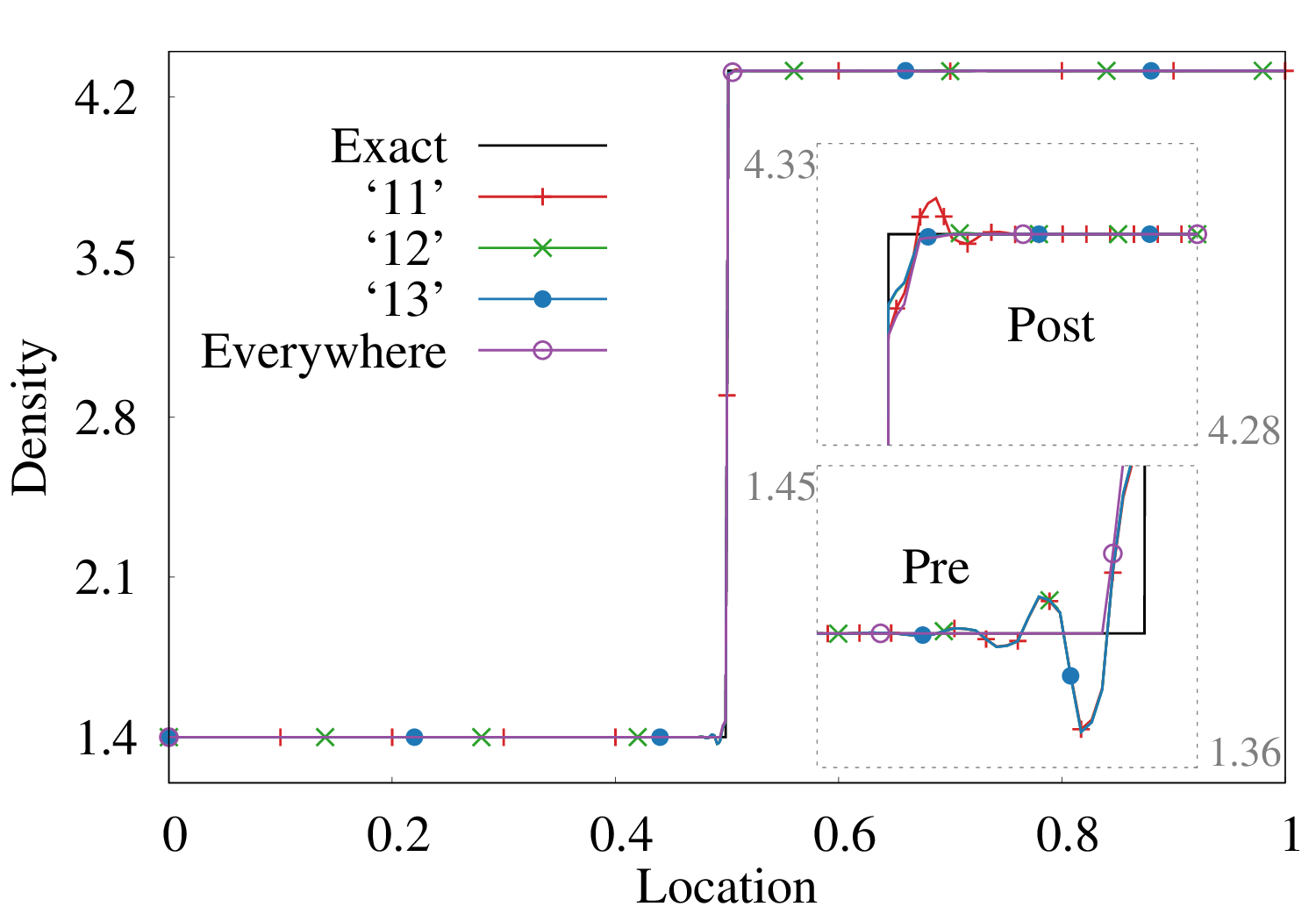}\label{fig:A50_density_pre1}}
\subfloat[\ang{50}]{\includegraphics[width=0.5\linewidth]{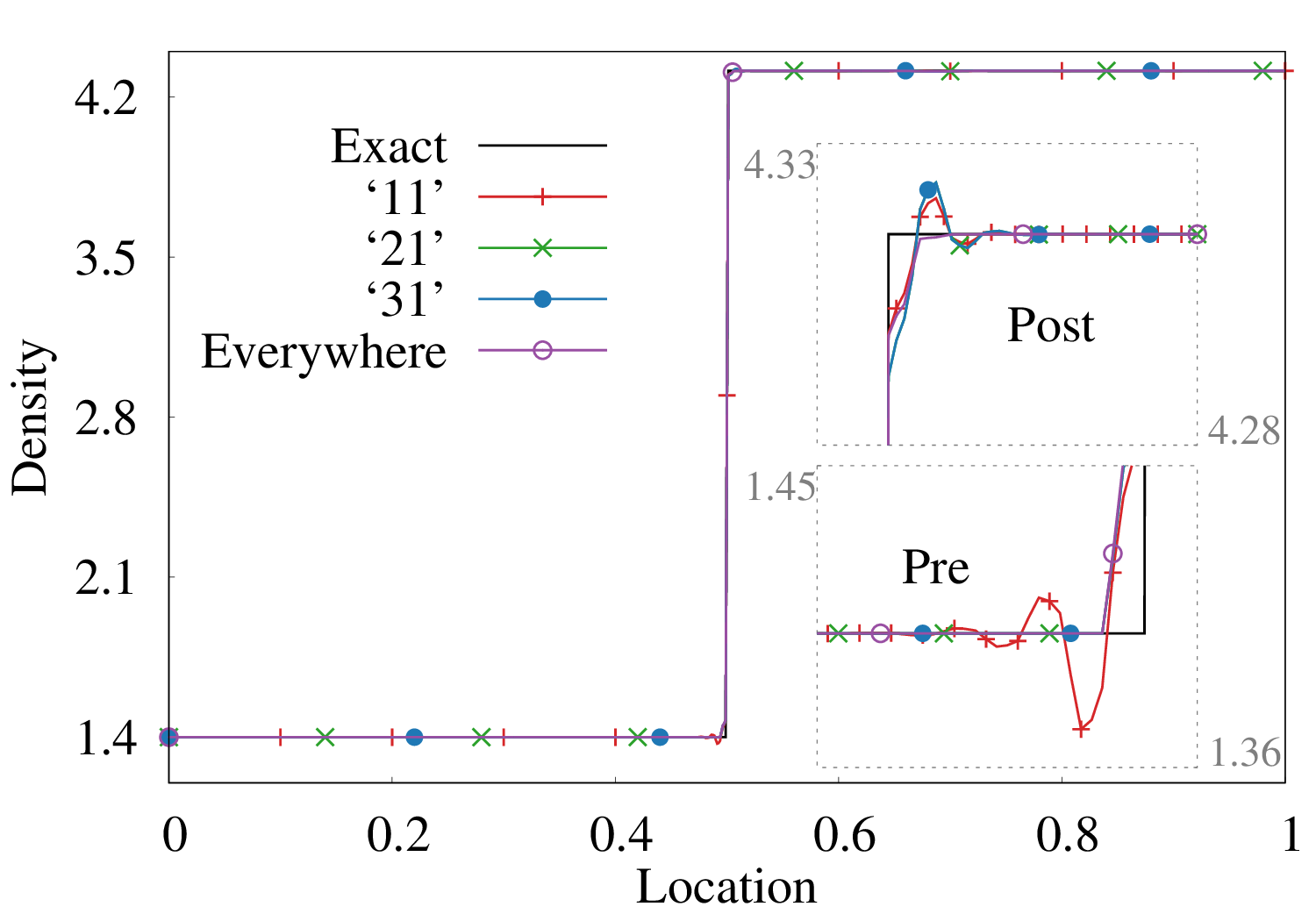}\label{fig:A50_density_post1}}
\caption{Aligned oblique shock with upstream Mach number 3. Density profiles along the line $y = 0.5$ for three different shock angles on the medium grid. Top row (a, b): \ang{30}. Middle row (c, d): \ang{40}. Bottom row (e, f): \ang{50}. Left column (a, c, e): configurations with a single troubled-cell in the pre-shock. Right column (b, d, f): configurations with a single troubled-cell in the post-shock.}
\label{fig:AOS_density_PreAndPost1}
\end{figure}

Density profiles of the cells just above the line $y = 0.5$ for three different oblique shock angles using the medium grid are presented in Figure (\ref{fig:AOS_density_PreAndPost1}) for the aligned shock case. These results represent configurations with a single troubled-cell either in pre- or post-shock region. When the limiter is applied only to a single troubled-cell in the pre-shock region (`1x'), the solution exhibits undershoots in the pre-shock region regardless of the number troubled-cells in the post-shock region. Similarly, when the limiter is applied only to a  single troubled-cell in the post-shock region (`x1'), the solution exhibits overshoots in the post-shock region regardless of the number troubled-cells in the pre-shock region. This pattern holds true for all the three shock angles, though the magnitude of the overshoots or undershoots varies.

\begin{figure}
\centering
\subfloat[\ang{30}]{\includegraphics[width=0.5\linewidth]{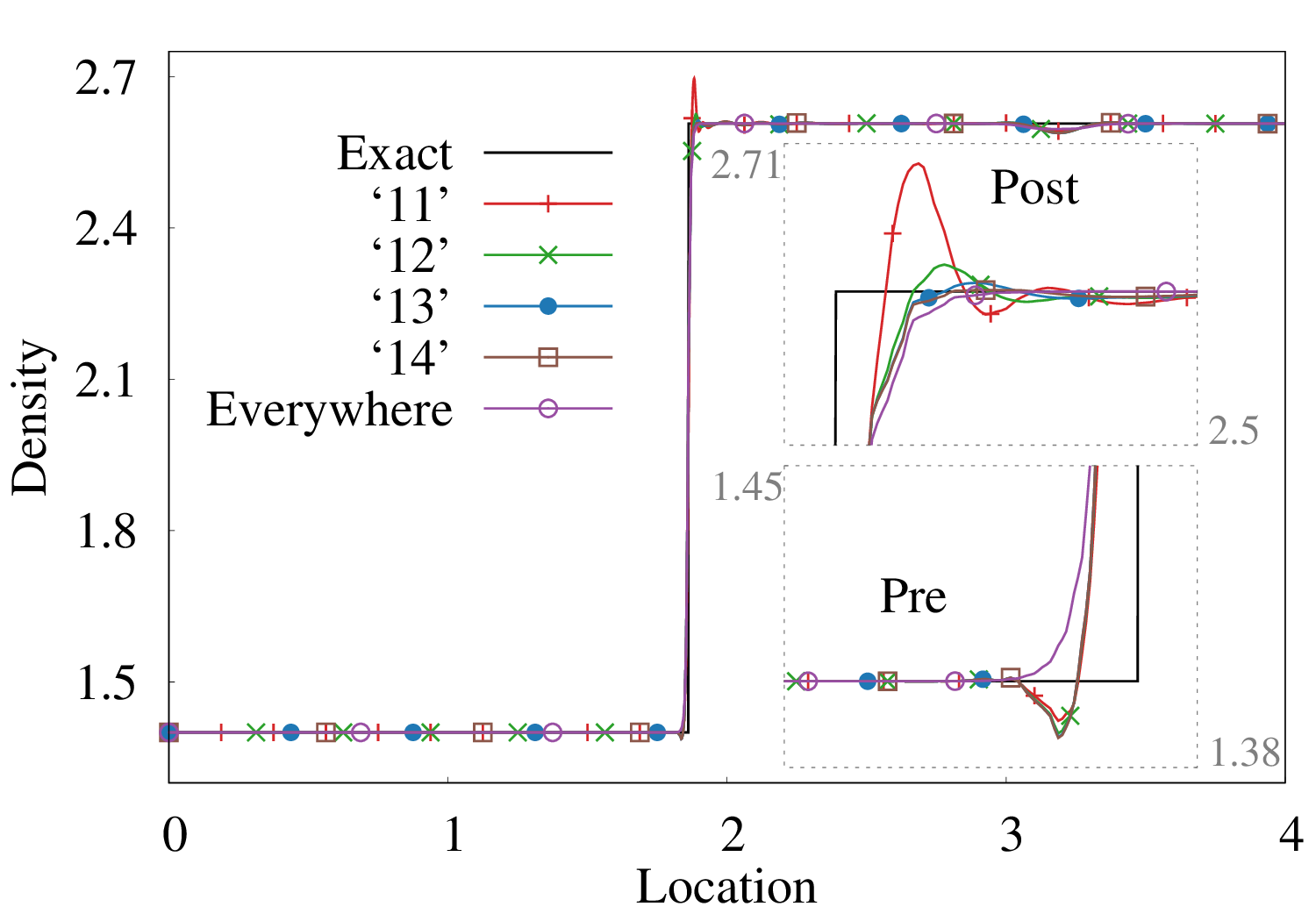}\label{fig:NA30_density_pre1}}
\subfloat[\ang{30}]{\includegraphics[width=0.5\linewidth]{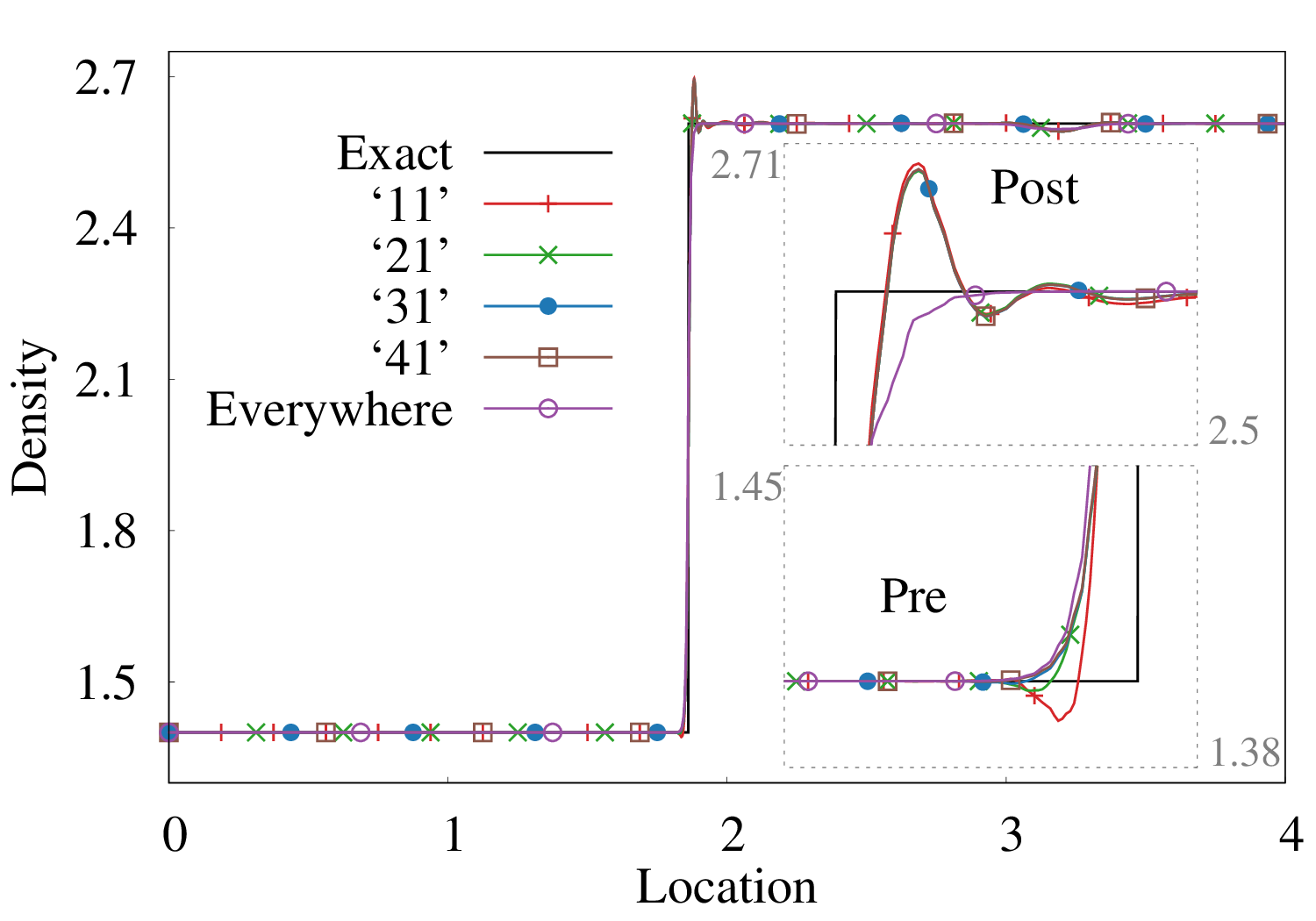}\label{fig:NA30_density_post1}}\\
\subfloat[\ang{40}]{\includegraphics[width=0.5\linewidth]{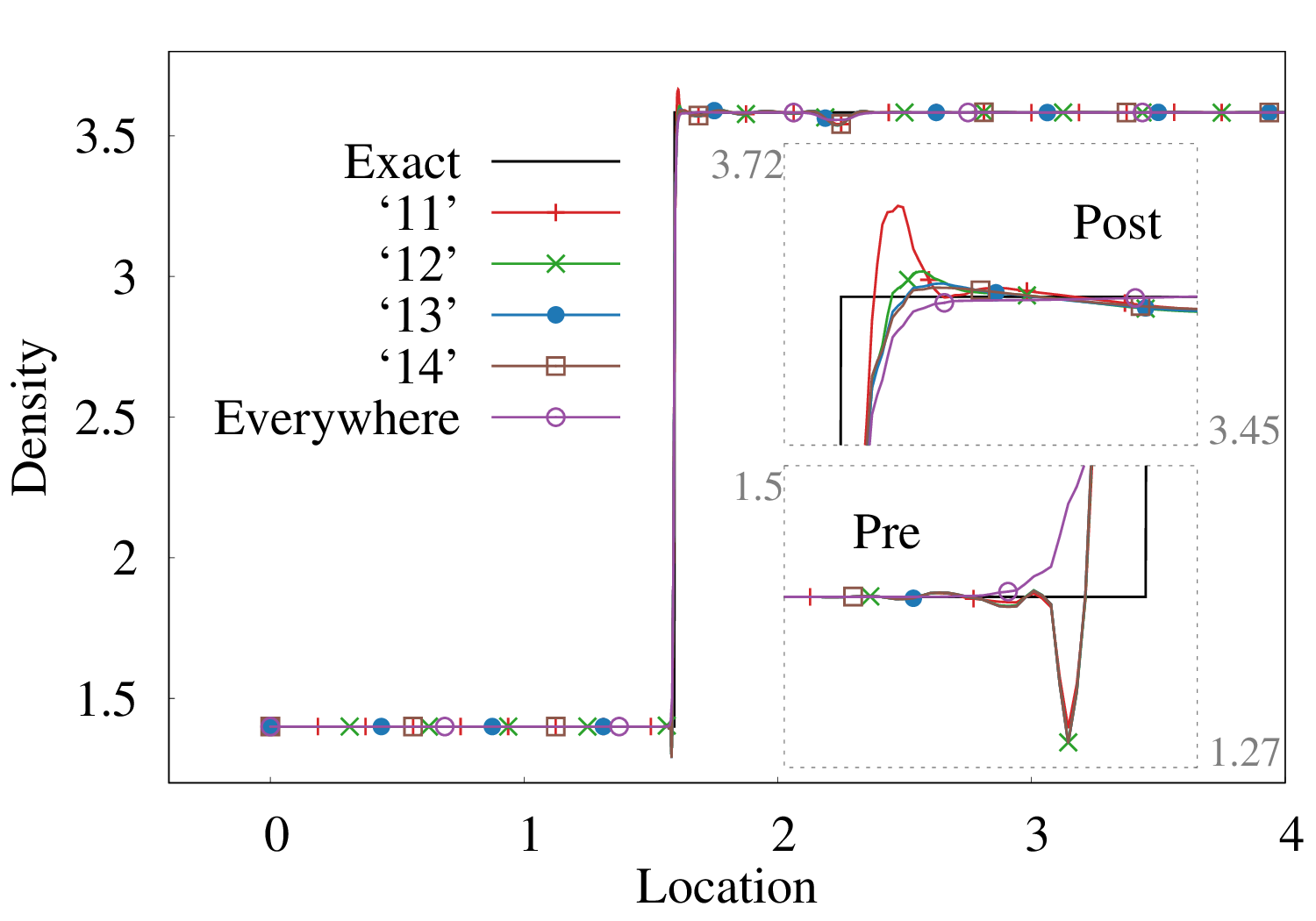}\label{fig:NA40_density_pre1}}
\subfloat[\ang{40}]{\includegraphics[width=0.5\linewidth]{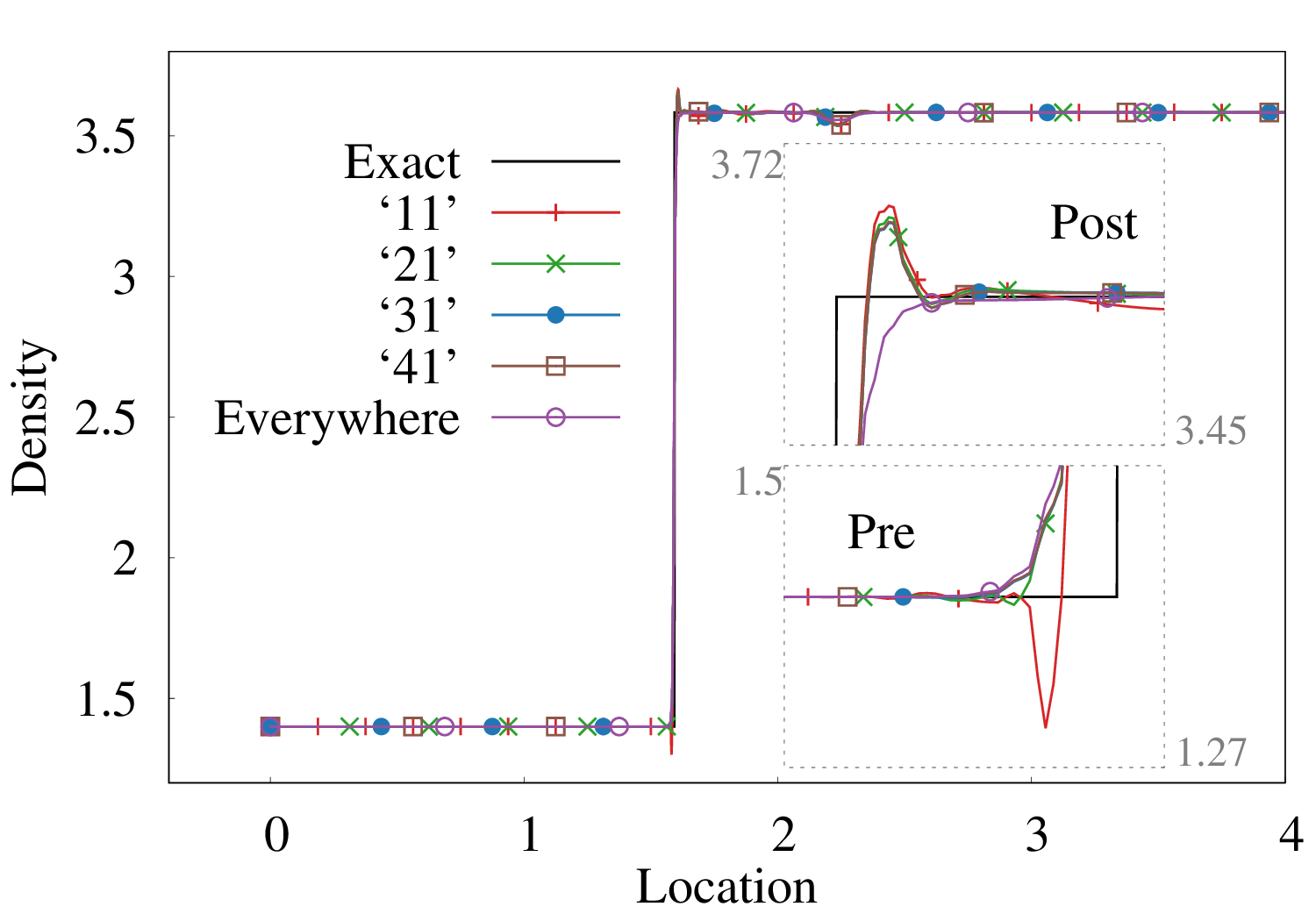}\label{fig:NA40_density_post1}}\\
\subfloat[\ang{50}]{\includegraphics[width=0.5\linewidth]{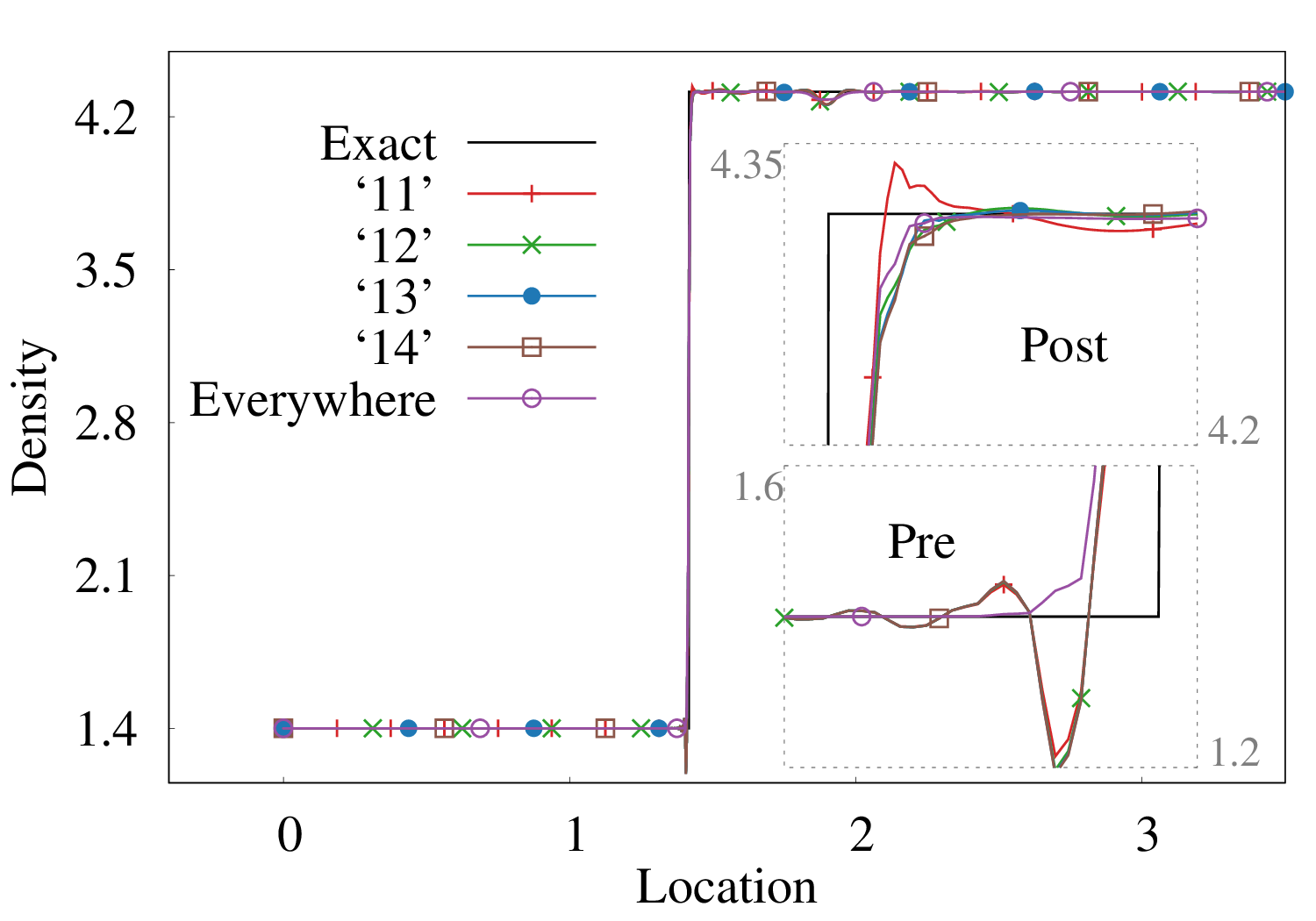}\label{fig:NA50_density_pre1}}
\subfloat[\ang{50}]{\includegraphics[width=0.5\linewidth]{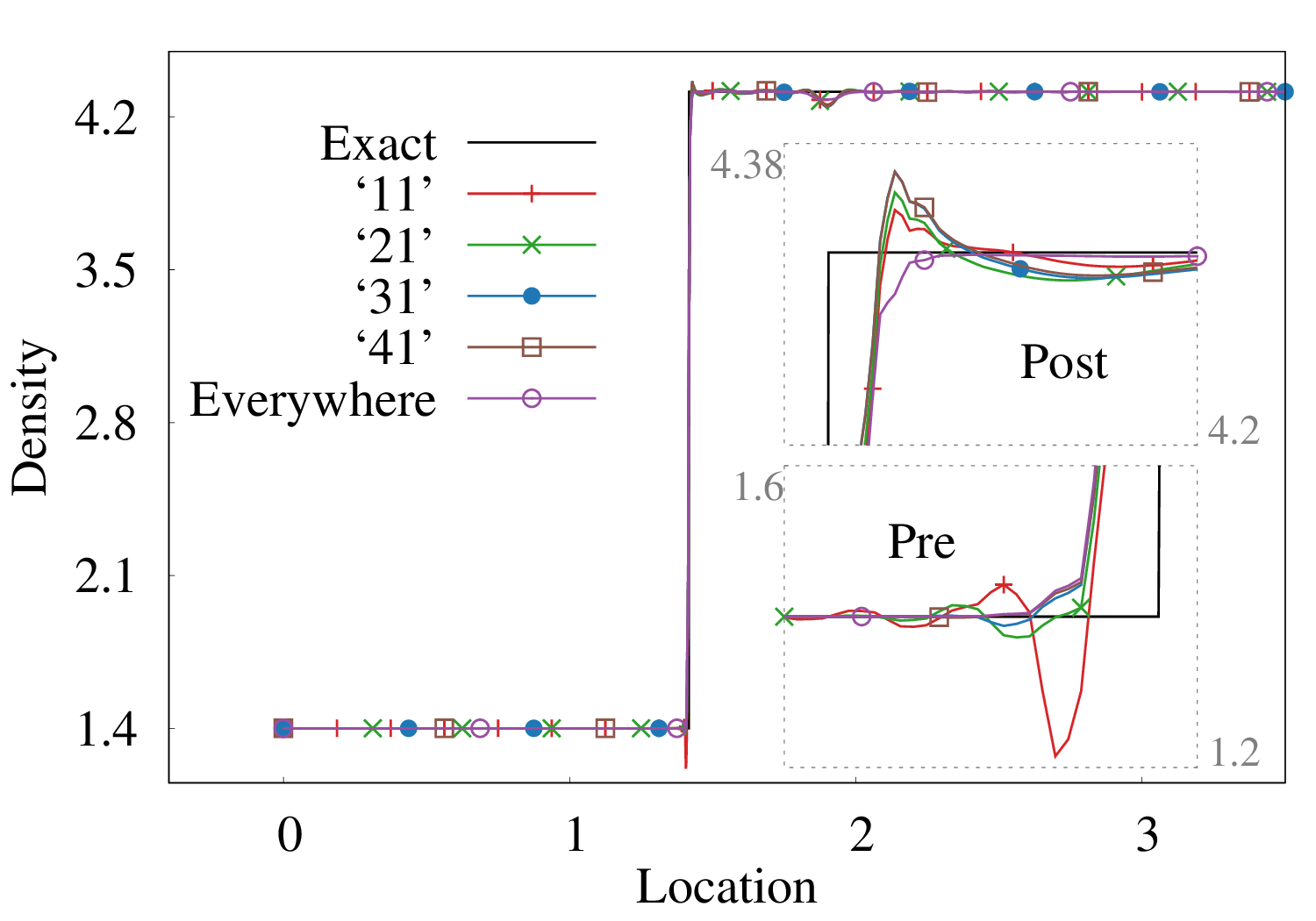}\label{fig:NA50_density_post1}}
\caption{Non-aligned oblique shock with upstream Mach number 3. Density profiles along the line $y = 0.5$ for three different shock angles on the medium grid. Top row (a, b): \ang{30}. Middle row (c, d): \ang{40}. Bottom row (e, f): \ang{50}. Left column (a, c, e): configurations with a single parallel line in the pre-shock. Right column (b, d, f): configurations with a single parallel line in the post-shock.}
\label{fig:NAOS_density_PreAndPost1}
\end{figure}

Similar results can be observed in Figure (\ref{fig:NAOS_density_PreAndPost1}) for nonaligned shocks when the limiter is applied only to cells identified by configurations with a single parallel line either in the pre- (`1x') or post-shock (`x1') region. Figures (\ref{fig:AOS_density_PreAndPost1}) and (\ref{fig:NAOS_density_PreAndPost1}) also show that adding extra troubled-cells in the post-shock region, while keeping one troubled-cell in the pre-shock region fixed, does not significantly affect the solution in pre-shock region. Similarly, adding extra troubled-cells in the pre-shock region, while keeping one troubled-cell in the post-shock region fixed, does not significantly affect the solution in the post-shock region.

\begin{figure}
\centering
\includestandalone[width=0.25\textwidth]{./6-Results/AlignOS/SameMach/legend}\\[-5ex]
\includegraphics[width=\textwidth]{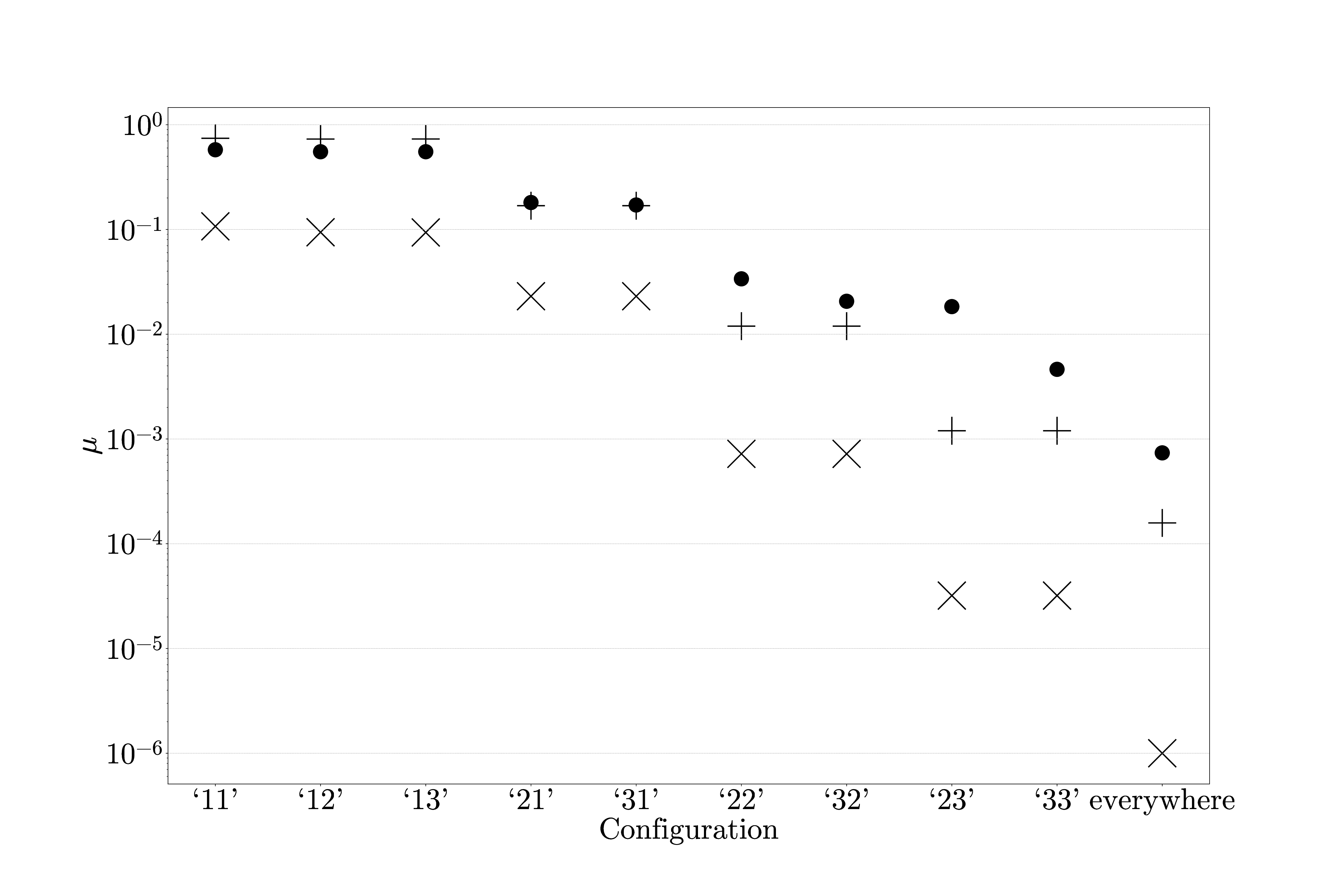}\label{fig:SameMachDensity}
\caption{Aligned oblique shock with upstream Mach number 3: $\mu$ for three different shock angles across all configurations, including the limiting-everywhere case. ($\mu = \text{TV} - L_{\infty}$)}
\label{fig:AOS_SameMach_TvLinfDiff}
\end{figure}

\begin{figure}
\centering
\includestandalone[width=0.25\textwidth]{./6-Results/NonAlignOS/SameMach/legend}\\[-5ex]
\includegraphics[width=\textwidth]{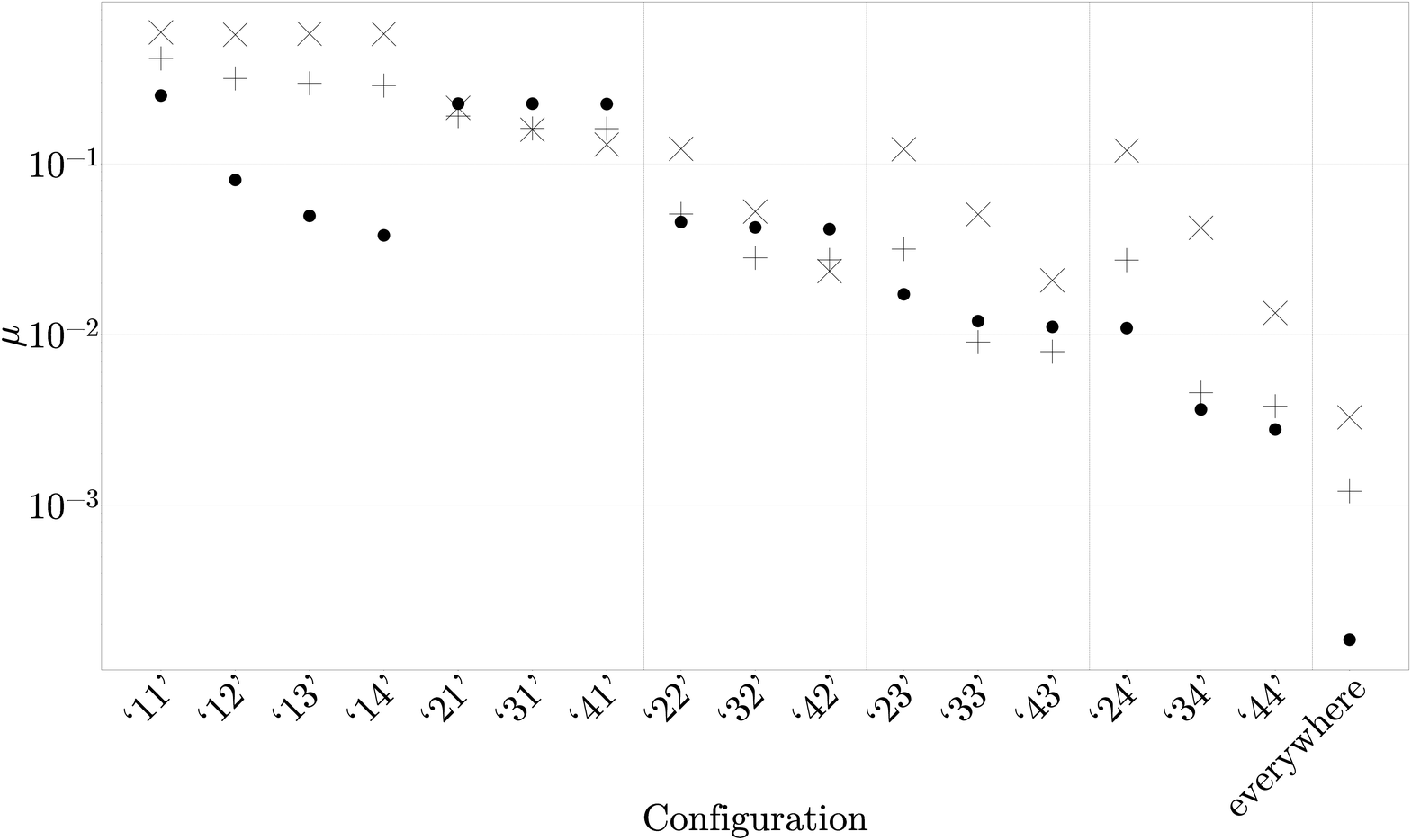}\label{fig:NAOS_SameMachDensity}
\caption{Non-aligned oblique shock with upstream Mach number 3: $\mu$ for three different shock angles across all configurations including the limiting everywhere case. ($\mu = \text{TV} - L_{\infty}$)}
\label{fig:NAOS_SameMach_TvLinfDiff}
\end{figure}

Figures (\ref{fig:AOS_SameMach_TvLinfDiff}) and (\ref{fig:NAOS_SameMach_TvLinfDiff}) further validate these observations by showing the monotonicity parameter ($\mu$) values for all configurations, including the limiting everywhere approach, for both aligned and nonaligned shocks, respectively. For configurations `1x', $\mu$ is significantly higher compared to other configurations, indicating significant oscillations in the solution for both aligned shocks and nonaligned shocks. The same holds for configurations `x1', although $\mu$ is smaller compared to the configurations `1x'. For all other configurations, with at least two troubled-cells in both the pre- and post-shock regions, $\mu$ is of the order $10 ^{-2}$ or smaller. This indicates that the solution for these configurations exhibits fewer oscillations and close to being monotonic.

\begin{figure}
\centering
\subfloat[\ang{30}]{\includegraphics[width=0.5\linewidth]{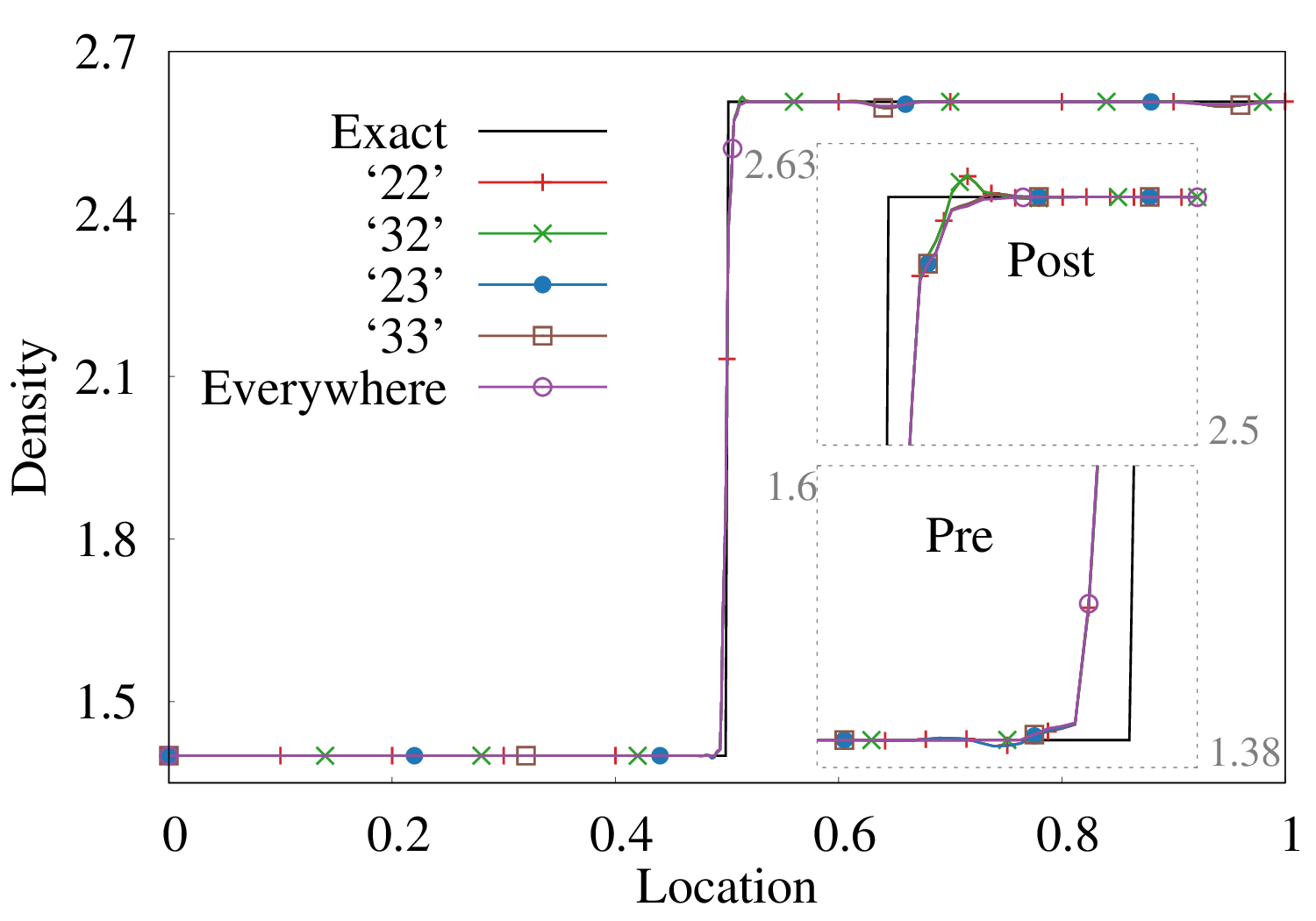}\label{fig:A30_density_2and3}}
\subfloat[\ang{40}]{\includegraphics[width=0.5\linewidth]{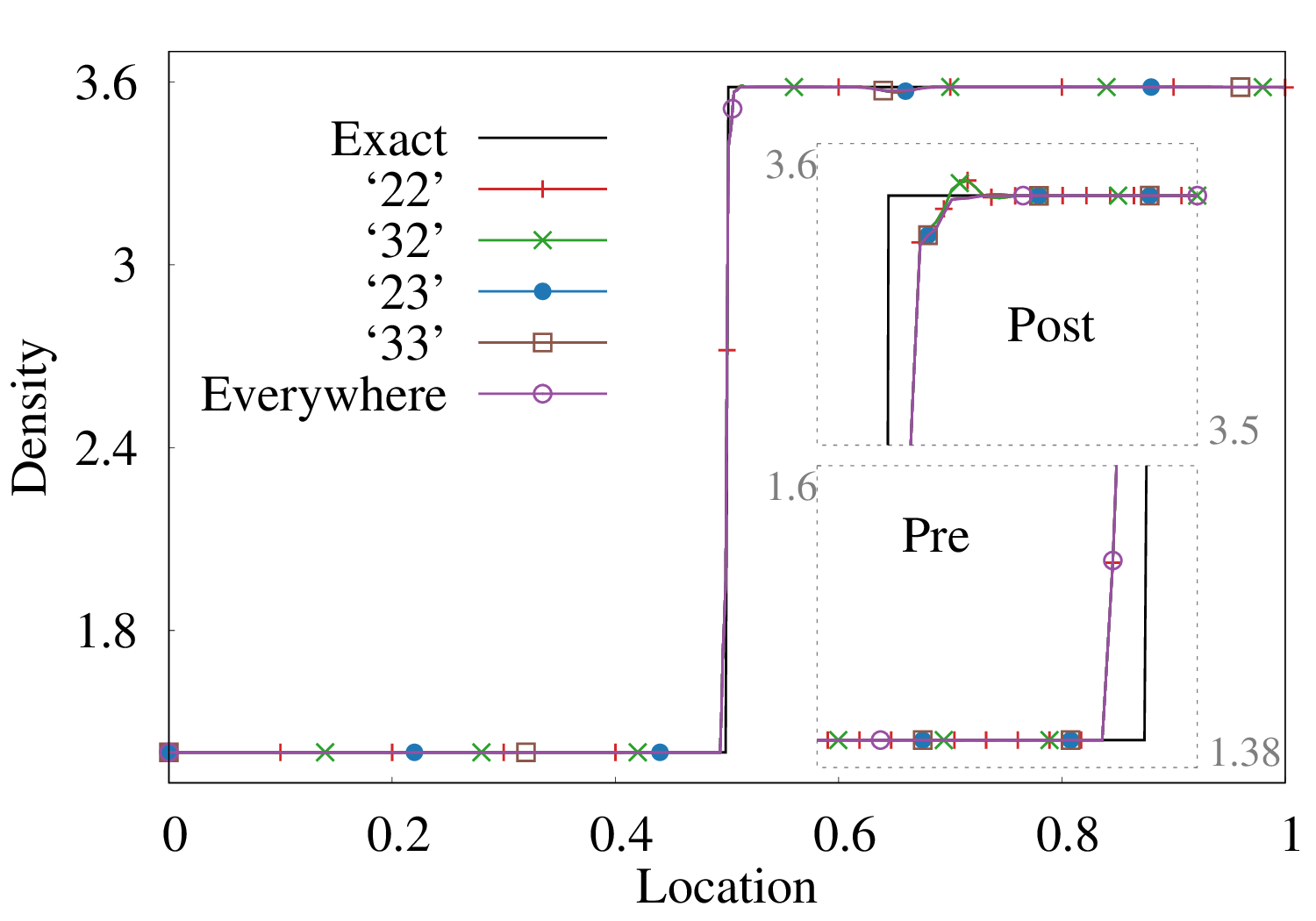}\label{fig:A40_density_2and3}}\\
\subfloat[\ang{50}]{\includegraphics[width=0.5\linewidth]{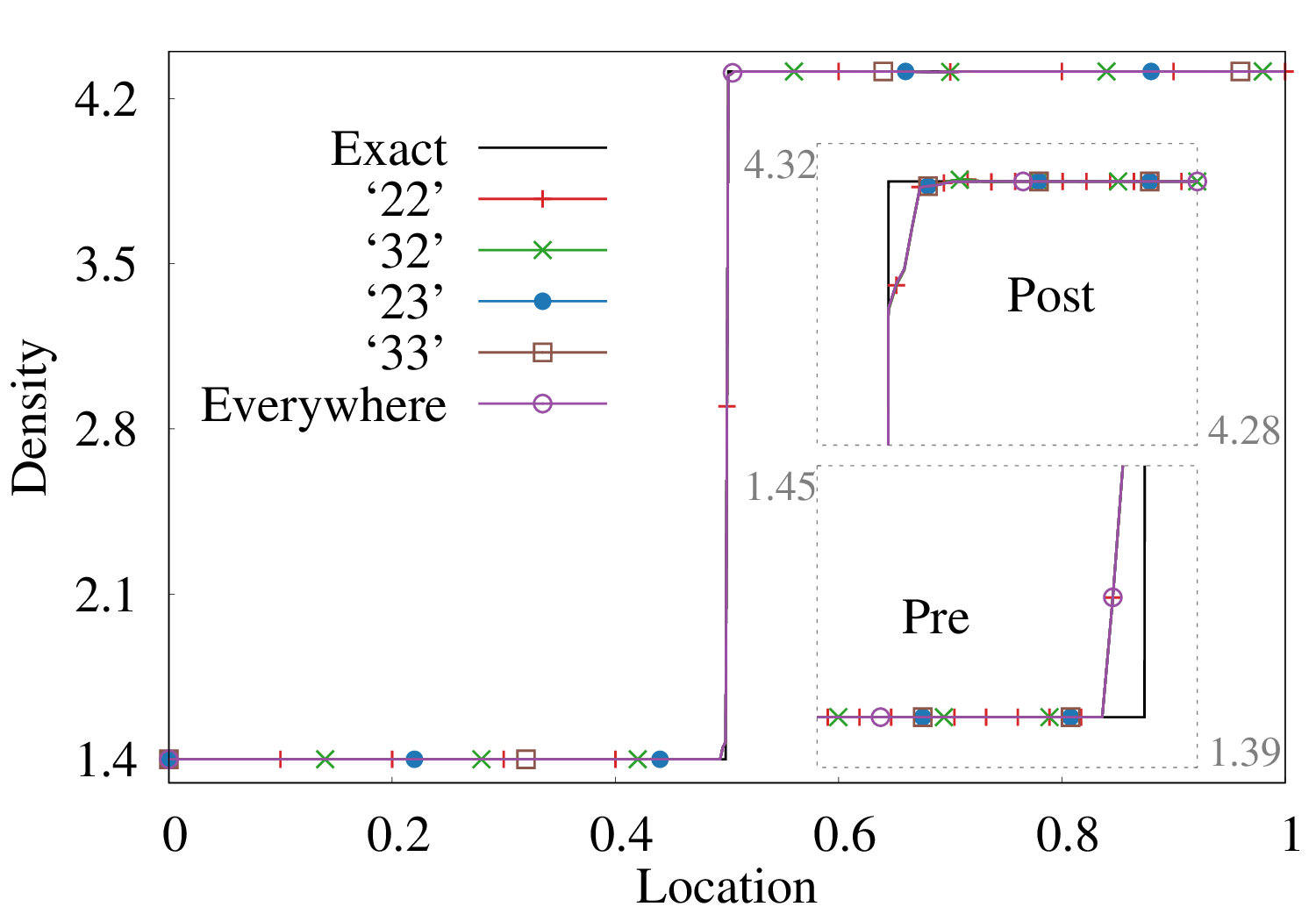}\label{fig:A50_density_2and3}}
\caption{Aligned oblique shock with upstream Mach number 3. Density profiles along the line $y = 0.5$ for three different shock angles on the medium grid for the configurations `22', `32', `23', and `33'.}
\label{fig:AOS_density_2and3}
\end{figure}

\begin{figure}
\centering
\subfloat[\ang{30}]{\includegraphics[width=0.5\linewidth]{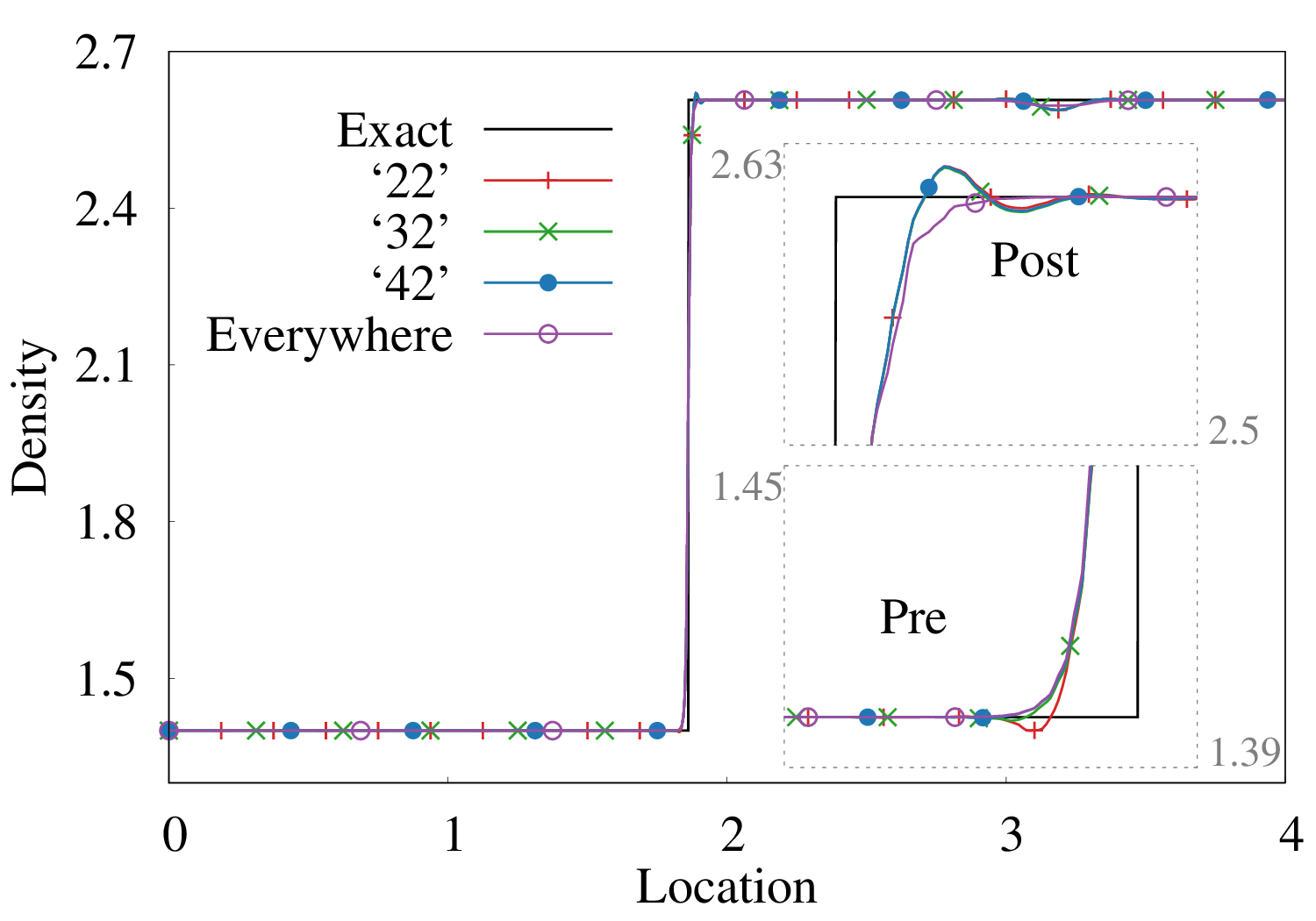}\label{fig:NA30_density_post2}}
\subfloat[\ang{30}]{\includegraphics[width=0.5\linewidth]{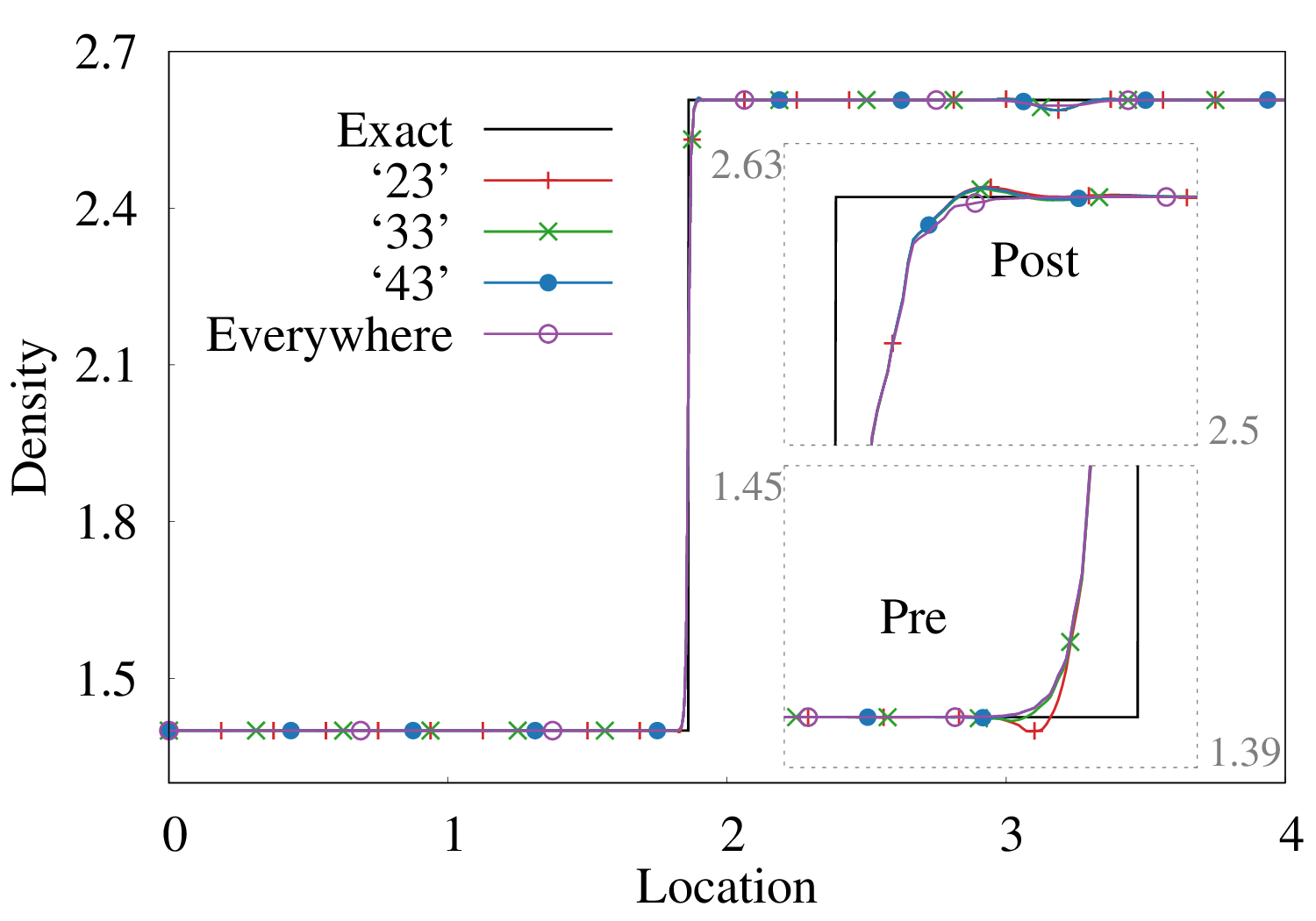}\label{fig:NA30_density_post3}}\\
\subfloat[\ang{40}]{\includegraphics[width=0.5\linewidth]{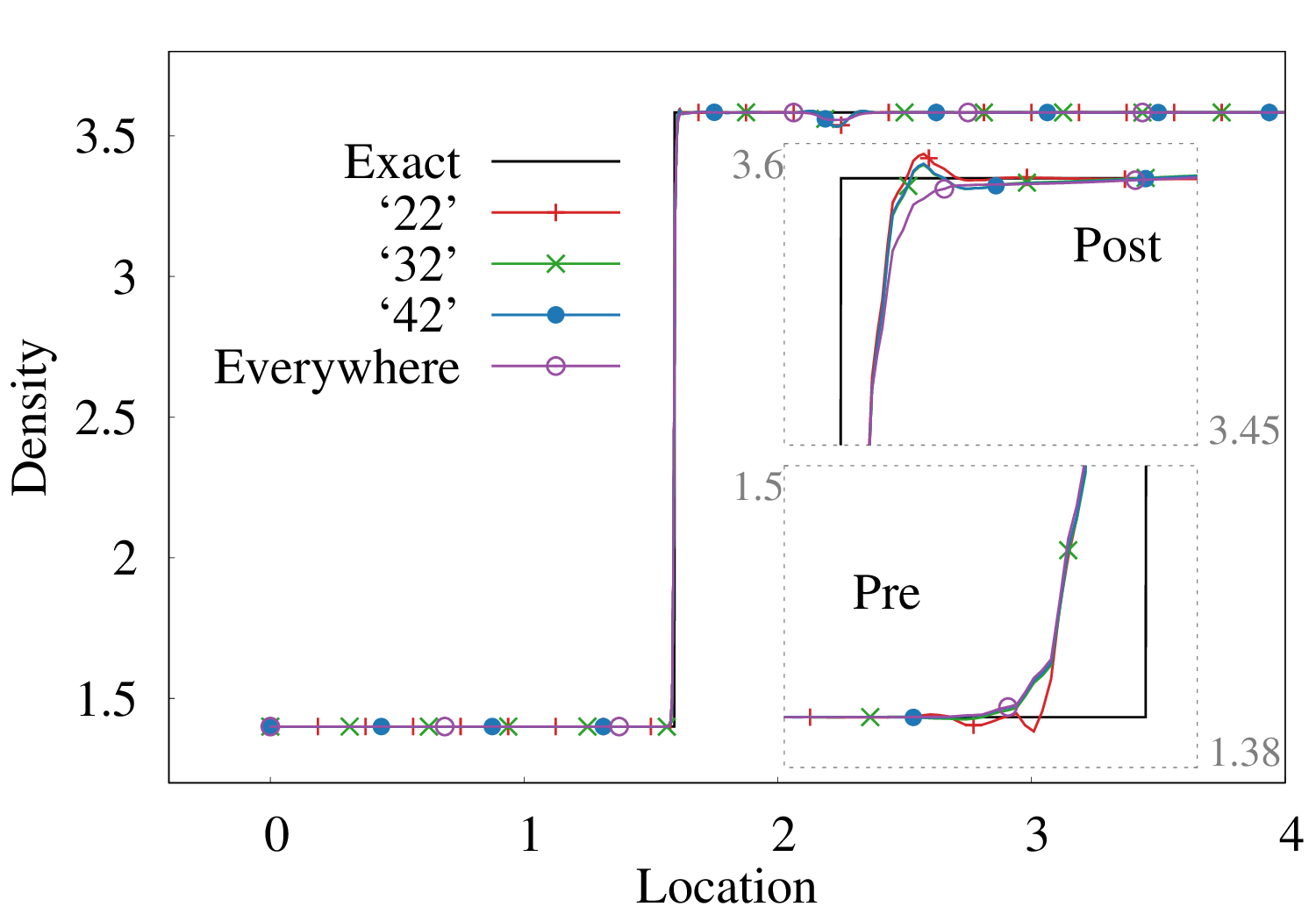}\label{fig:NA40_density_post2}}
\subfloat[\ang{40}]{\includegraphics[width=0.5\linewidth]{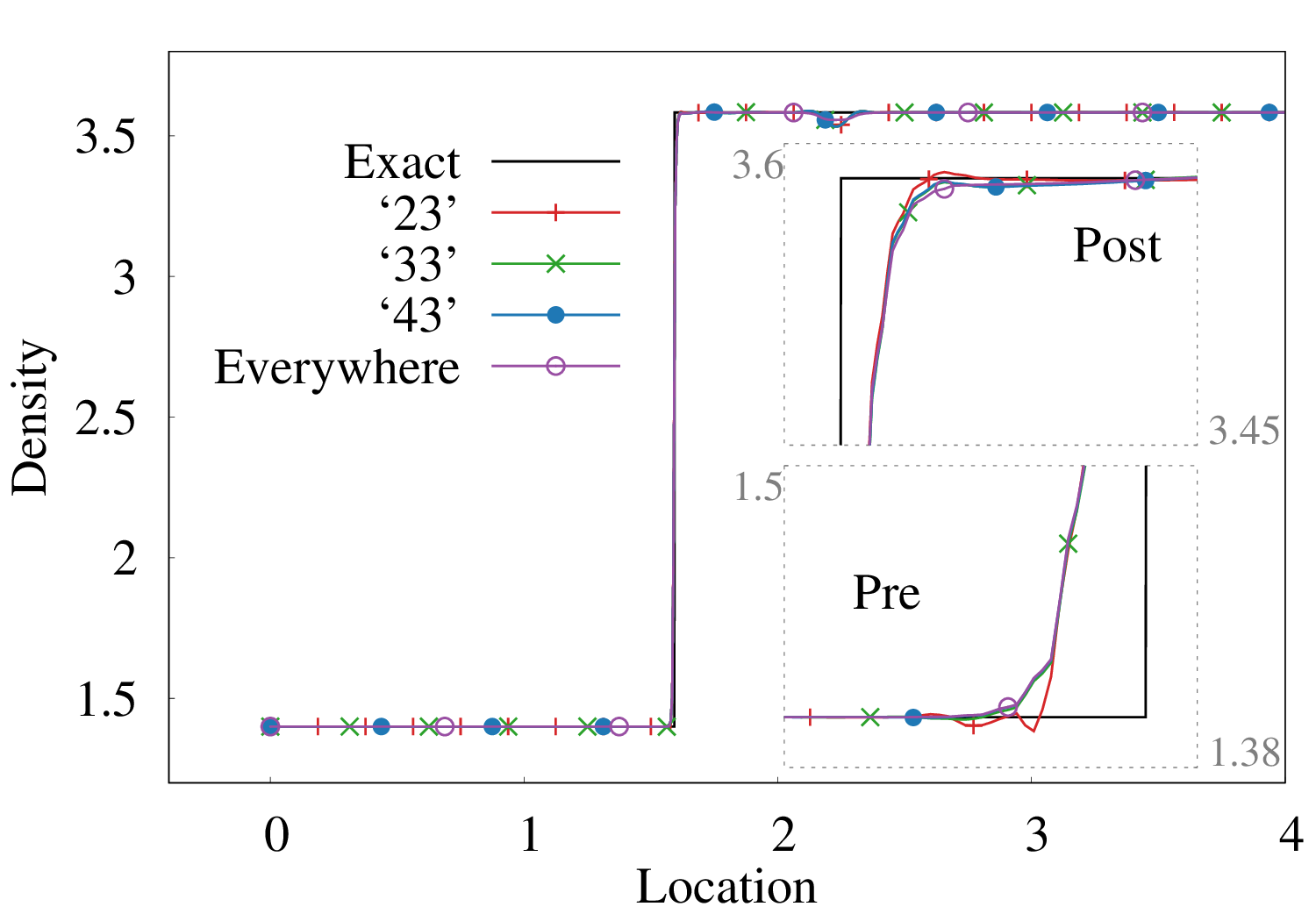}\label{fig:NA40_density_post3}}\\
\subfloat[\ang{50}]{\includegraphics[width=0.5\linewidth]{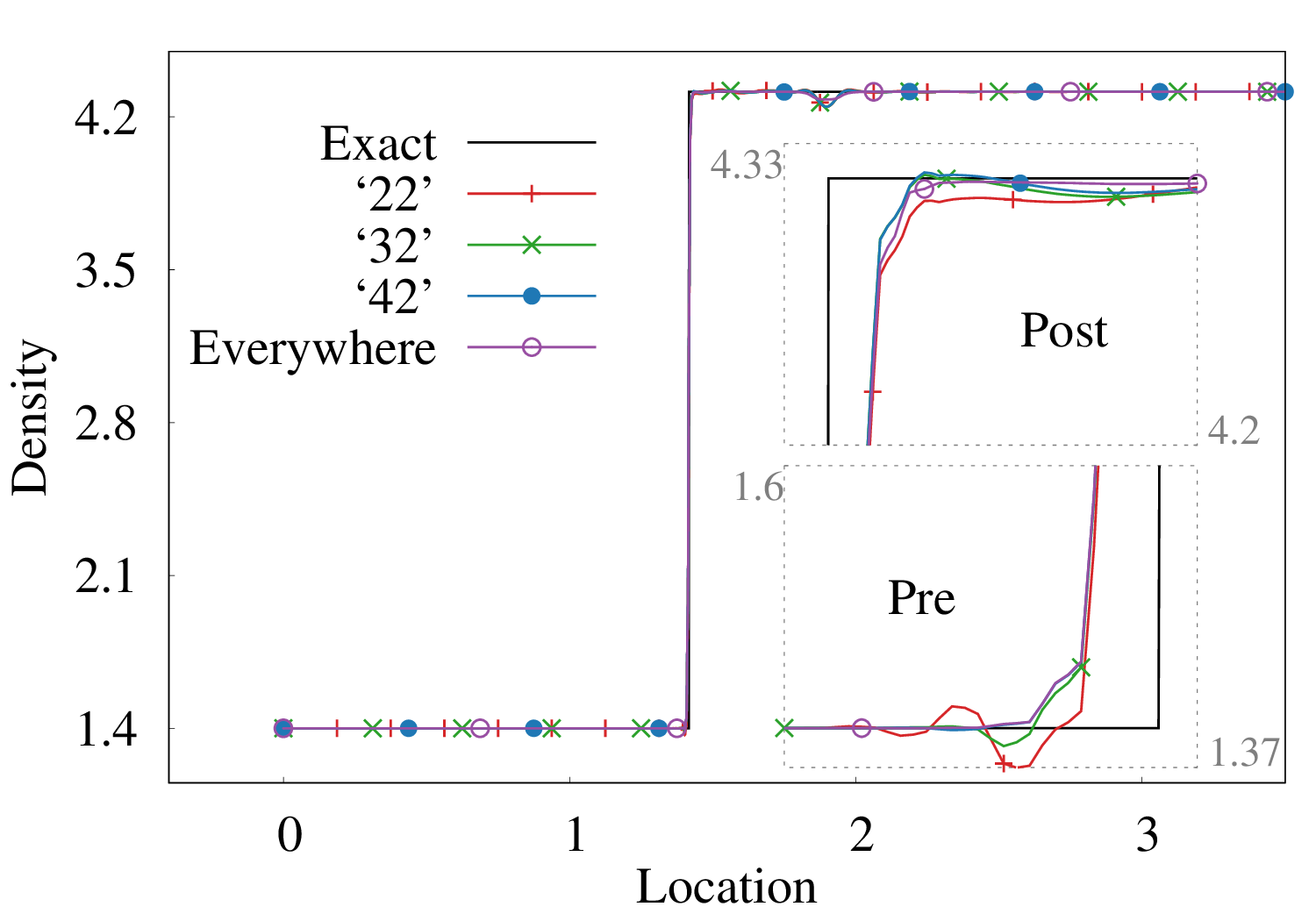}\label{fig:NA50_density_post2}}
\subfloat[\ang{50}]{\includegraphics[width=0.5\linewidth]{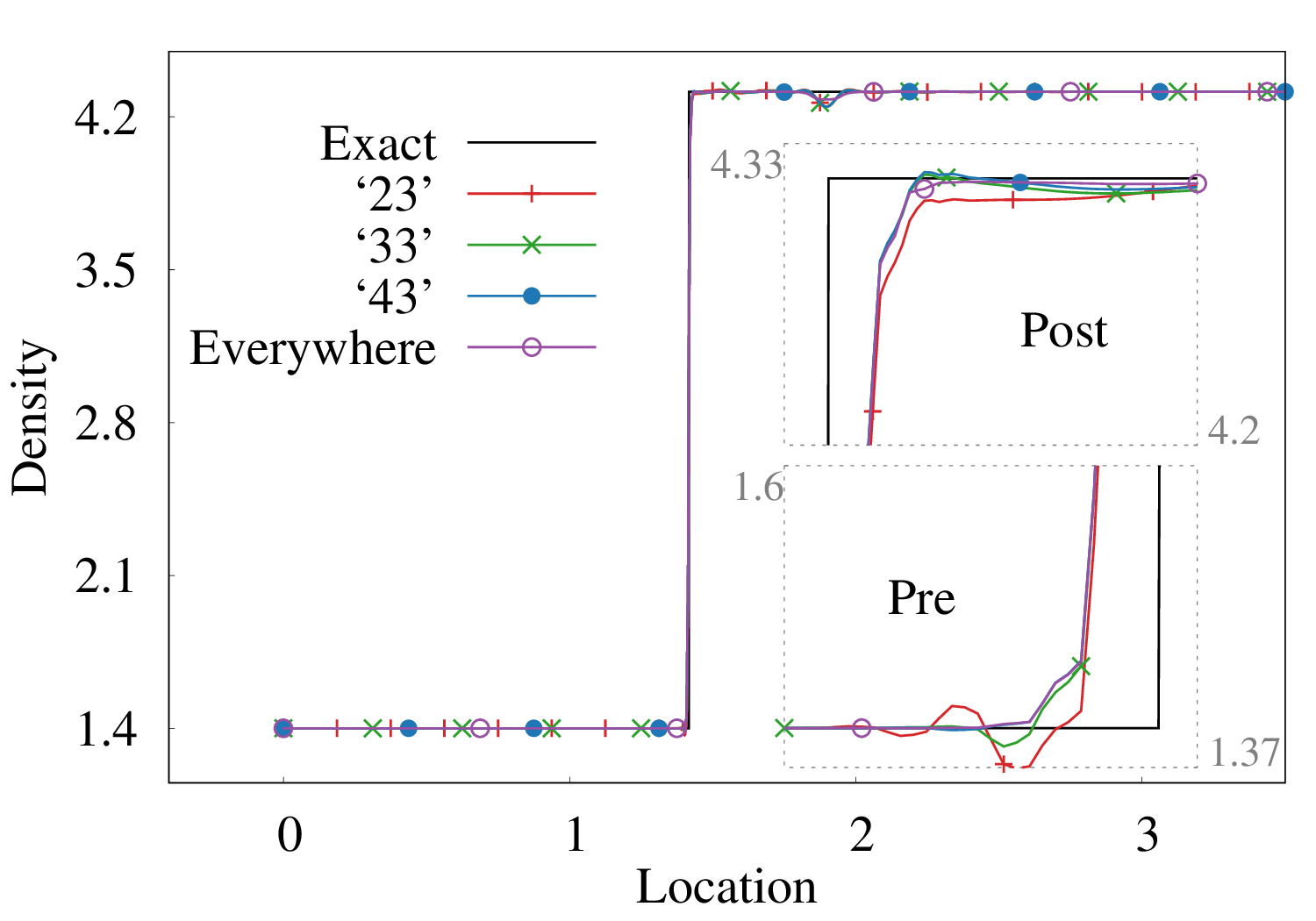}\label{fig:NA50_density_post3}}
\caption{Non-aligned oblique shock with upstream Mach number 3. Density profiles along the line $y = 0.5$ for three different shock angles on the medium grid. Top row (a, b): \ang{30}. Middle row (c, d): \ang{40}. Bottom row (e, f): \ang{50}. Left column (a, c, e): configurations with two parallel lines in the post-shock. Right column (b, d, f): configurations with three parallel lines in the post-shock.}
\label{fig:NAOS_density_Post2And3}
\end{figure}

Density profiles of the cells along the same line are shown in Figures (\ref{fig:AOS_density_2and3}) and (\ref{fig:NAOS_density_Post2And3}) for the remaining configurations of aligned and nonaligned shocks, respectively. There are no significant overshoots or undershoots in the solutions of these configurations compared to the configurations
`1x' or `x1'.

For aligned shocks, slight overshoots are observed in the solutions of configurations with two troubled-cells in the post-shock region , especially for smaller shock angles. The solutions of configurations with three troubled-cells in the post-shock region closely match the solution of limiting everywhere approach. Notably, configurations `22' and `32' yield identical solutions, as do configurations `23' and `33' for shock angles of \ang{40} and \ang{50}. However, for the \ang{30} shock angle, slight improvement in the solution from the configuration `22' to `32', as well as from `23' to `33' is observed.

These observations are also reflected in Figure (\ref{fig:AOS_SameMach_TvLinfDiff}), where $\mu$ is exactly the same for configurations `22' and `32' for shock angles of \ang{40} and \ang{50}, as well as for configurations `23' and `33'. For a shock angle of \ang{30}, $\mu$ is slightly reduced from configuration `22' to `23' as well as from configuration `23' to `33'. Overall, two troubled-cells in the pre-shock region seem to be enough for higher shock angles. However, for smaller shock angles, limiting in more than two troubled-cells has the scope to improve the solution (in terms of making it monotonic). In the post-shock region, limiting in more and more troubled-cells improves the solution to be better monotonic, though three troubled-cells seem to be enough as it closely matches the limiting everywhere approach.

For nonaligned shocks, slight undershoots and overshoots are present in configurations with two parallel lines in the pre-shock region or two parallel lines in the post-shock region. More oscillations are observed in the pre-shock region for higher shock angles, and in the post-shock region for smaller shock angles. This suggests that a greater number of troubled-cells are needed in the pre-shock region for higher shock angles, while a greater number of troubled-cells are required in the post-shock region for smaller shock angles.

These observations can also be verified from Figure (\ref{fig:NAOS_SameMach_TvLinfDiff}). For \ang{50} shocks, there is not much change in the $\mu$ for configurations `22', `23', and `24' and similarly for configurations `32', `33', and `34'.  For \ang{30} shocks, there is a slight change in the $\mu$ from configurations with 2 in the pre-shock region to 3, but not much change from configurations with 3 in the pre-shock region to 4. Overall, for higher shock angles more number of troubled-cells are needed in the pre-shock region and for smaller shock angles more number of troubled-cells are needed in the post-shock region.

\begin{figure}
\centering
\subfloat[\ang{30}]{\includegraphics[width=0.5\linewidth]{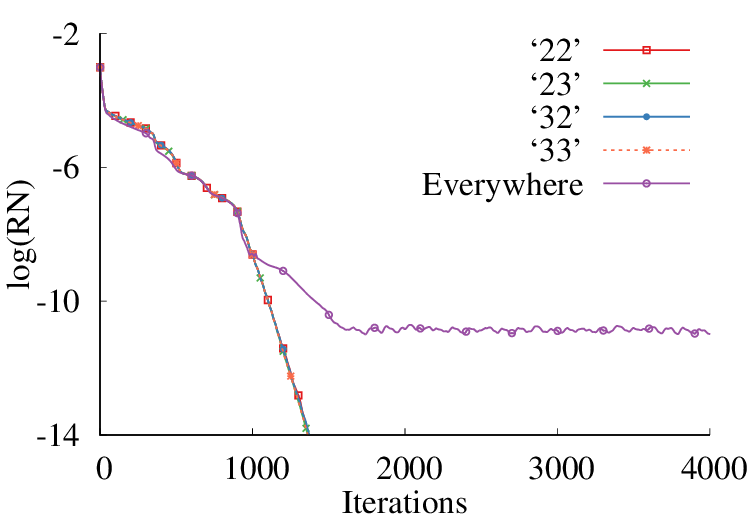}\label{fig:AOS_RN_30}}
\subfloat[\ang{40}]{\includegraphics[width=0.5\linewidth]{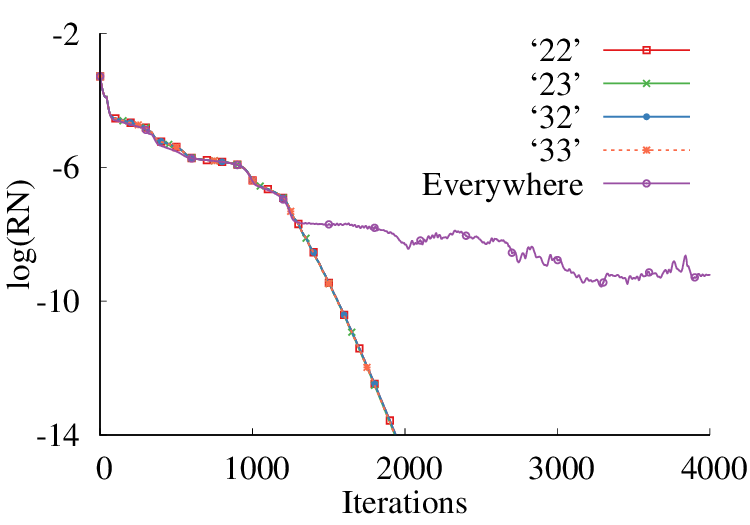}\label{fig:AOS_RN_40}}\\
\subfloat[\ang{50}]{\includegraphics[width=0.5\linewidth]{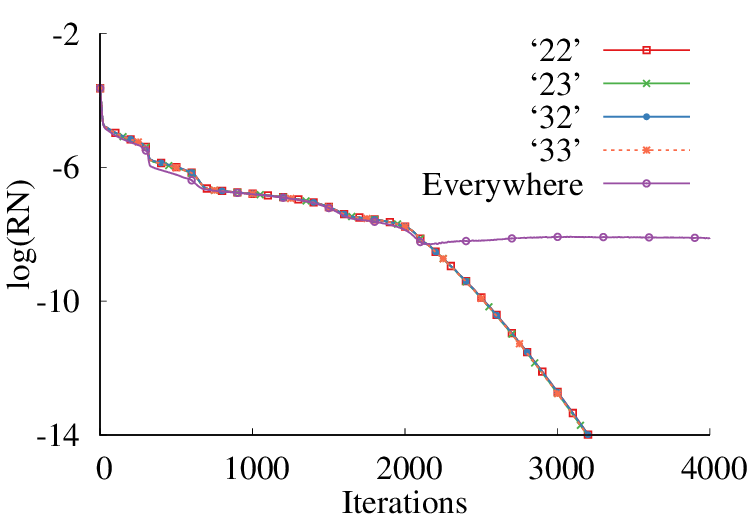}\label{fig:AOS_RN_50}}
\caption{Aligned oblique shock with inlet Mach 3. The convergence history of the residual norm as a function of number of iterations for configurations `22', `23', `32', and `33', and for limiting everywhere case.}
\label{fig:AlignedOS_M3_RN_Main}
\end{figure}

\begin{figure}
\centering
\subfloat[\ang{30}]{\includegraphics[width=0.5\linewidth]{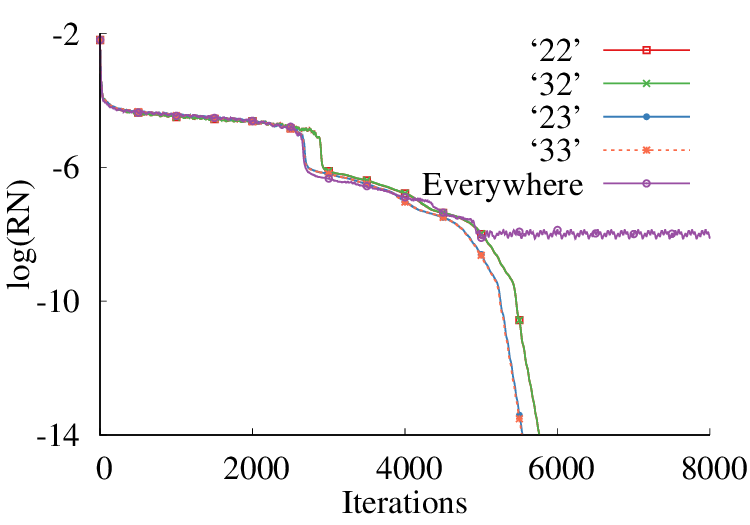}\label{fig:NAOS_RN_30}}
\subfloat[\ang{40}]{\includegraphics[width=0.5\linewidth]{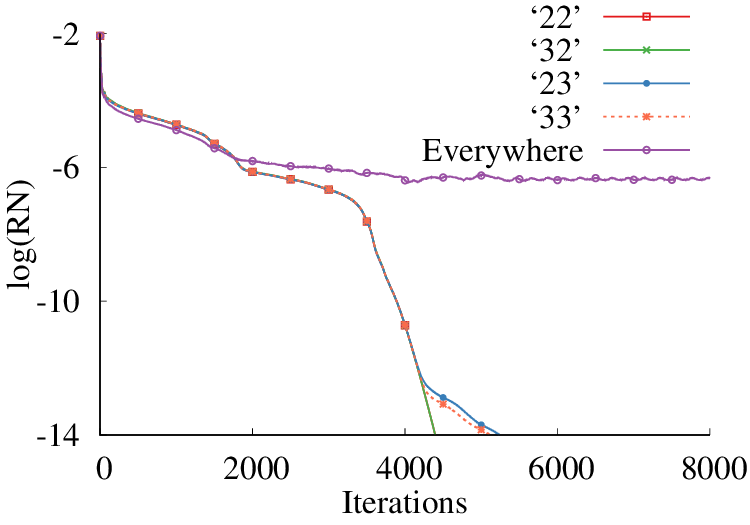}\label{fig:NAOS_RN_40}}\\
\subfloat[\ang{50}]{\includegraphics[width=0.5\linewidth]{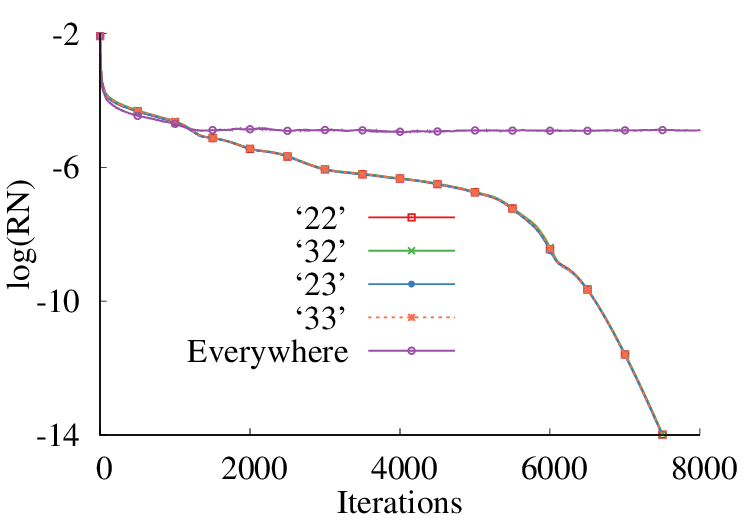}\label{fig:NAOS_RN_50}}
\caption{Non-aligned oblique shock with inlet Mach 3. The convergence history of the residual norm as a function of number of iterations for configurations `22', `23', `32', and `33', and for limiting everywhere case.}
\label{fig:NonAlignedOS_M3_RN}
\end{figure}

From here on, we present the results only for configurations with notation where the number is 2 or more both in the pre- and post-shock regions. The convergence history of the residual norm $RN$ (Eq. \ref{eq:RN}) versus number of iterations is plotted in Figures (\ref{fig:AlignedOS_M3_RN_Main}) and (\ref{fig:NonAlignedOS_M3_RN}) for configurations `22', `23', `32', and `33', including limiting everywhere approach for both aligned and nonaligned shocks. It can be observed that the solution fully converges (i.e., $RN < 10^{-14}$) to steady-state solution in nearly the same number of iterations for all these configurations and is significantly better than the limiting everywhere approach. These results show that limiting only in the vicinity of the shock certainly improves convergence and also gives a solution that closely matches the limiting everywhere approach if that region is large enough.

To verify whether the earlier conclusions about the number of troubled-cells in the pre- and post-shock regions actually depend on the shock angle ($\beta$) with respect to the flow or any other shock properties such as shock strength (pressure ratio across the shock), we varied the Mach numbers for a given shock angle. We selected upstream Mach numbers based on the pressure ratios of previous test cases, i.e., the pressure ratios of an upstream Mach number of 3 and shock angles \ang{30}, \ang{40}, and \ang{50}. For an upstream Mach number of 3, the pressure ratio across the shock depends on the shock angle. For \ang{30}, it is approximately 2.46; for \ang{40}, it is approximately 4.17; and for \ang{50}, it is approximately 5.99 (actual values are given in Table \ref{tab:AOS_SameAngle}). Now, for the \ang{30} shock angle, we determined upstream Mach numbers so as to get pressure ratios of 4.17 and 5.99. Similarly, we determined upstream Mach numbers for \ang{40} and \ang{50} shock angles. For clarity, the shock angles, Mach numbers, and the corresponding pressure ratios are shown in Table (\ref{tab:AOS_SameAngle}).

\begin{table}
\centering
\caption{For a given shock angle, different Mach numbers and the pressure ratios across the shock for that combination.}
\begin{tabular}{c c c}
\toprule
Shock Angle & Mach Number & Pressure Ratio ($P_2/P_1$) \\
\midrule
\multirow{3}{*}{\ang{30}} & $M_{PR_1}$ = 3 & $PR_1$ = 2.45833333\\
& $M_{PR_2}$ = 3.85672566 & $PR_2$ = 4.17168040 \\
& $M_{PR_3}$ = 4.59626666 & $PR_3$ = 5.99498626 \\
\midrule
\multirow{3}{*}{\ang{40}} & $M_{PR_1}$ = 2.33358574 & $PR_1$ = 2.45833333\\
& $M_{PR_2}$ = 3 & $PR_2$ = 4.17168040 \\
& $M_{PR_3}$ = 3.57526078 & $PR_3$ = 5.99498626 \\
\midrule
\multirow{3}{*}{\ang{50}} & $M_{PR_1}$ = 1.958110935 & $PR_1$ = 2.45833333 \\
& $M_{PR_2}$ = 2.517298895 & $PR_2$ = 4.17168040 \\
& $M_{PR_3}$ = 3 & $PR_3$ = 5.99498626 \\
\bottomrule
\end{tabular}%
\label{tab:AOS_SameAngle}
\end{table}

\begin{figure}
\centering
\includestandalone[width=0.3\textwidth]{./6-Results/AlignOS/SameAngle/legend}\\[-5ex]
\subfloat[\ang{30}]{\includegraphics[width=0.5\textwidth]{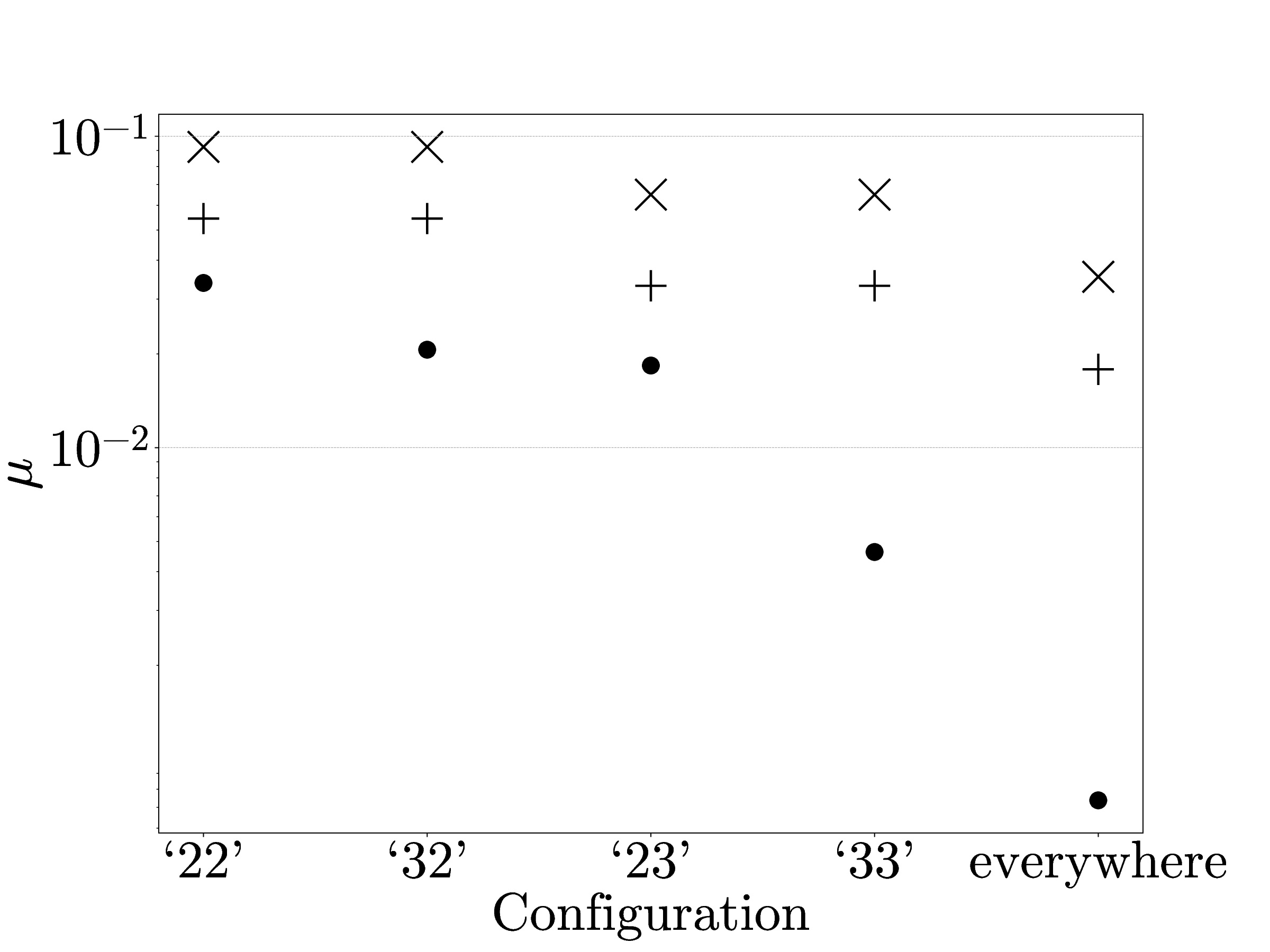}\label{fig:SameAngleDensity_30}}
\subfloat[\ang{40}]{\includegraphics[width=0.5\textwidth]{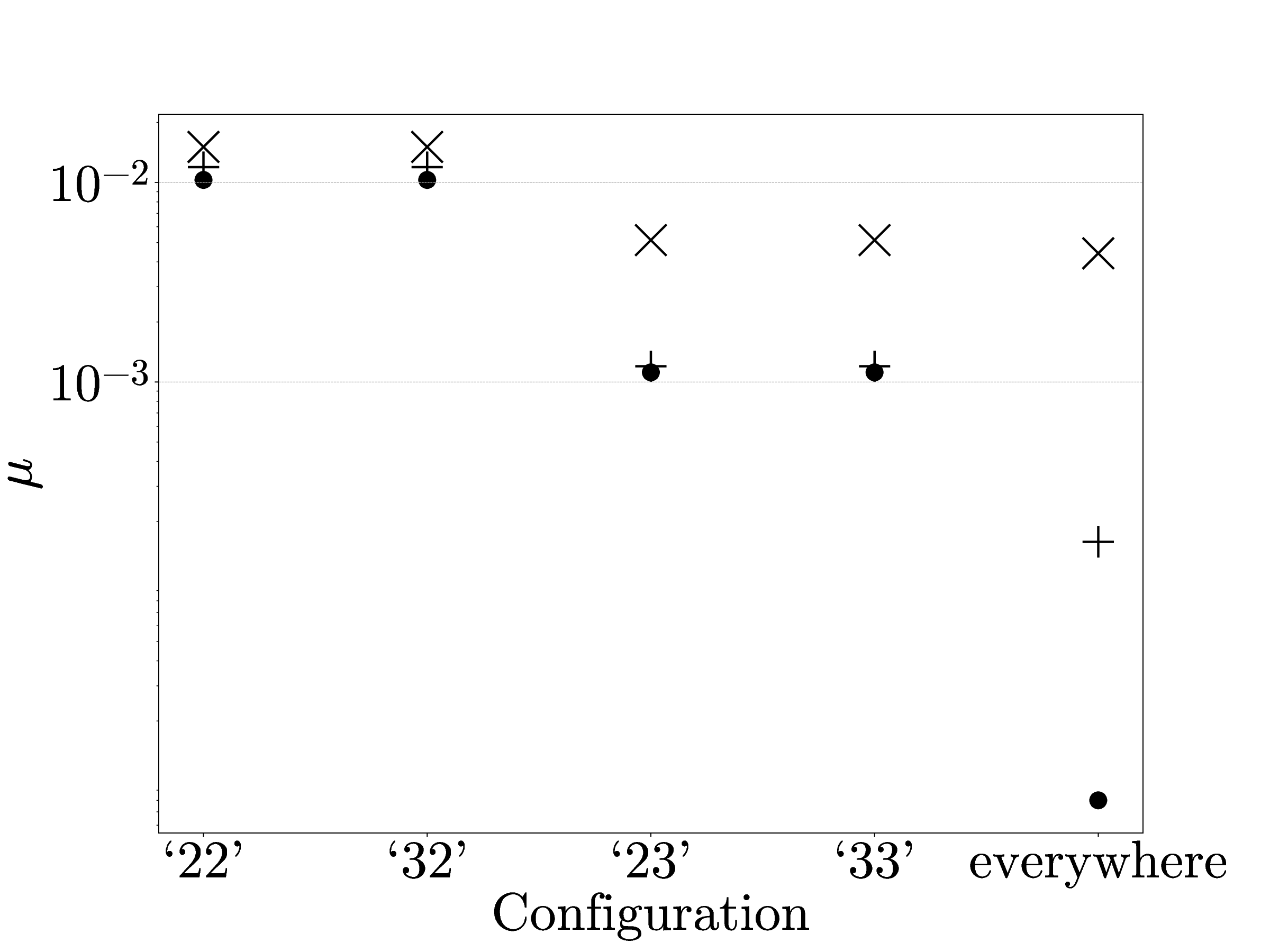}\label{fig:SameAngleDensity_40}}\\[-2ex]
\subfloat[\ang{50}]{\includegraphics[width=0.5\textwidth]{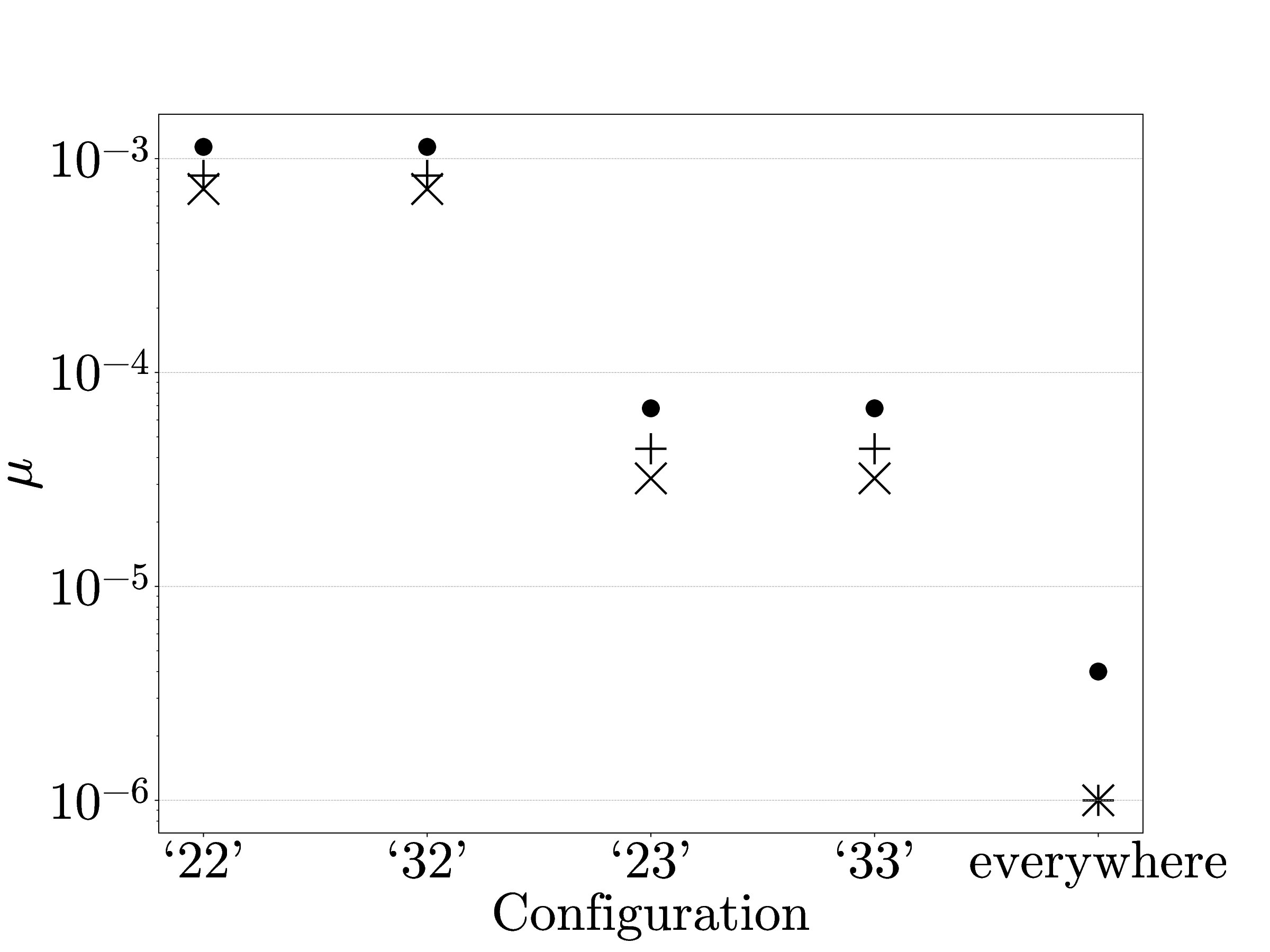}\label{fig:SameAngleDensity_50}}
\caption{Aligned oblique shock: $\mu$ for three different upstream Mach numbers for a given shock angle. For upstream Mach numbers $M_{PR_1}$, $M_{PR_2}$, and $M_{PR_3}$ (in ascending order), refer Table (\ref{tab:AOS_SameAngle}). ($\mu = \text{TV} - L_{\infty}$)}
\label{fig:AOS_SameAngle_TvLinfDiff}
\end{figure}

For aligned shocks, Figure (\ref{fig:AOS_SameAngle_TvLinfDiff}) presents the $\mu$ for configurations `22', `32', `23', and `33' for three different upstream Mach numbers for a given shock angle. For \ang{30} shock (Figure \ref{fig:SameAngleDensity_30}), except for Mach 3, adding extra troubled-cells in the pre-shock region gives exactly the same $\mu$ for other two higher Mach numbers. This indicates that as the upstream Mach number increases for a given shock angle, i.e., pressure ratio across the shock increases, extra troubled-cells in the pre-shock region does not affect the solution. However, for \ang{40} (Figure \ref{fig:SameAngleDensity_40}) and \ang{50} (Figure \ref{fig:SameAngleDensity_50}) shocks, even for smaller pressure ratios, the earlier conclusions holds true: adding extra troubled-cells in the pre-shock does not effect the solution and adding extra troubled-cells in the post-shock region improves the solution. Overall, the requirement for the number of troubled-cells in the pre-and post-shock depends on both the shock angle and the pressure ratio across the shock.

\begin{figure}
\centering
\includestandalone[width=0.3\textwidth]{./6-Results/AlignOS/SameAngle/legend}\\[-2ex]
\subfloat[\ang{30}]{\includegraphics[width=0.5\textwidth]{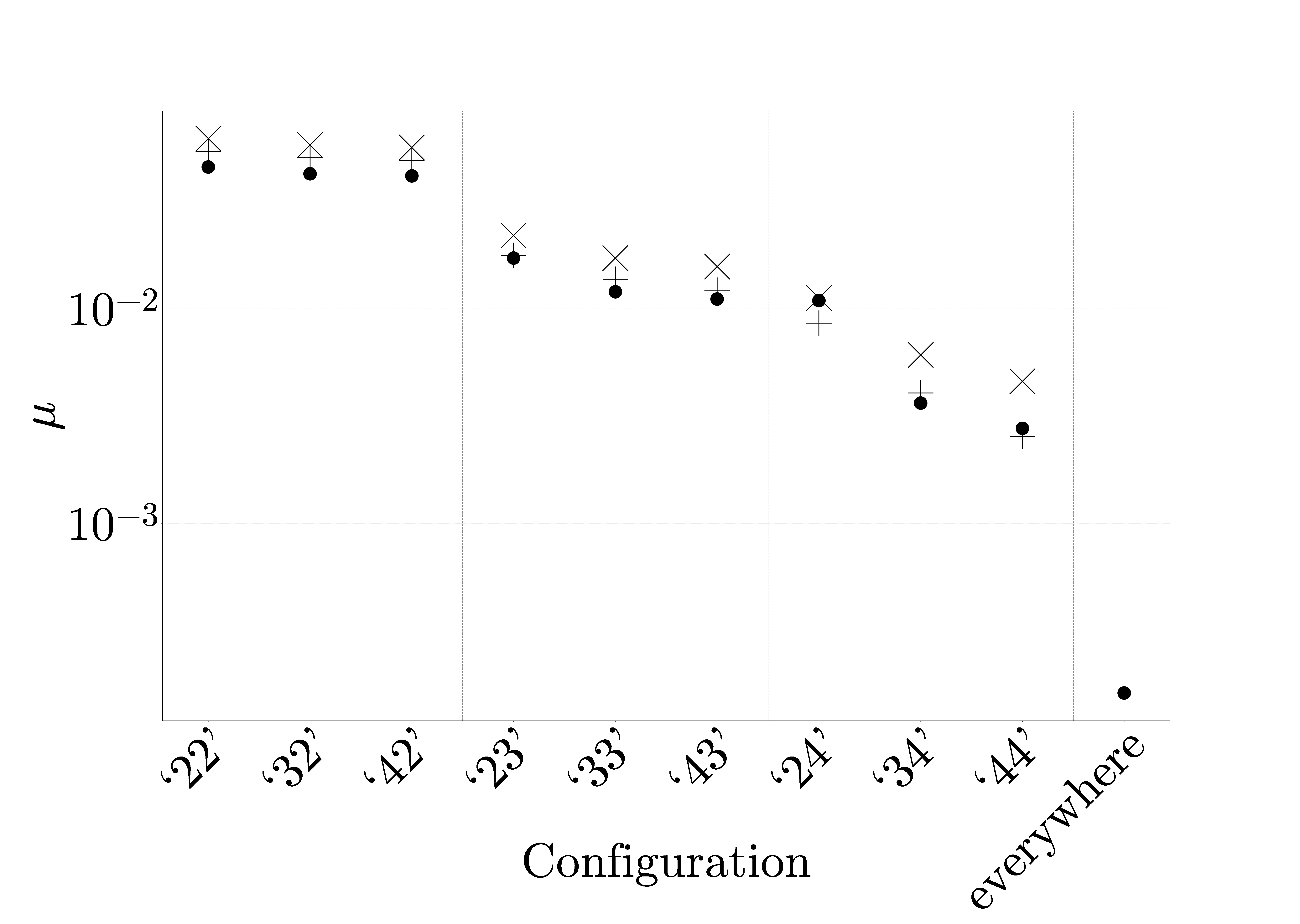}\label{fig:NAOS_SameAngleDensity_30}}
\subfloat[\ang{40}]{\includegraphics[width=0.5\textwidth]{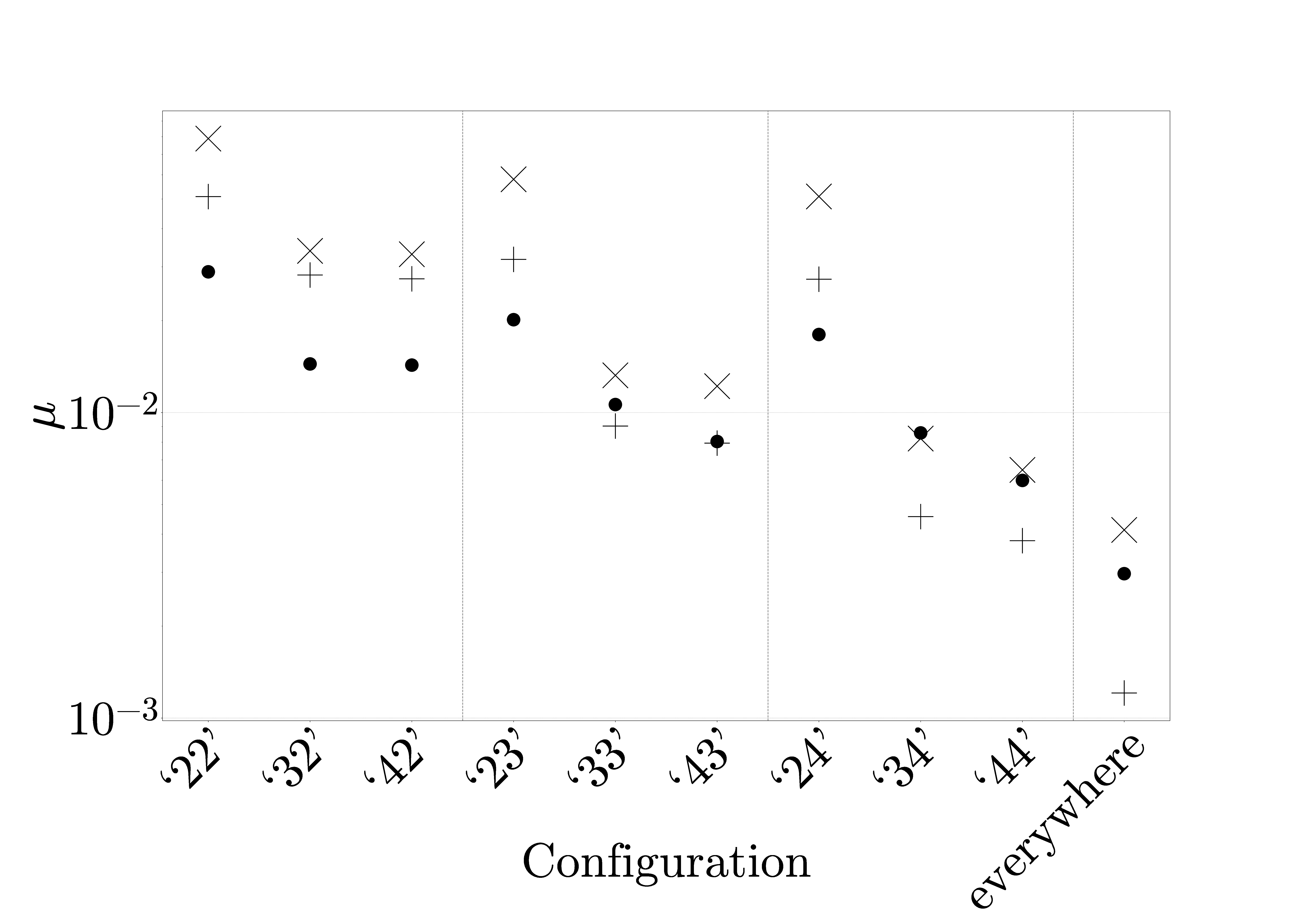}\label{fig:NAOS_SameAngleDensity_40}}\\[-2ex]
\subfloat[\ang{50}]{\includegraphics[width=0.5\textwidth]{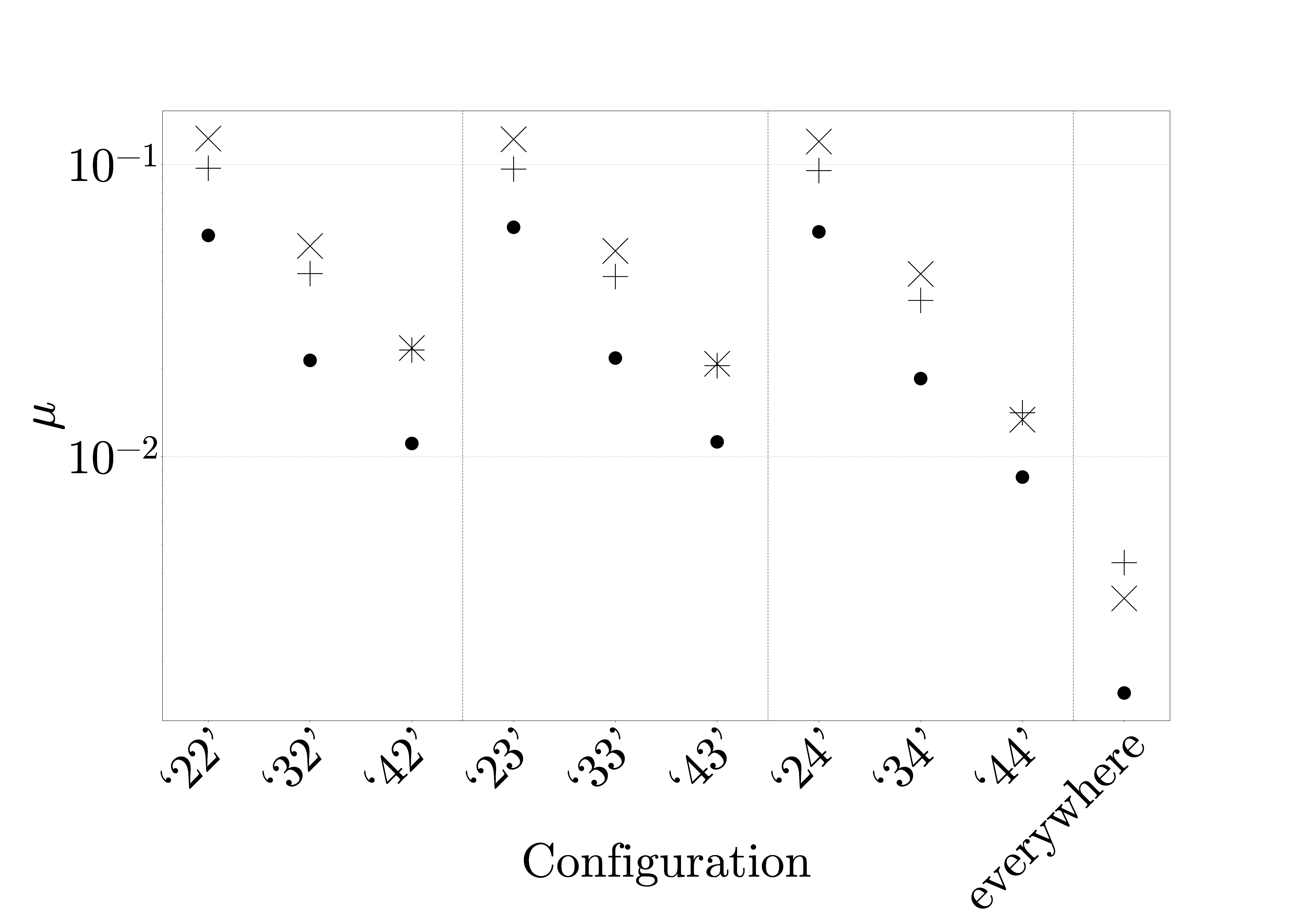}\label{fig:NAOS_SameAngleDensity_50}}
\caption{Non-aligned oblique shock: $\mu$ for three different upstream Mach numbers for a given shock angle. For upstream Mach numbers $M_{PR_1}$, $M_{PR_2}$, and $M_{PR_3}$ (in ascending order), refer Table (\ref{tab:AOS_SameAngle}). ($\mu = \text{TV} - L_{\infty}$)}
\label{fig:NAOS_SameAngle_TvLinfDiff}
\end{figure}

For nonaligned shocks, Figure (\ref{fig:NAOS_SameAngle_TvLinfDiff}) presents the $\mu$ for all configurations except `1x' and `x1' for three different upstream Mach numbers for a given shock angle. For each shock angle, the solutions of the other two Mach numbers follows same behaviour as Mach 3. For nonaligned shocks, it is clearly the shock angle that matters. For the smaller shock angles, more troubled-cells are needed in the post-shock region, whereas for higher shock angles, more troubled-cells are needed in the pre-shock region.

Next, we fixed the shock at a particular angle with respect to the grid to determine whether the number of troubled-cells depends on the actual shock angle or the shock angle with respect to the grid. Table (\ref{tab:FixedMisAlignment_SameMach}) presents the upstream flow angles for three different shock angles for a fixed nonalignment of the shock with respect to the grid.


\begin{table}
\centering
\caption{For a fixed nonalignment with respect to the grid, different shock angles and the corresponding upstream flow angles. Here, the counterclockwise (CCW) direction taken as the negative angular direction and the clockwise (CW) direction as the positive angular direction.}
\begin{tabular}{c c c}
\toprule
 \makecell{Fixed nonalignment\\w.r.t grid} & \makecell{Upstream flow\\angle(in degrees)} & \makecell{Actual shock \\angle (in degrees)}\\
\midrule
\multirow{3}{*}{\ang{30}} & 0 & 30 \\
& -10 & 40 \\
& -20 & 50 \\
\midrule
\multirow{3}{*}{\ang{40}} & 10 & 30 \\
& 0 & 40 \\
& -10 & 50 \\
\midrule
\multirow{3}{*}{\ang{50}} & 20 & 30 \\
& 10 & 40 \\
& 0 & 50 \\
\bottomrule
\end{tabular}%
\label{tab:FixedMisAlignment_SameMach}
\end{table}

\begin{figure}
\centering
\includestandalone[width=0.4\textwidth]{./6-Results/NonAlignOS/FixedMisAlignment/SameMach/legend}\\[-5ex]
\subfloat[\ang{30} fixed misalignment w.r.t. grid]{\includegraphics[width=0.5\textwidth]{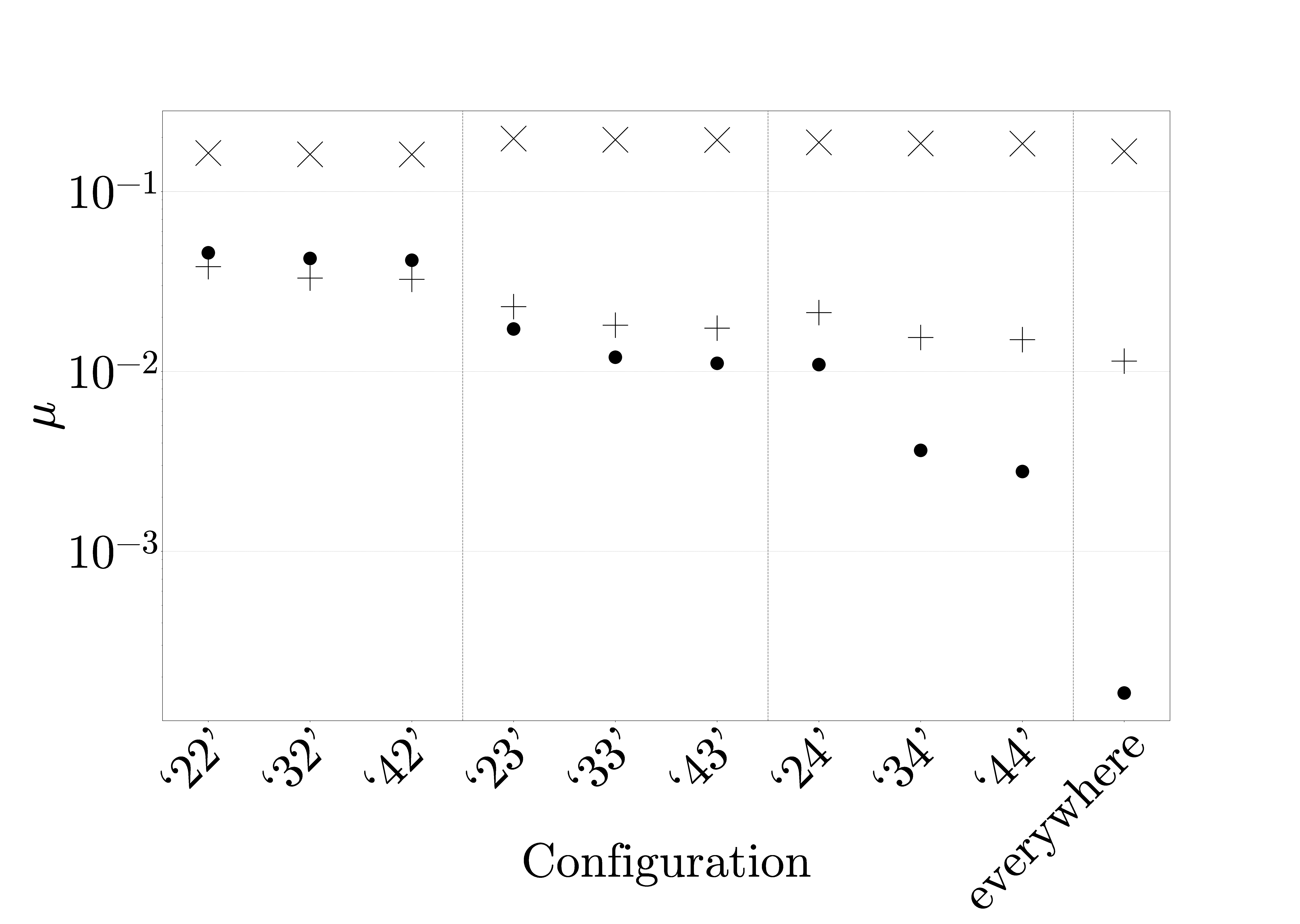}\label{fig:FMA_SameMach_30}}
\subfloat[\ang{40} fixed misalignment w.r.t. grid]{\includegraphics[width=0.5\textwidth]{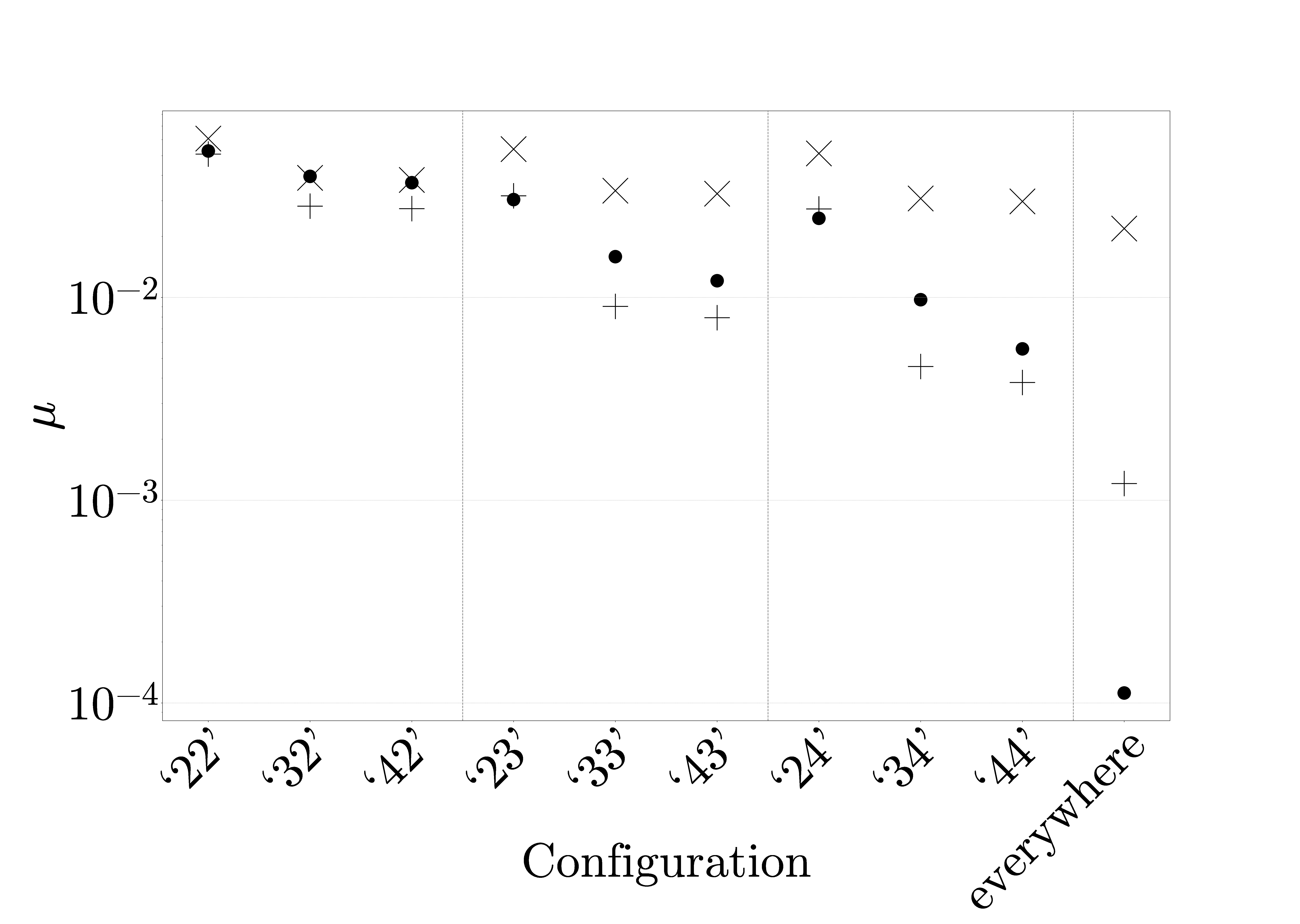}\label{fig:FMA_SameMach_40}}\\[-2ex]
\subfloat[\ang{50} fixed misalignment w.r.t. grid]{\includegraphics[width=0.5\textwidth]{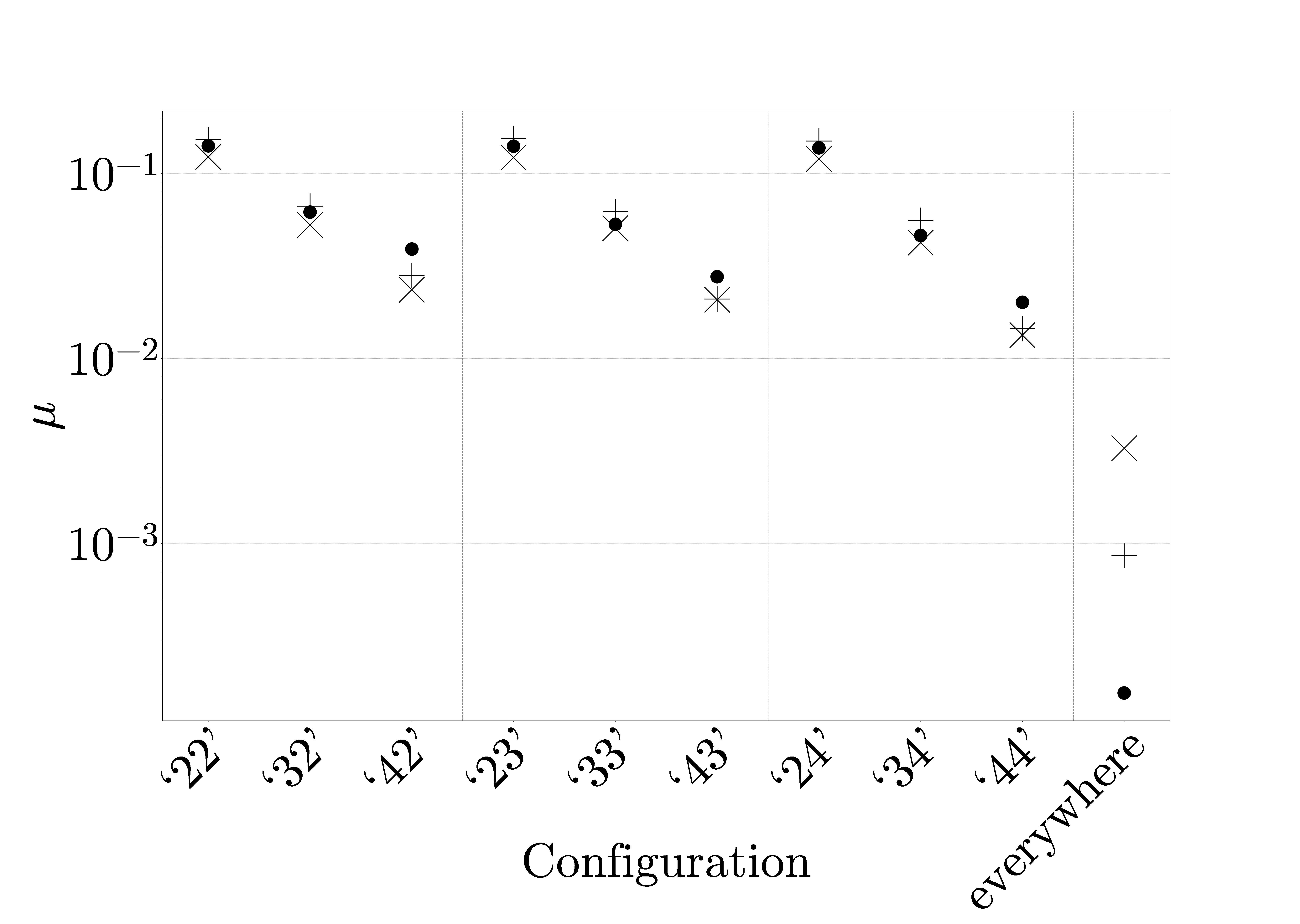}\label{fig:FMA_SameMach_50}}
\caption{Fixed misalignment of the shock with respect to the grid: $\mu$ for three different shock angles for a given fixed misalignment of the shock with respect to the grid given in Table (\ref{tab:FixedMisAlignment_SameMach}). ($\mu = \text{TV} - L_{\infty}$)}
\label{fig:FMA_SameMach_TvLinfDiff}
\end{figure}

Results are presented in Figure (\ref{fig:FMA_SameMach_TvLinfDiff}) for three different fixed nonalignment angles with respect to the grid. By examining the results of the \ang{50} shock nonalignment with respect to the grid in Figure (\ref{fig:FMA_SameMach_50}), we observe that for all three shock angles, the effect of the troubled-cells is the same i.e., adding extra troubled-cells in the pre-shock region improves the solution, but not in the post-shock region, although this is quite the opposite for smaller shock angles, as seen earlier. This indicates that the dependence of number of troubled-cells is actually on the shock angle with respect to the grid rather than the actual shock angle.

The same conclusions hold for the other two fixed nonalignment angles. However, when the actual shock angle is higher than the shock angle with respect to the grid, the solutions exhibits oscillations regardless of the number of troubled-cells. This can be verified from Figure (\ref{fig:FMA_SameMach_30}), where a \ang{50} shock angle is represented as a \ang{30} shock with respect to the grid. In this case, the $\mu$ values are higher and remain same, independent of the number of the troubled-cells.

Next, we present the results for flow over a ramp with an aligned shock test case illustrated in Figure (\ref{fig:AlignedRamp_setup}). Since the shock is aligned to the grid, pre-labelling the troubled-cells is similar to the aligned oblique shock case.

\begin{figure}
\centering
\includestandalone[width=0.25\textwidth]{./6-Results/AlignedRamp/SameMach/legend}\\[-5ex]
\includegraphics[width=\textwidth]{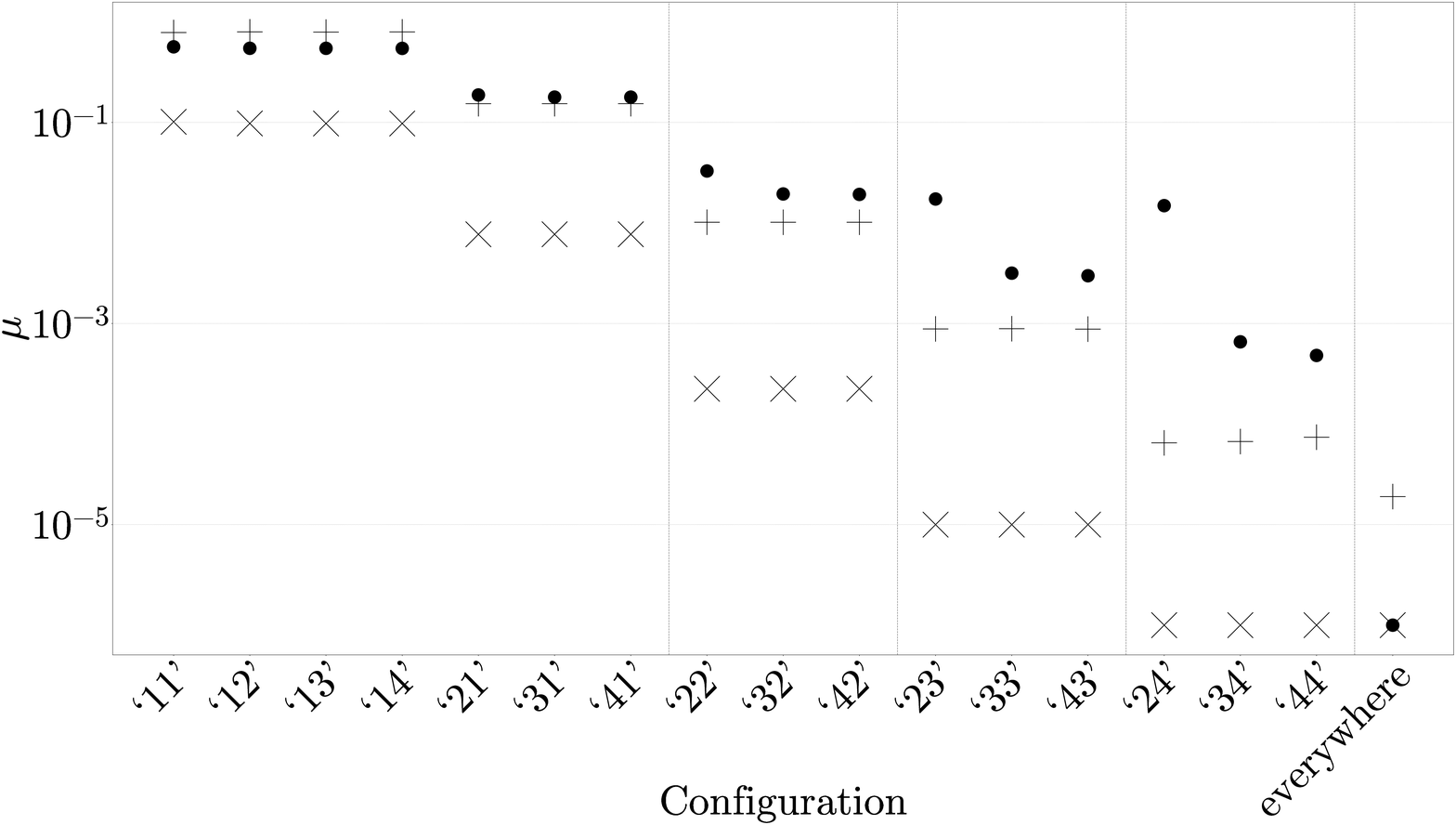}\label{fig:SameMachDensity}
\caption{Flow over a ramp with inlet Mach 3. $\mu$ for three different shock angles across all configurations, including the limiting-everywhere case. ($\mu = \text{TV} - L_{\infty}$)}
\label{fig:AR_SameMach_TvLinfDiff}
\end{figure}

Figure (\ref{fig:AR_SameMach_TvLinfDiff}) presents the $\mu$ for all configurations, including limiting everywhere approach. These results exhibit the same behaviour as aligned shocks mentioned earlier, since the shock for this case is also aligned to the grid.
However, the solution did not fully converge and stalled at a level higher than $10^{-14}$, as seen in Figure (\ref{fig:AR_RN}) for both the best configuration for that shock angle and the limiting everywhere approach. The skewness of the grid might be the reason for this convergence stalling, and as mentioned earlier, we are not going to delve deeper into it.

\begin{figure}
\centering
\includegraphics[width=0.7\linewidth]{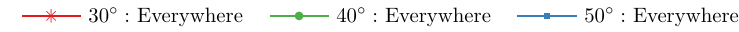}
\includegraphics[width=0.55\linewidth]{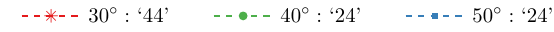}\\
\includegraphics[width=0.6\linewidth]{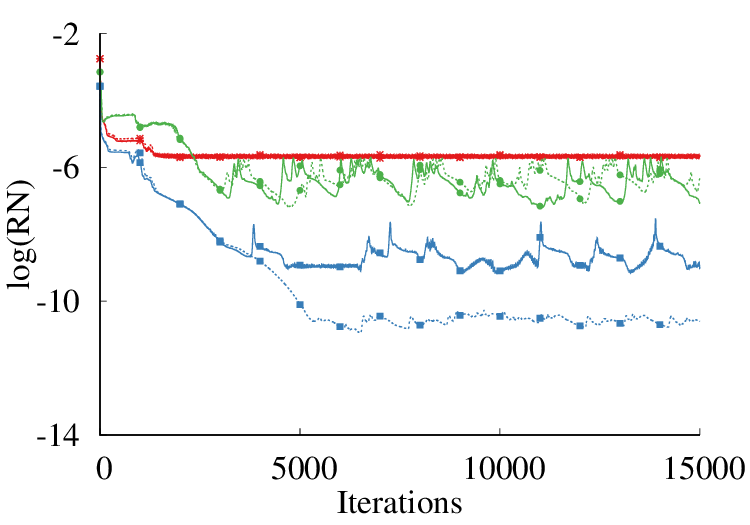}\label{fig:AR_RN}
\caption{Flow over a ramp with inlet Mach 3. The convergence history of the residual norm as a function of number of iterations for both the best configuration and the limiting everywhere approach for each angle.}
\label{fig:AR_RN}
\end{figure}

\section{Comparison of the troubled-cell indicator with the optimal configuration}
\label{sec:Comparison}

In this section, we compare the results of the limiting restricted region approach for two different sets of troubled-cells: (a) identified by the troubled-cell indicator given in equation (\ref{eq:TCI}) for a threshold constant $K = 0.05$ (b) identified by the optimal configuration i.e., `33' for the aligned shocks, `44' for the non-aligned shocks.

\begin{figure}
\centering
\subfloat[\ang{30}]{\includegraphics[width=0.3\linewidth]{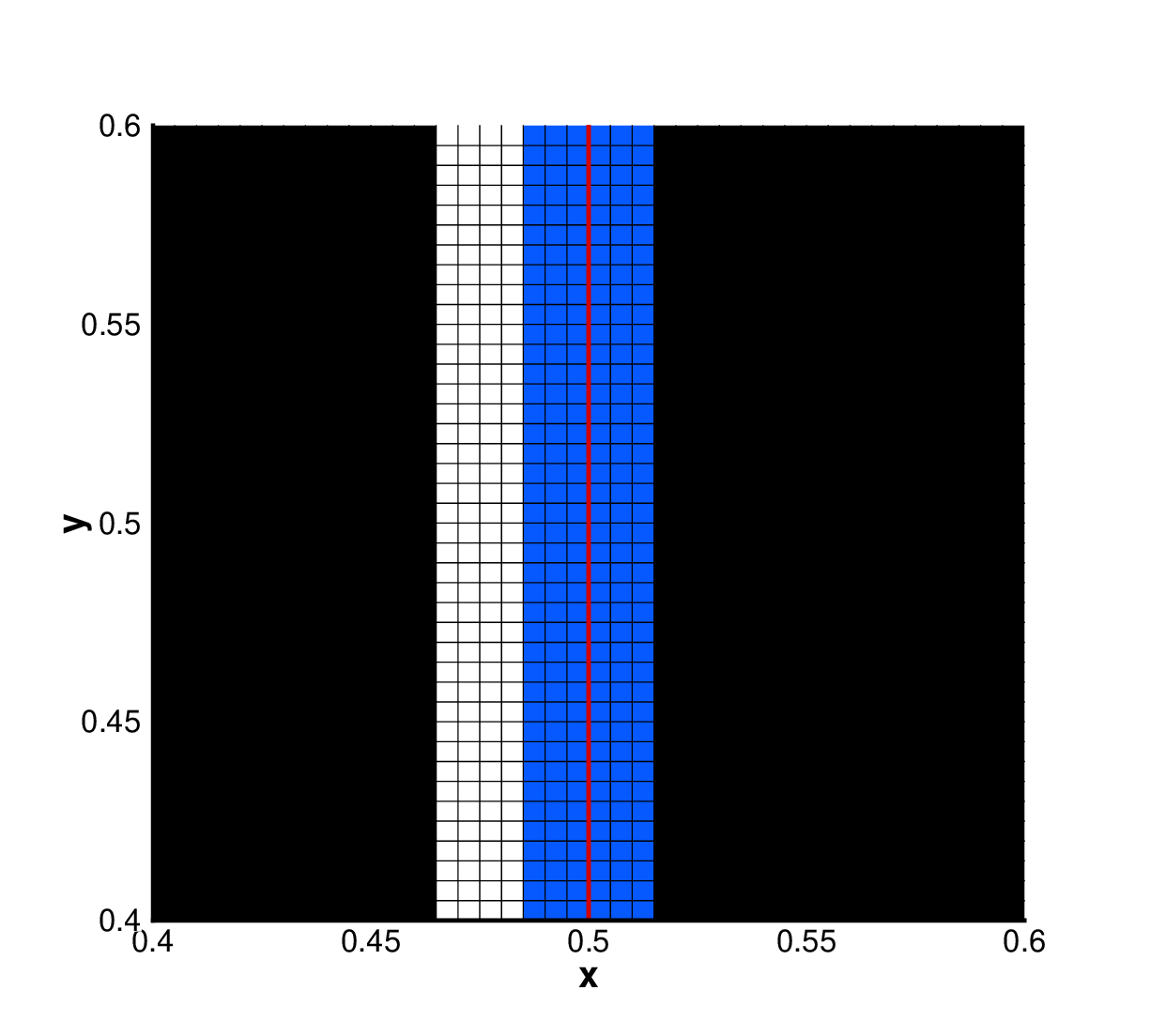}\label{fig:AOS_comp_TC_30}}
\subfloat[\ang{40}]{\includegraphics[width=0.3\linewidth]{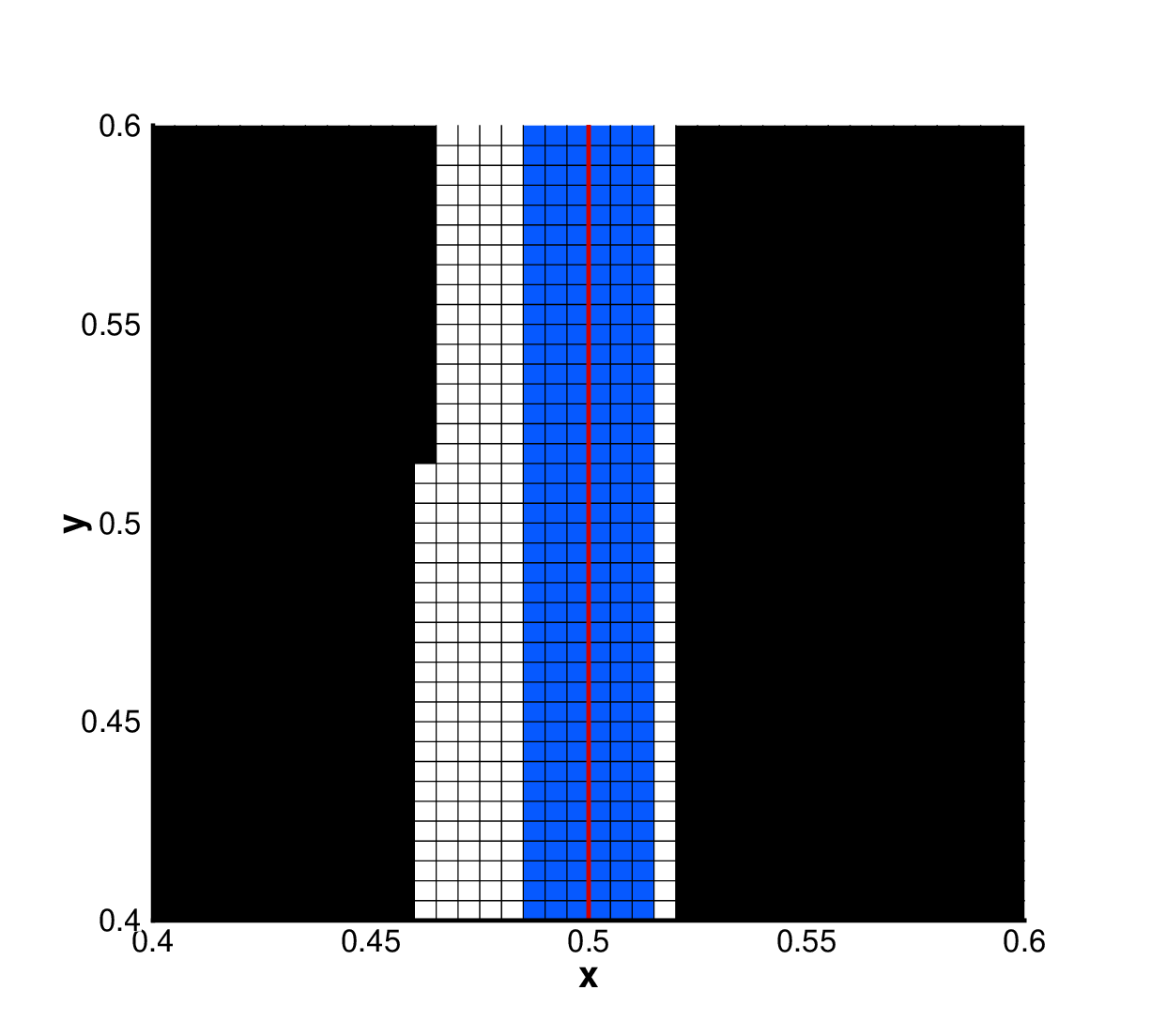}\label{fig:AOS_comp_TC_40}}
\subfloat[\ang{50}]{\includegraphics[width=0.3\linewidth]{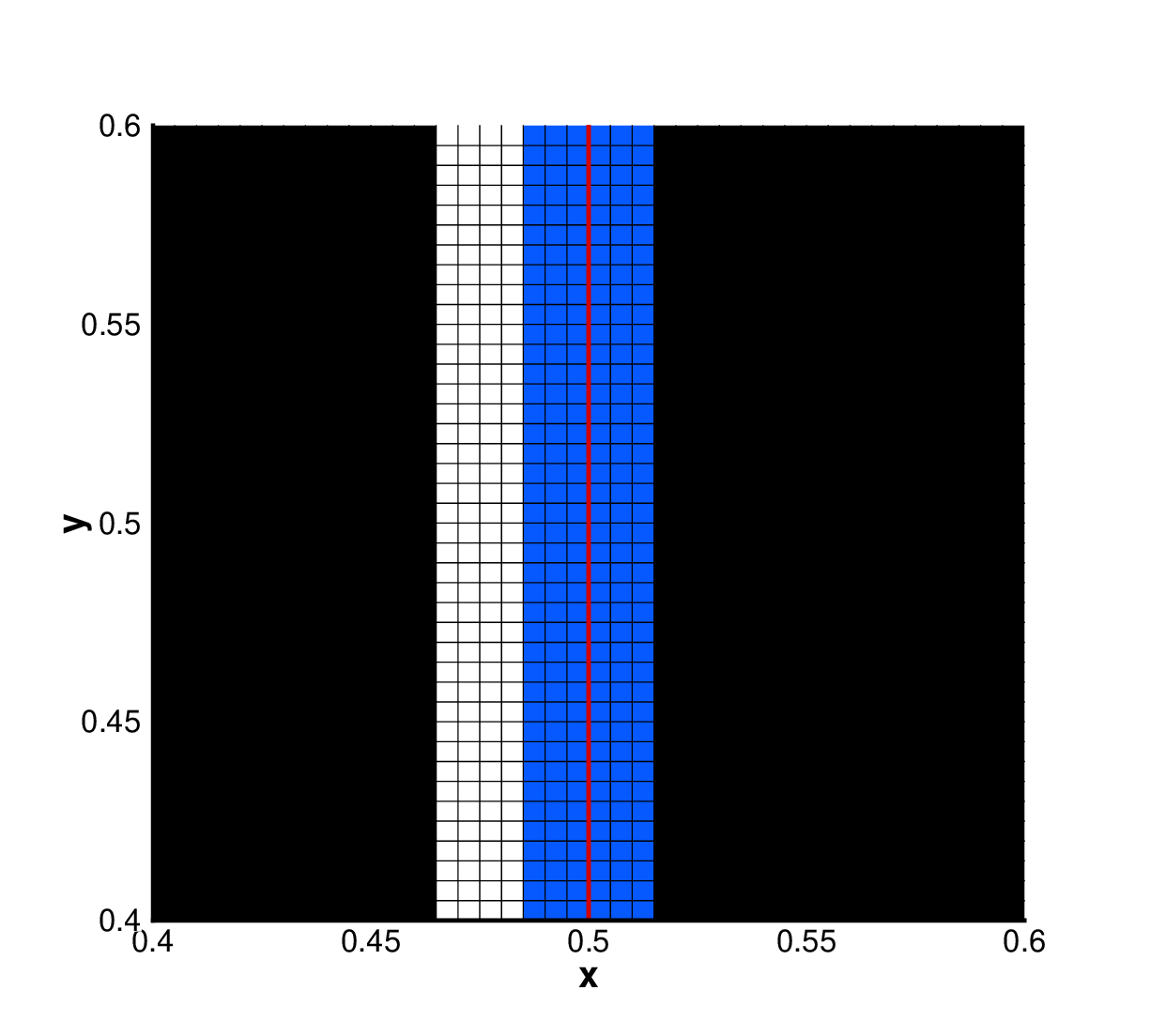}\label{fig:AOS_comp_TC_50}} \\
\subfloat[\ang{30}]{\includegraphics[width=0.3\linewidth]{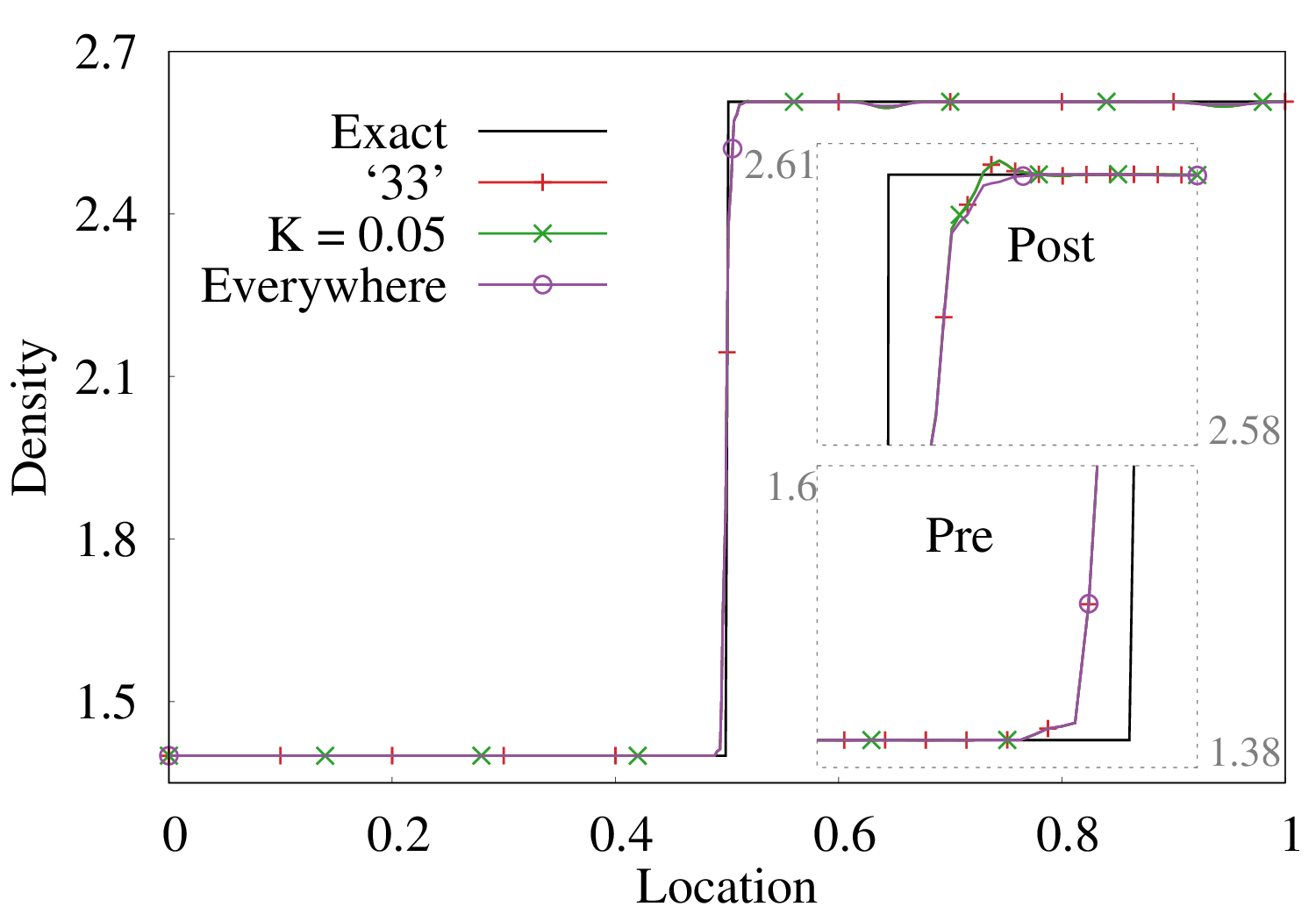}\label{fig:AOS_comp_density_30}}
\subfloat[\ang{40}]{\includegraphics[width=0.3\linewidth]{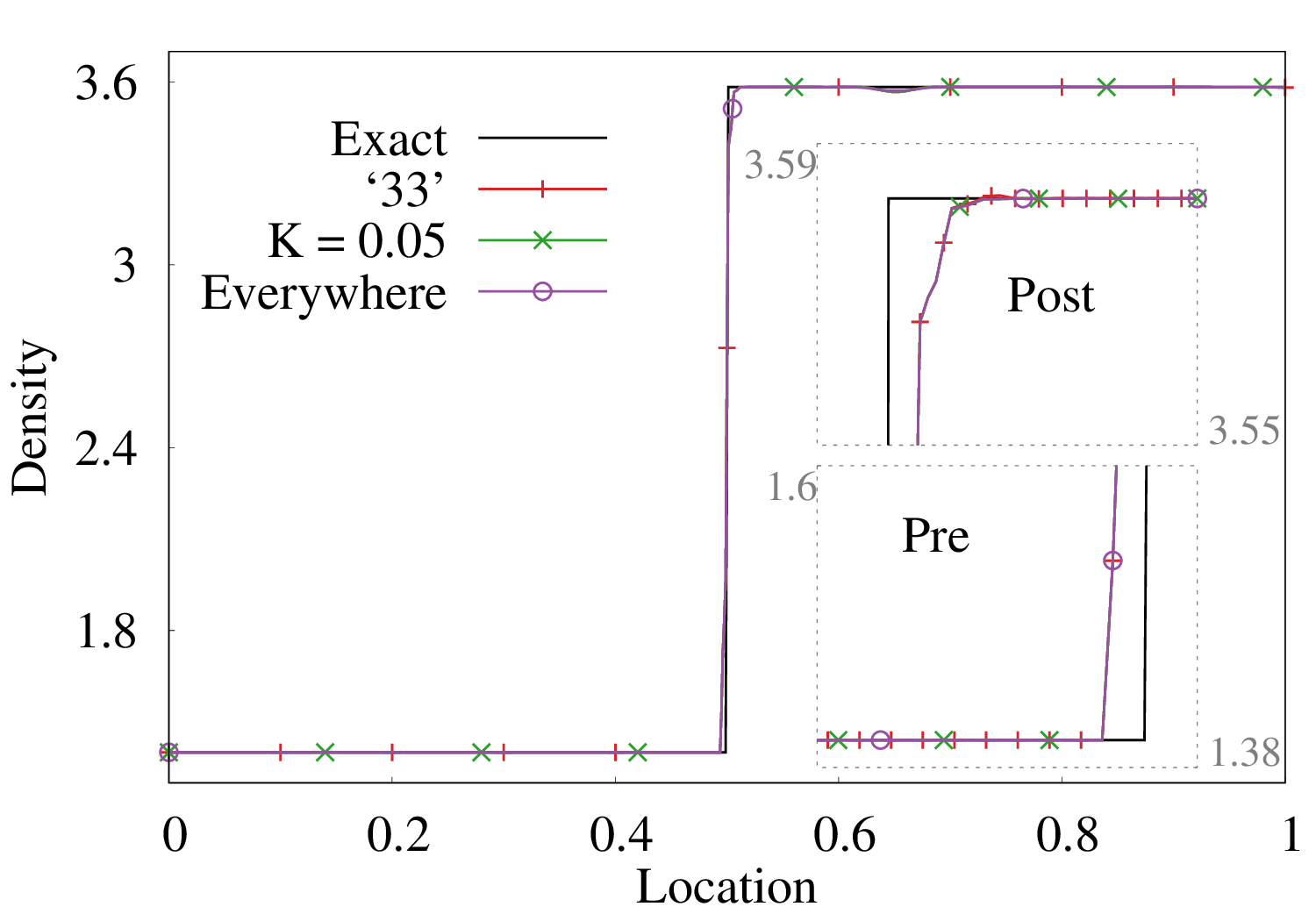}\label{fig:AOS_comp_density_40}}
\subfloat[\ang{50}]{\includegraphics[width=0.3\linewidth]{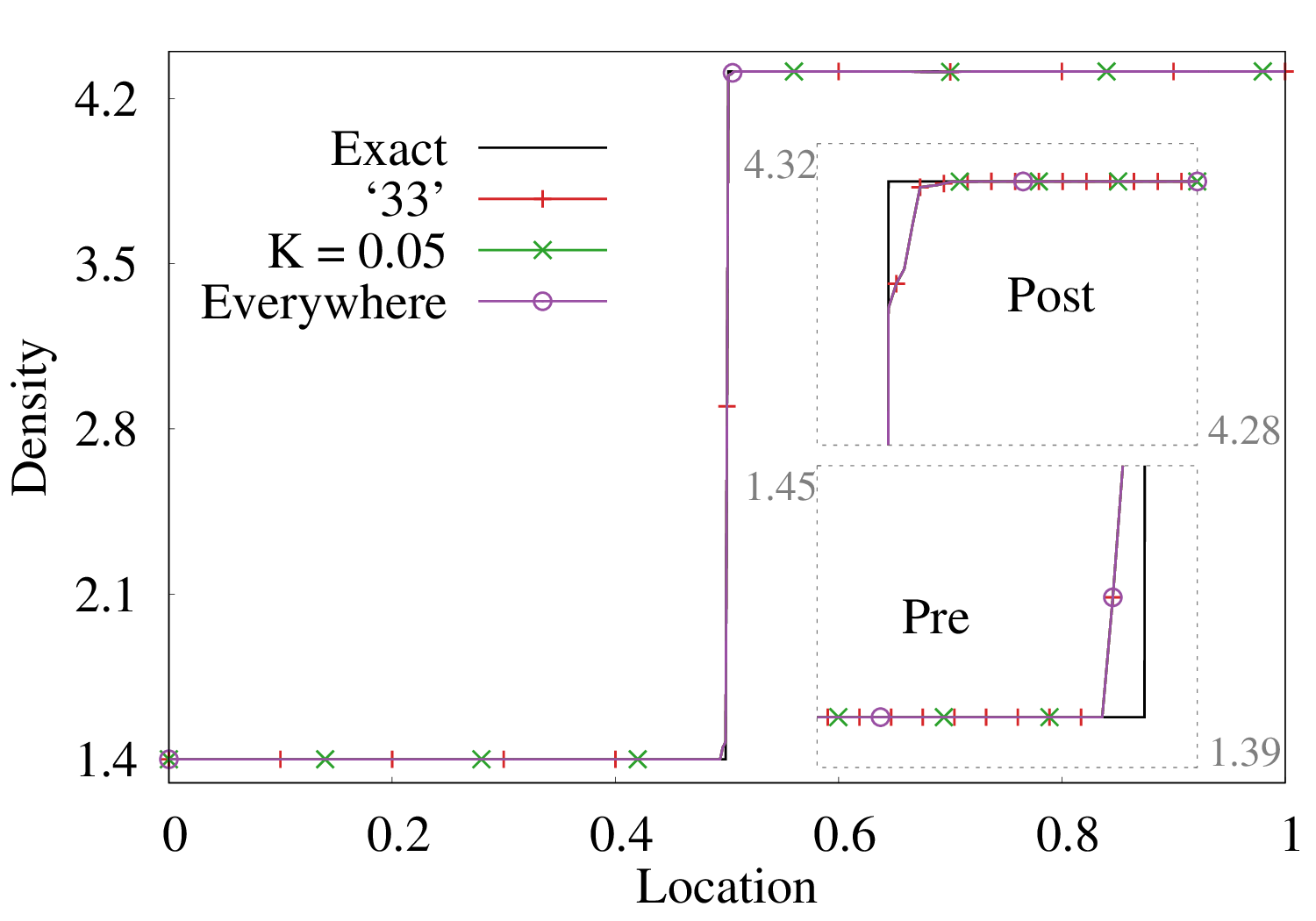}\label{fig:AOS_comp_density_50}} \\
\subfloat[\ang{30}]{\includegraphics[width=0.3\linewidth]{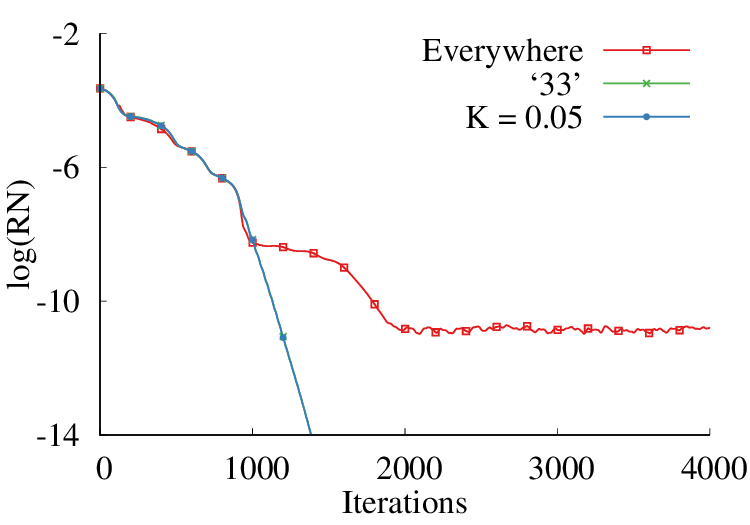}\label{fig:AOS_comp_RN_30}}
\subfloat[\ang{40}]{\includegraphics[width=0.3\linewidth]{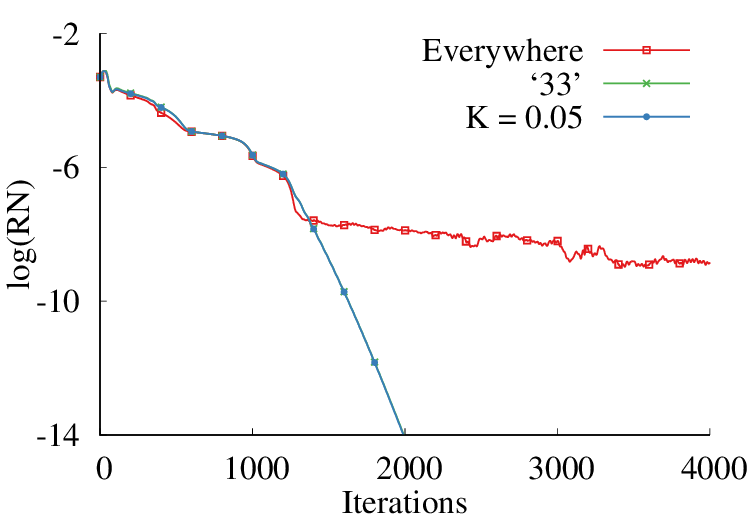}\label{fig:AOS_comp_RN_40}}
\subfloat[\ang{50}]{\includegraphics[width=0.3\linewidth]{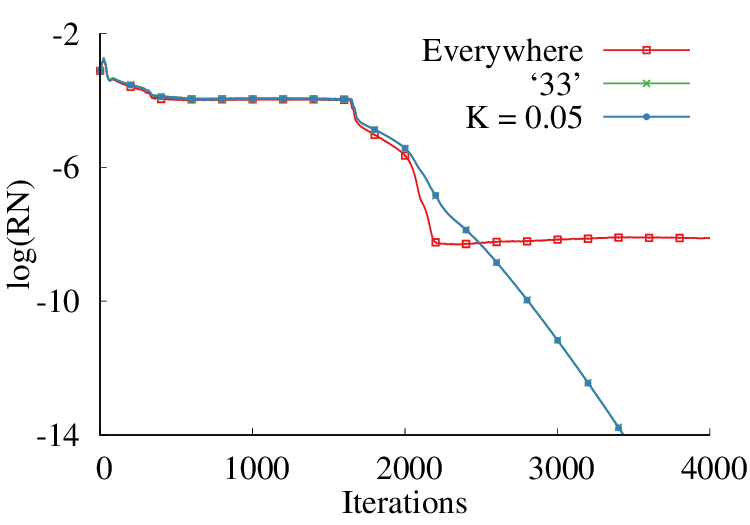}\label{fig:AOS_comp_RN_50}} \\
\caption{Aligned oblique shock. Top row (a, b, c): Zoomed-in view of troubled-cells. Blue color cells are identified by both $K = 0.05$ and configuration `33'. White color cells are identified by the indicator for $K = 0.05$. Red line represents the exact shock. Middle row (d, e, f): Density profiles along the line $y = 0.5$. Bottom row (g, h, i): The convergence history of the residual norm as a function of number of iterations.}
\label{fig:Aligned_Compare}
\end{figure}

\begin{figure}
\centering
\subfloat[\ang{30}]{\includegraphics[width=0.3\linewidth]{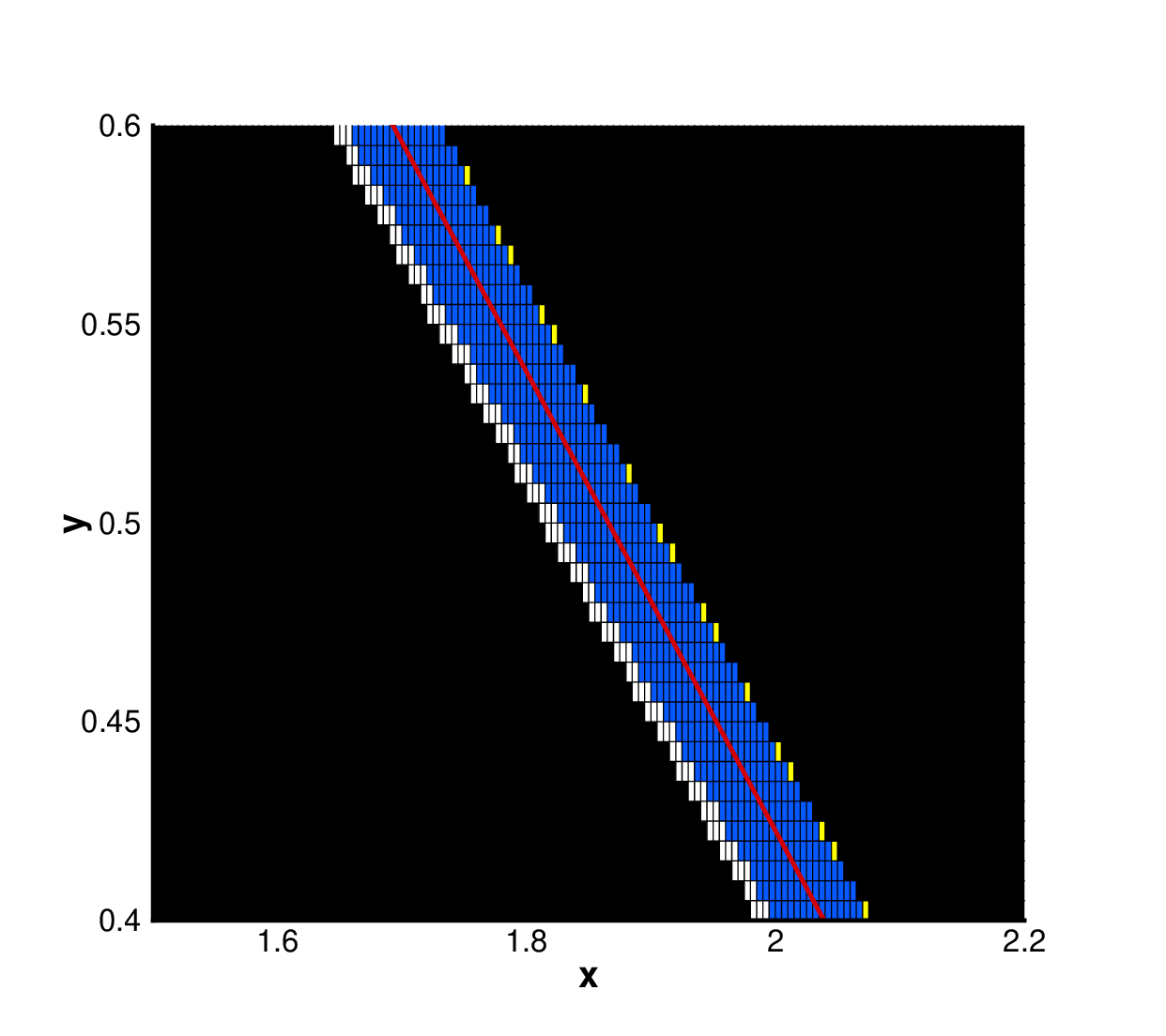}\label{fig:NAOS_comp_TC_30}}
\subfloat[\ang{40}]{\includegraphics[width=0.3\linewidth]{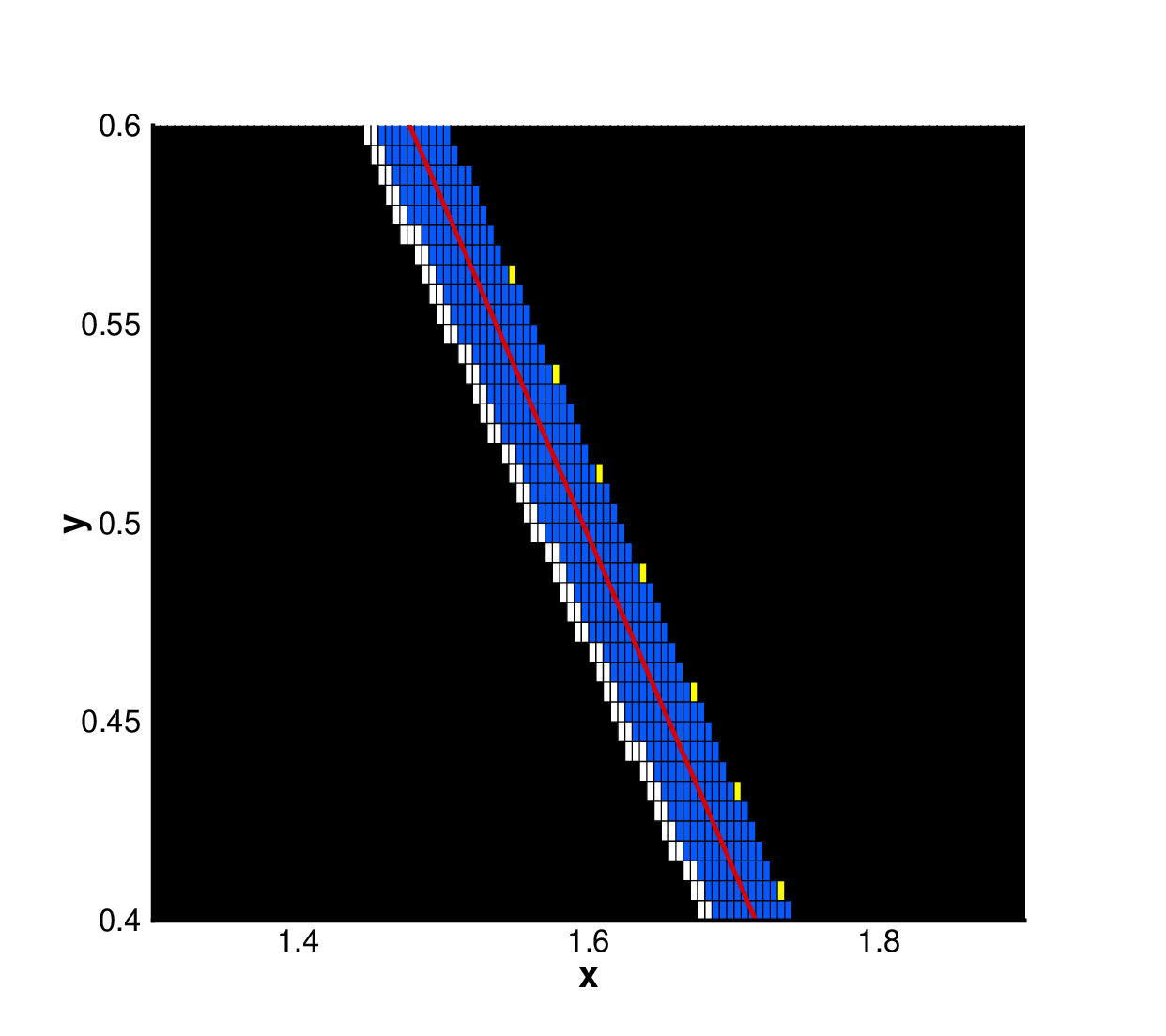}\label{fig:NAOS_comp_TC_40}}
\subfloat[\ang{50}]{\includegraphics[width=0.3\linewidth]{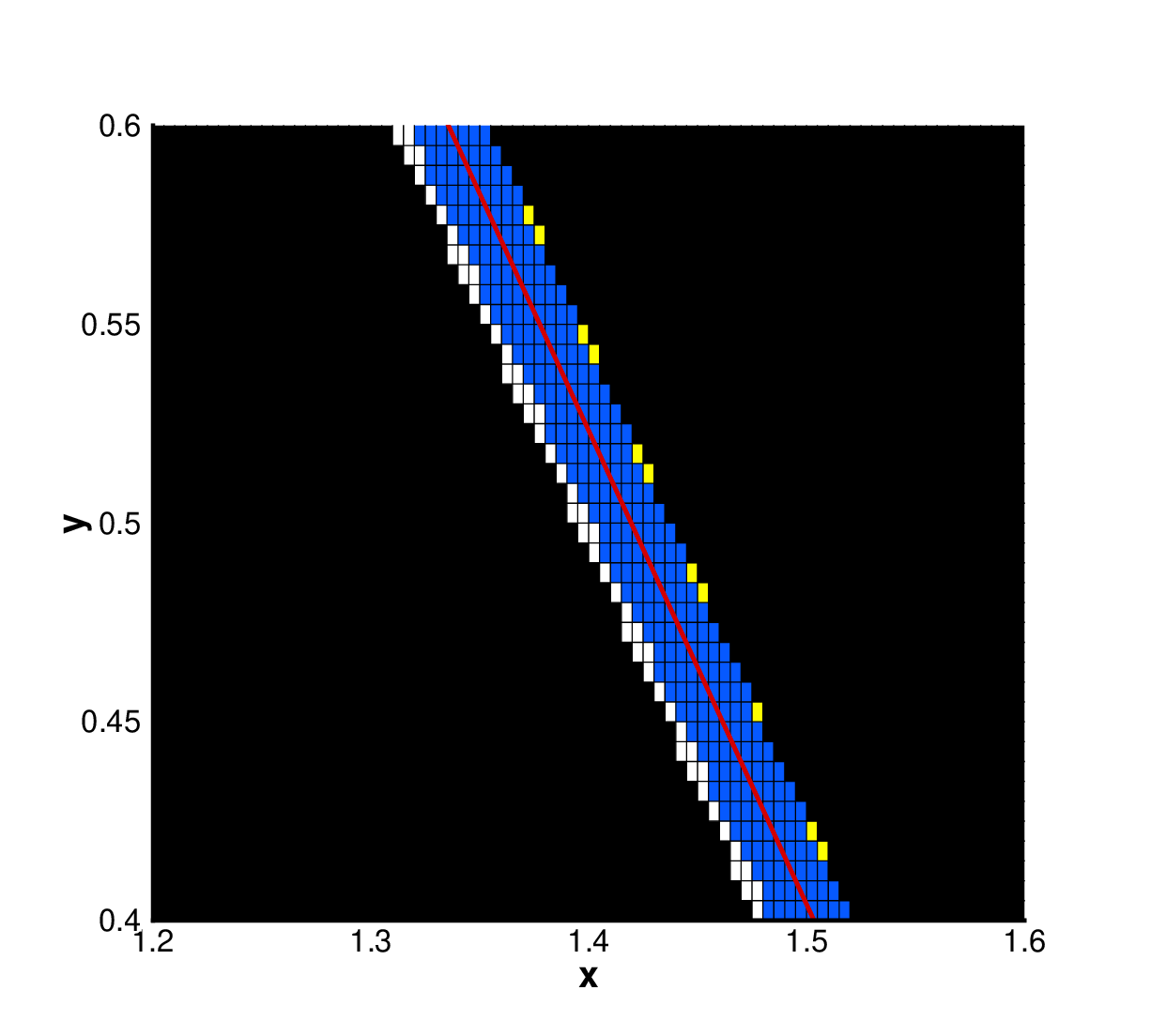}\label{fig:NAOS_comp_TC_50}} \\
\subfloat[\ang{30}]{\includegraphics[width=0.3\linewidth]{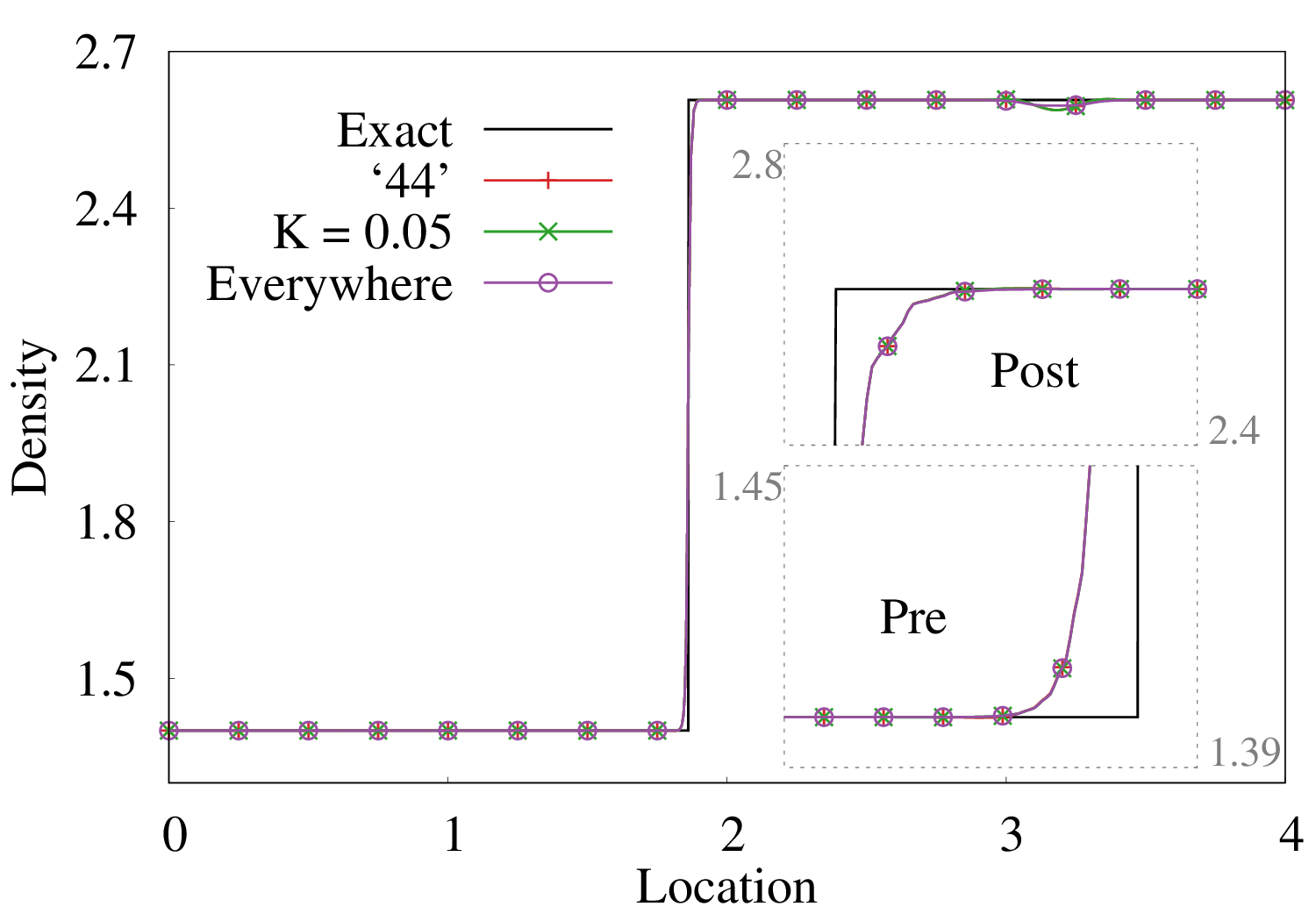}\label{fig:NAOS_comp_density_30}}
\subfloat[\ang{40}]{\includegraphics[width=0.3\linewidth]{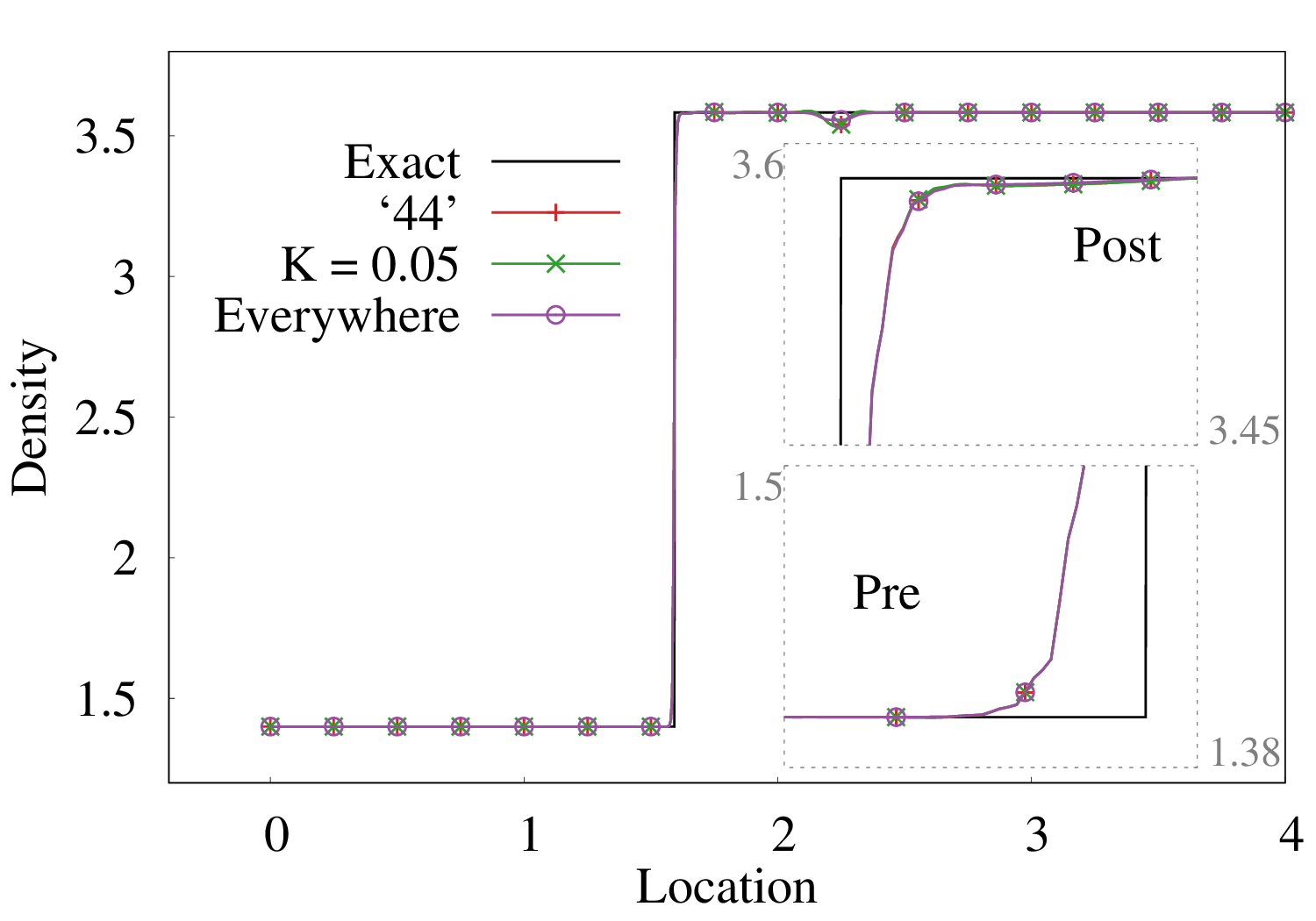}\label{fig:NAOS_comp_density_40}}
\subfloat[\ang{50}]{\includegraphics[width=0.3\linewidth]{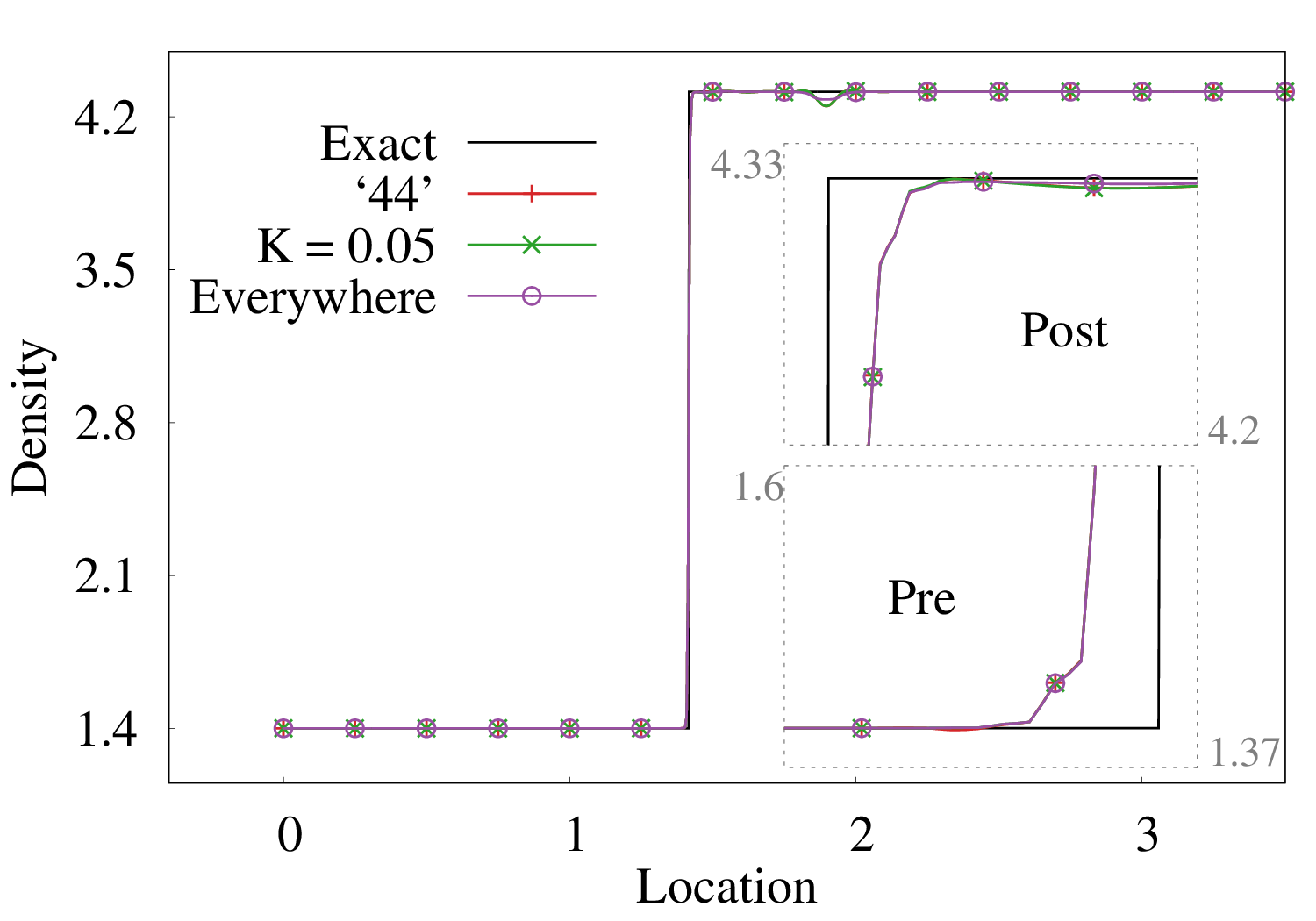}\label{fig:NAOS_comp_density_50}} \\
\subfloat[\ang{30}]{\includegraphics[width=0.3\linewidth]{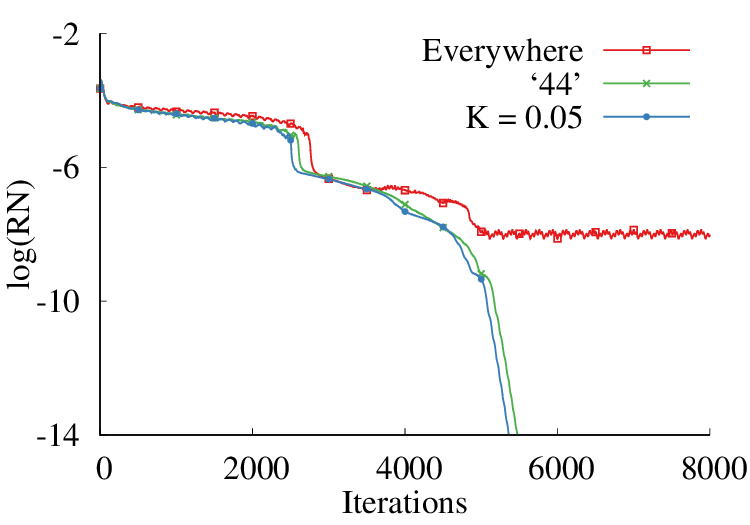}\label{fig:NAOS_comp_RN_30}}
\subfloat[\ang{40}]{\includegraphics[width=0.3\linewidth]{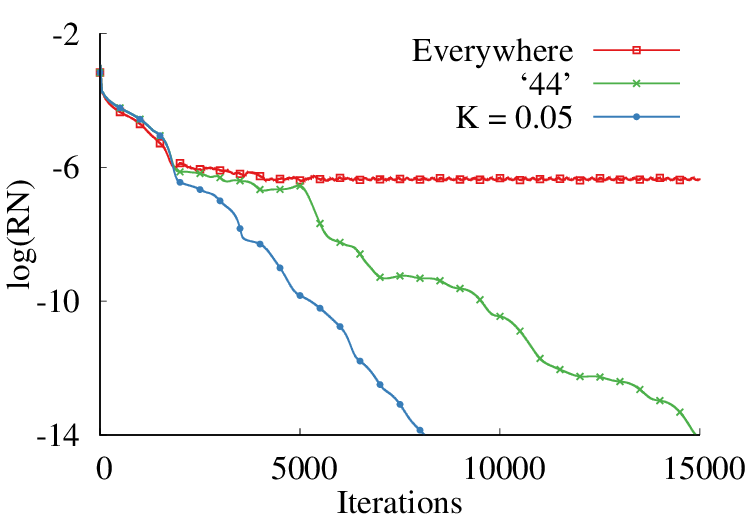}\label{fig:NAOS_comp_RN_40}}
\subfloat[\ang{50}]{\includegraphics[width=0.3\linewidth]{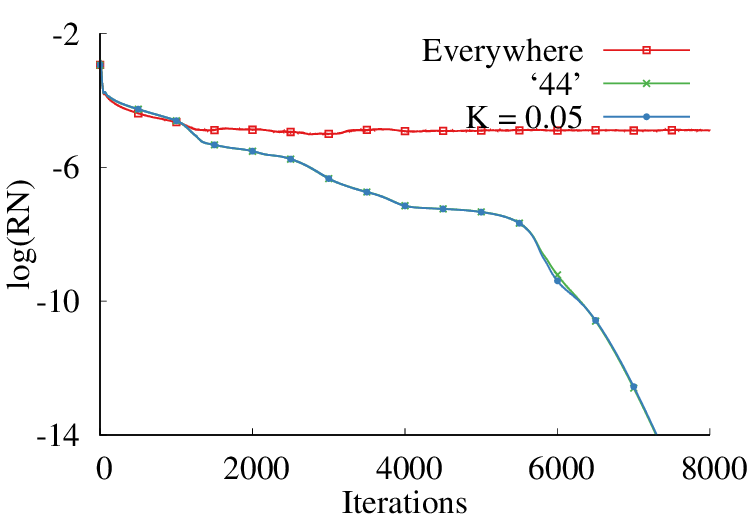}\label{fig:NAOS_comp_RN_50}} \\
\caption{Nonaligned oblique shock. Top row (a, b, c): Zoomed-in view of troubled-cells. Blue color cells are identified by both $K = 0.05$ and configuration `44'. White color cells are identified by the indicator for $K = 0.05$. Yellow color cells are identified by the configuration `44'. Red line represents the exact shock. Middle row (d, e, f): Density profiles along the line $y = 0.5$. Bottom row (g, h, i): The convergence history of the residual norm as a function of number of iterations.}
\label{fig:Nonaligned_Compare}
\end{figure}

For case (a), we use the first-order solutions obtained using the Lax-Friedrichs or AUSM+ flux schemes. For both the aligned and non-aligned shocks, the troubled-cell indicator identifies all the troubled-cells which are identified by the optimal configuration of the respective shocks, except for a few in the post-shock region particularly for the non-aligned shocks. However, the indicator identifies extra troubled-cells in the pre-shock region for the Lax-Friedrichs solution compared to the AUSM+ solution as the earlier solution is more dissipative in nature. For aligned shocks, the indicator identifies only one or two troubled-cells in both pre- and post-shock regions for the AUSM+ solution as it is more accurate than the Lax-Friedrichs solution. Here, we present the results of Lax-Friedrichs solution for the aligned shocks, and AUSM+ solution for non-aligned shocks. The top rows of Figures (\ref{fig:Aligned_Compare}) and (\ref{fig:Nonaligned_Compare}) present the zoomed-in view of troubled-cells identified by the indicator for $K = 0.05$ and the respective optimal configuration for the aligned and non-aligned shocks, respectively.

To obtain the high-order solutions, we initialise the computational domain with the first-order Lax-Friedrich solution for the aligned shocks, and with the first-order AUSM+ solution for non-aligned shocks, along with the information of troubled-cells. The middle rows of Figures (\ref{fig:Aligned_Compare}) and (\ref{fig:Nonaligned_Compare}) present the density profiles along the $y = 0.5$ line for aligned and non-aligned shocks, respectively. The solutions obtained by limiting only in troubled-cells identified by the indicator and by the respective optimal configuration are quite similar and they match with the solutions of limiting everywhere approach closely. The results can be verified by the values of the monotonicity parameter for these solutions which are presented in Table (\ref{tab:mu_Forall}).

The bottom rows of Figures (\ref{fig:Aligned_Compare}) and (\ref{fig:Nonaligned_Compare}) present the convergence history of the residual norms for aligned and non-aligned shocks, respectively. The convergence of residual norm is similar for both the cases of limiting restricted region approach and far better than the limiting everywhere approach.

\begin{table}
\centering
\caption{Overall $L_{\infty}$ norm of the density error, the total variation of the density error and the monotonicity parameter for both limiting approaches.}
\resizebox{\textwidth}{!}{%
\begin{tabular}{c c c c c c c}
\toprule
\makecell{Test case} & \makecell{Angle} & \makecell{$K$} & \makecell{$L_2$} & \makecell{$TV$} & \makecell{$L_{\infty}$} & $\mu$\\
\midrule
\multirow{9}{*}{Aligned shock} & \multirow{3}{*}{\ang{30}} & `33' & 0.110258 & 0.692538 &	 0.687915 & 4.623e-03\\
& & $0.05$ & 0.110257 & 0.692678 & 0.687912 & 4.766e-03\\
 & & \cellcolor{lightgray}Everywhere & \cellcolor{lightgray}0.110428	& \cellcolor{lightgray} 0.689646 & \cellcolor{lightgray} 0.688909 & \cellcolor{lightgray}  7.370e-04\\
\cmidrule(){2-7}
& \multirow{3}{*}{\ang{40}} & `33' & 0.121666 & 0.762058 & 0.760859 & 1.199e-03 \\
& & $0.05$ & 0.121745 & 0.761714 & 0.761382 & 3.320e-04\\
& & \cellcolor{lightgray} Everywhere & \cellcolor{lightgray} 0.121747 &	 \cellcolor{lightgray} 0.761551 & \cellcolor{lightgray} 0.761393 & \cellcolor{lightgray} 1.580e-04 \\
\cmidrule(){2-7}
& \multirow{3}{*}{\ang{50}} & `33' & 0.015424 & 0.093142 & 0.093110 & 3.200e-05\\
& & $0.05$ & 0.015424 & 0.093142 & 0.093110 & 3.200e-05\\
& & \cellcolor{lightgray} Everywhere & \cellcolor{lightgray} 0.015422  &	 \cellcolor{lightgray} 0.093099 & \cellcolor{lightgray} 0.093099 & \cellcolor{lightgray} 0 \\
\midrule
\multirow{9}{*}{Non-aligned shock} & \multirow{3}{*}{\ang{30}} & `44' & 0.227137 & 1.149538 & 1.146761 & 2.777e-03\\
& & $0.05$ & 0.226998 & 1.149307 & 1.146582 & 2.725e-03\\
& & \cellcolor{lightgray} Everywhere & \cellcolor{lightgray} 0.227162	& \cellcolor{lightgray} 1.147028 & \cellcolor{lightgray} 1.146865 & \cellcolor{lightgray} 1.630e-04 \\
\cmidrule(){2-7}
& \multirow{3}{*}{\ang{40}} & `44' & 0.341932 & 2.019927 & 2.016120 & 3.807e-03\\
& & $0.05$ & 0.341941 & 2.019600 & 2.016197 & 3.403e-03\\
&  & \cellcolor{lightgray} Everywhere & \cellcolor{lightgray} 0.342349 & \cellcolor{lightgray} 2.017959 & \cellcolor{lightgray} 2.016751 & \cellcolor{lightgray} 1.208e-03\\
\cmidrule(){2-7}
& \multirow{3}{*}{\ang{50}} & `44' & 0.389068 & 2.665835 & 2.652448 & 1.339e-02 \\
& & $0.05$ & 0.389010 & 2.662210 & 2.652388 & 9.822e-03\\
&  & \cellcolor{lightgray} Everywhere & \cellcolor{lightgray} 0.388665 & \cellcolor{lightgray} 2.654769  & \cellcolor{lightgray} 2.651495 & \cellcolor{lightgray} 3.274e-03 \\
\bottomrule
\end{tabular}%
}
\label{tab:mu_Forall}
\end{table}

\begin{figure}
\centering
\subfloat[]{\includegraphics[width=0.4\textwidth]{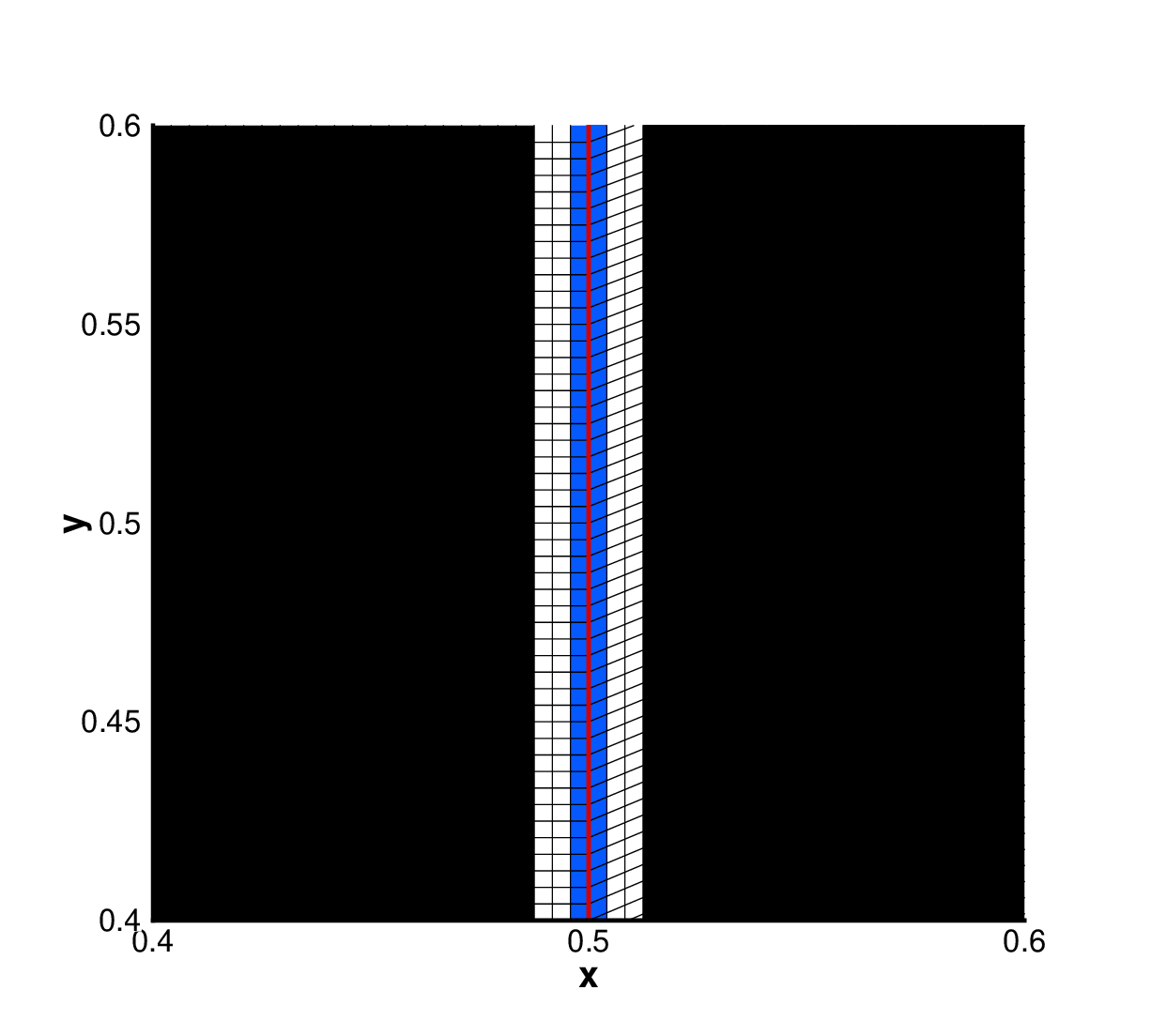}\label{fig:dmr_ramp_tc}}
\subfloat[]{\includegraphics[width=0.4\linewidth]{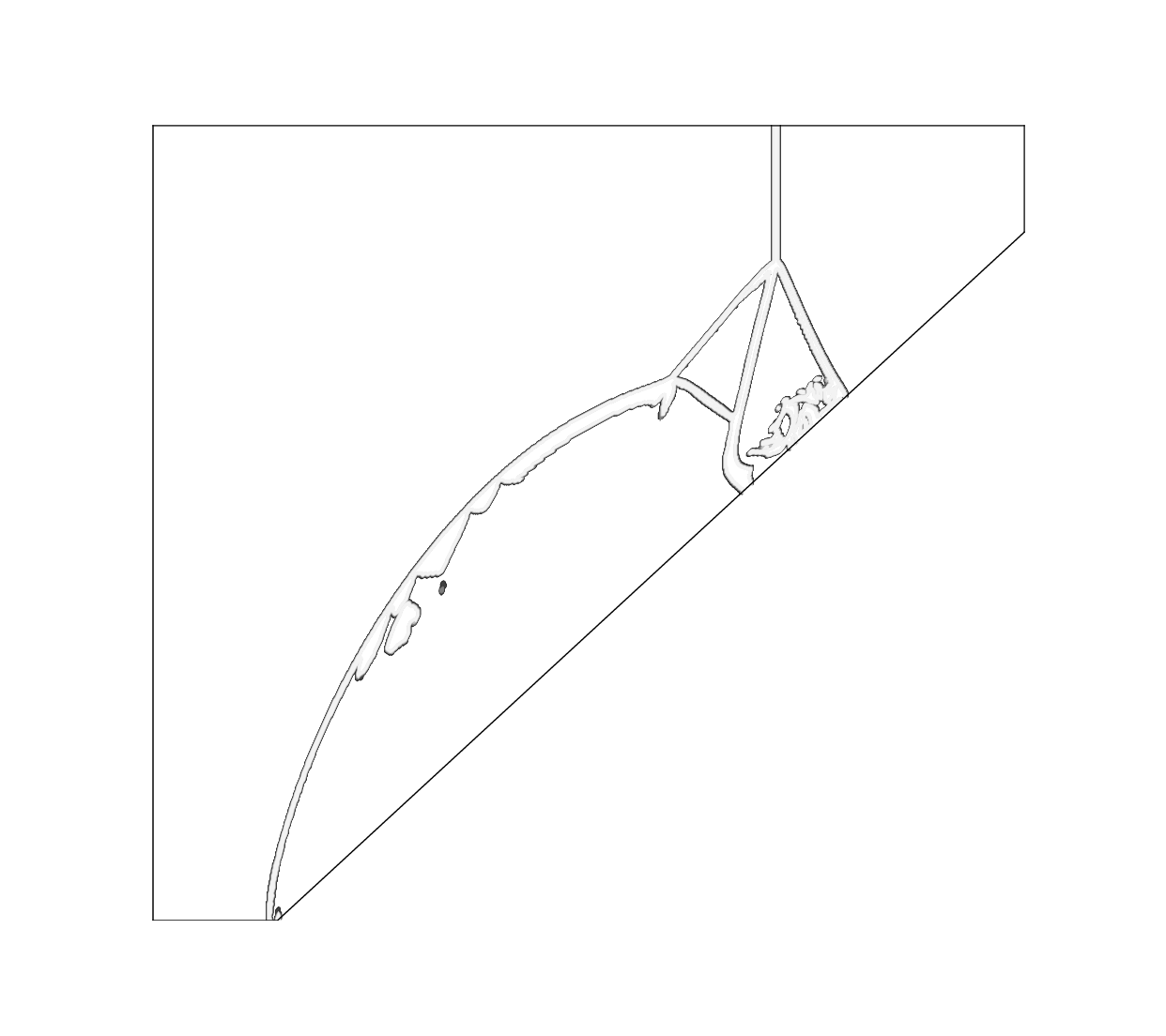}\label{fig:dmr_ramp_tc_final}}\\
\subfloat[]{\includegraphics[width=0.45\textwidth]{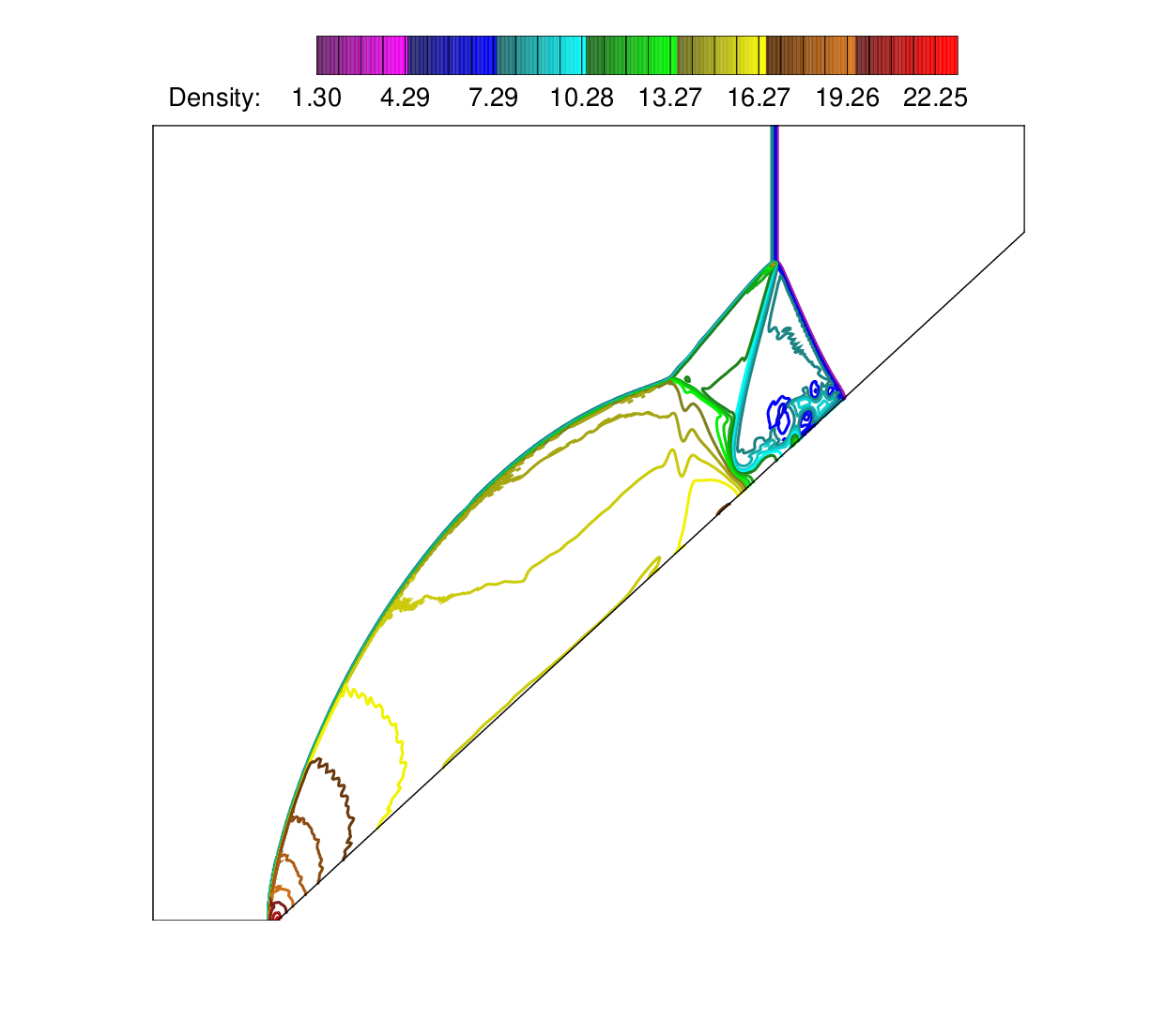}\label{fig:dmr_ramp_density}}
\subfloat[]{\includegraphics[width=0.45\textwidth]{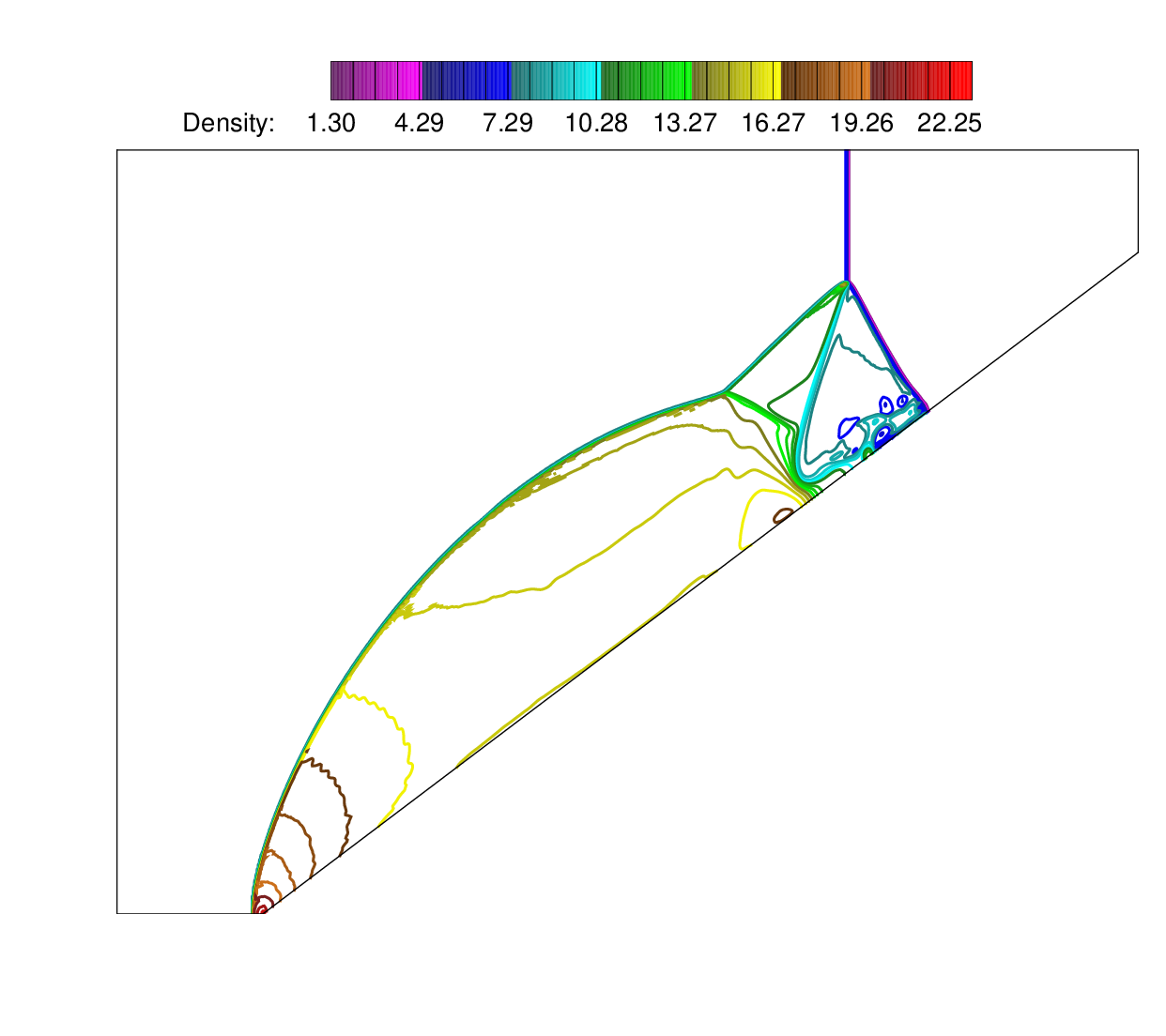}\label{fig:dmr_ramp_density}}
\caption{Double Mach Reflection - Aligned. (a) Zoomed-in view of troubled-cells. Blue color cells are identified by both indicator and optimal configuration. White color cells are identified by the optimal configuration. Red line represents the exact shock. (b)
Troubled-cell region identified by the indicator at the last time step. (c) Density contours of the numerical solutions obtained from limiting restricted region approach. (d) Density contours of the numerical solutions obtained from limiting everywhere approach.}
\label{fig:DMR_Ramp}
\end{figure}

\begin{figure}
\centering
\subfloat[]{\includegraphics[width=0.4\textwidth]{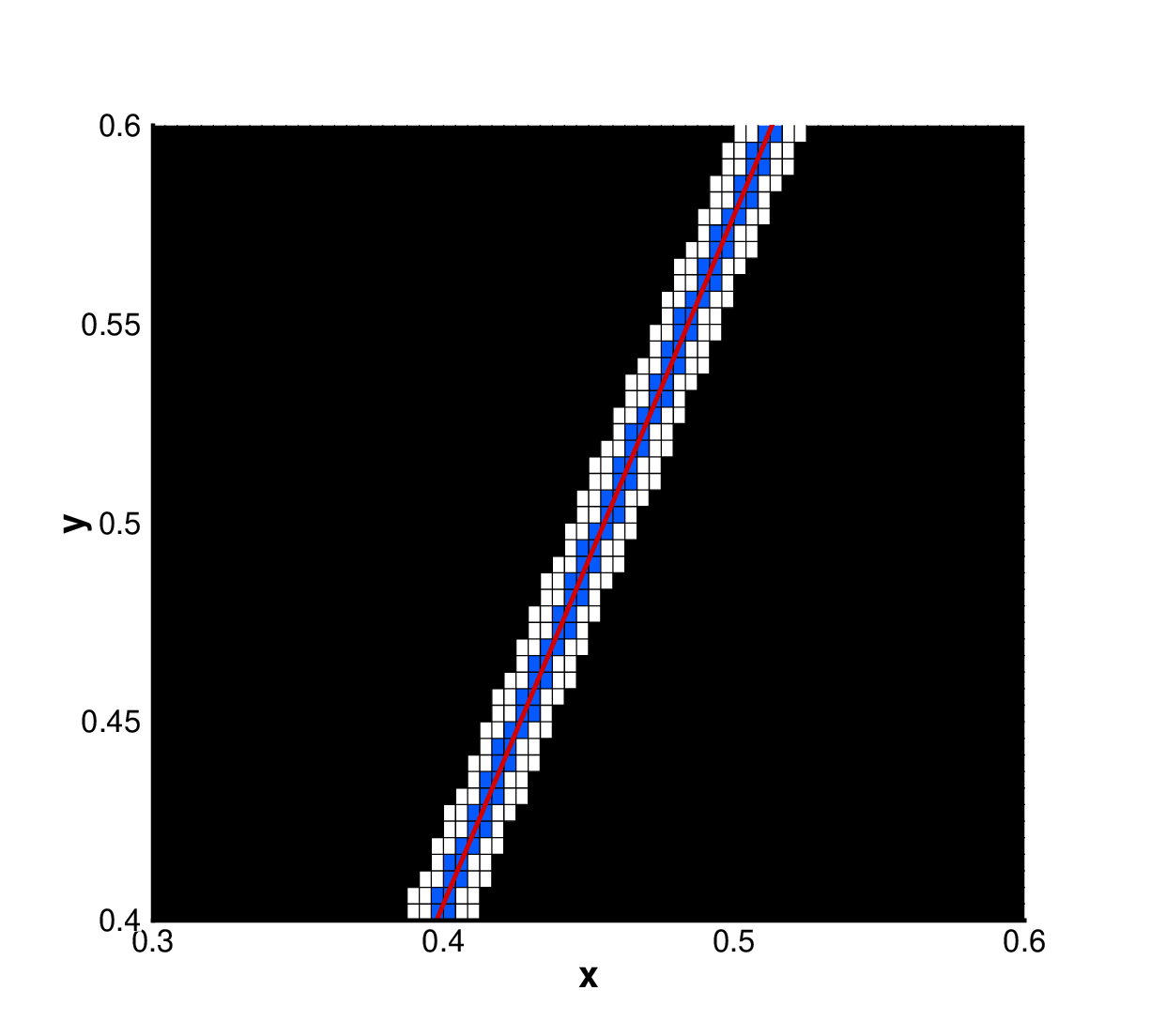}\label{fig:dmr_tc}}
\subfloat[]{\includegraphics[width=0.4\linewidth]{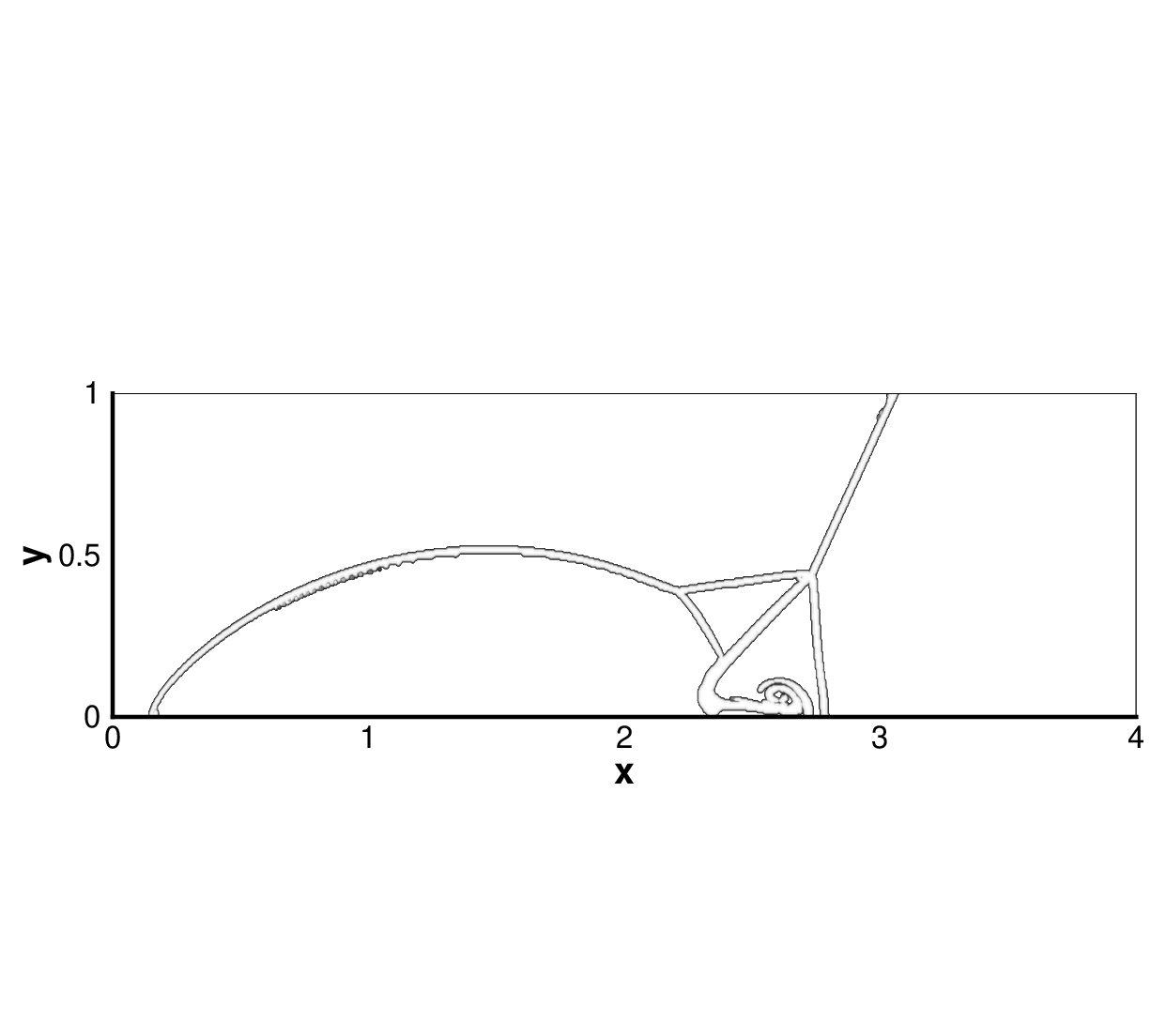}\label{fig:dmr_tc_final}}\\
\subfloat[]{\includegraphics[width=0.45\textwidth]{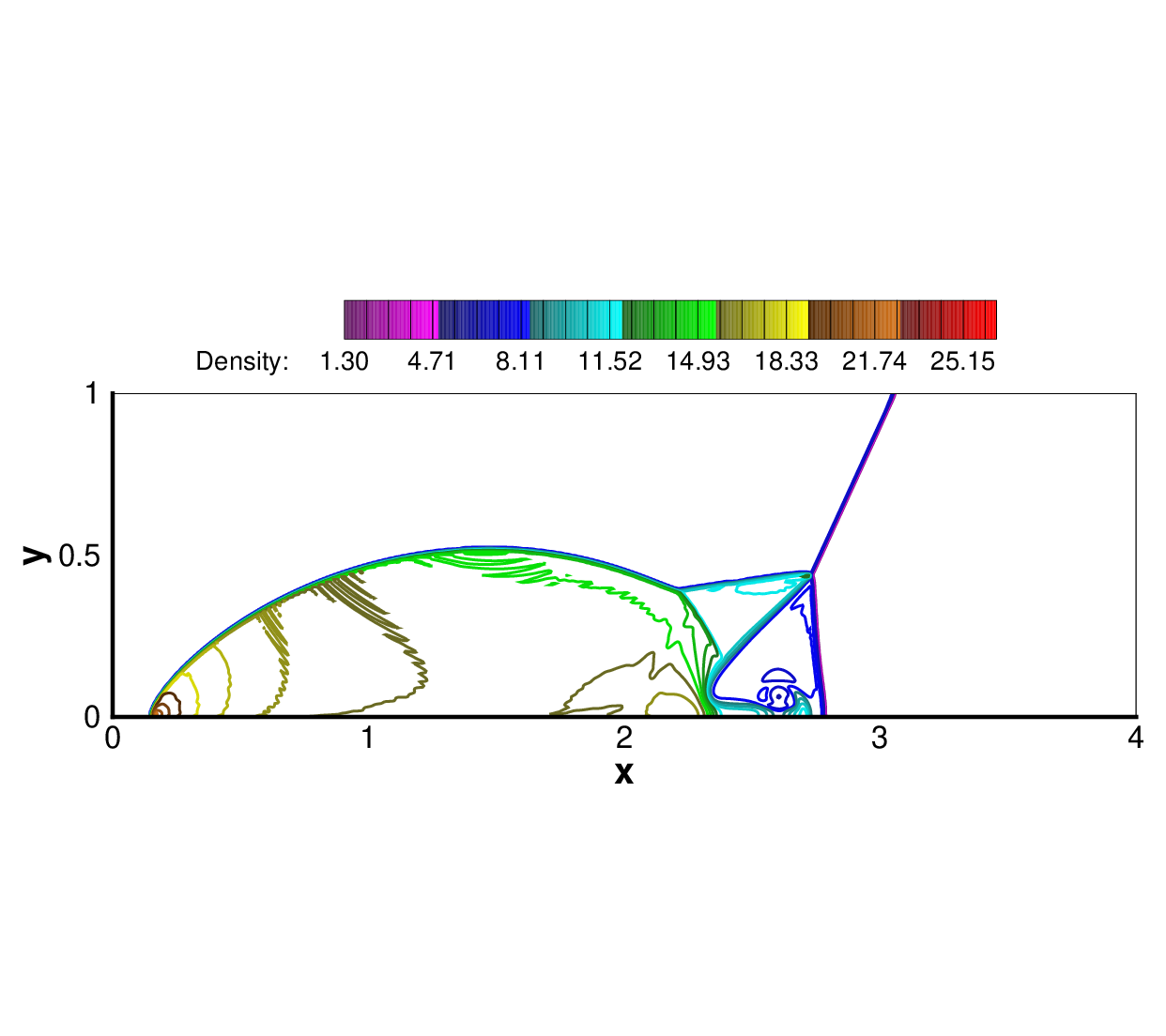}
\label{fig:dmr_density}}
\subfloat[]{\includegraphics[width=0.45\textwidth]{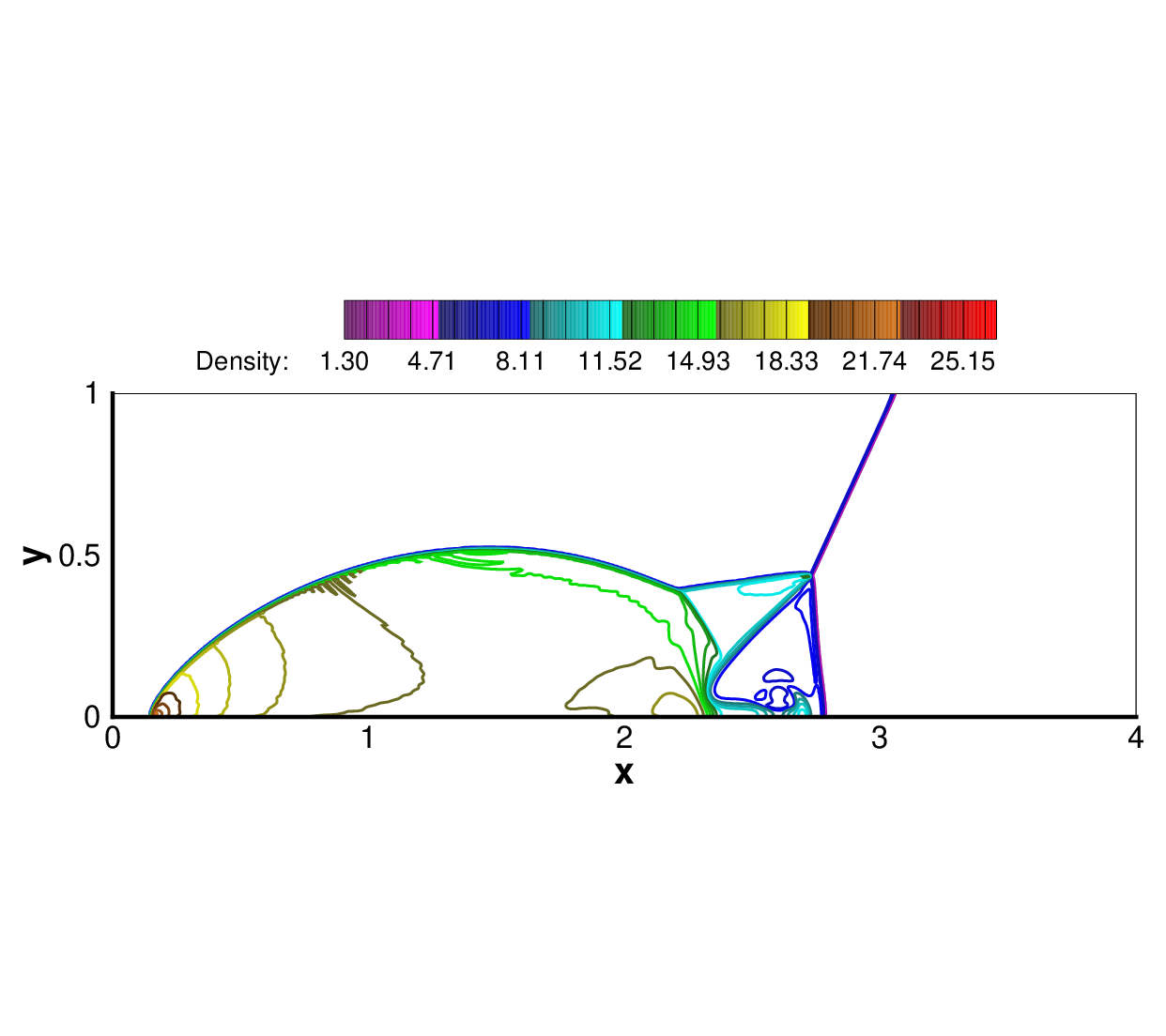}
\label{fig:dmr_density}}
\caption{Double Mach Reflection - Nonaligned. (a) Zoomed-in view of troubled-cells. Blue color cells are identified by both indicator and optimal configuration. White color cells are identified by the optimal configuration. Red line represents the exact shock. (b)
Troubled-cell region identified by the indicator at the last time step. (c) Density contours of the numerical solutions obtained from limiting restricted region approach. (d) Density contours of the numerical solutions obtained from limiting everywhere approach.}
\label{fig:DMR}
\end{figure}

To assess the effectiveness of the troubled-cell indicator for test cases involving complex flow features, we investigate the double Mach reflection for both aligned and non-aligned setups, as shown in Figures (\ref{fig:DMR_ramp_setup}) and (\ref{fig:DMR_setup}), respectively. Here, aligned and nonaligned setups mean that the alignment of the shock with respect to the grid at the initial condition. Since this is a unsteady test case, the troubled-cell indicator need to be applied at every time step to identify the troubled-cells at that particular time step. Figures (\ref{fig:dmr_ramp_tc}) and (\ref{fig:dmr_tc}) present the zoomed-in view of troubled-cells at the first time step for both aligned and non-aligned setups, respectively. In these figures, the cells highlighted in blue represent those identified by the indicator. Notably, only one cell before and after the shock is identified by the indicator because of the exact solution initialisation at $t = 0$. Applying the limiter only in these identified cells resulted in the solution diverging at the end of first time step. As established in the previous section, limiting in fewer number of troubled-cells produces overshoots and overshoot and the higher shock strength worsen these oscillations which ultimately lead to numerical instability.

To address this issue, we use the optimal configuration of aligned (i.e., `33') and non-aligned shocks (i.e., `44') to identify the troubled-cells at the first time step. In Figures (\ref{fig:dmr_ramp_tc}) and (\ref{fig:dmr_tc}), the cells highlighted in both white and blue represent those identified by these respective optimal cases. Subsequently, the troubled-cell indicator is applied at each time step to update the set of troubled cells. By applying the limiter in these troubled cells, a stable solution is successfully obtained.

Figures (\ref{fig:DMR_Ramp}) and (\ref{fig:DMR}) present the density contours of the numerical solutions obtained using both limiting approaches and the troubled-cell region identified by the indicator at the final time step for both aligned and nonaligned setups of double Mach reflection test case, respectively. The results indicate that both the limiting approaches capture the primary features of the solution with comparable accuracy. However, the contour lines for the limiting restricted region approach are less smooth than those observed for the limiting everywhere approach.

Overall, for steady-state test cases, the troubled-cell indicator identifies enough number of troubled-cells to obtain a solution comparable to that of limiting everywhere approach for $K = 0.05$. For the unsteady test cases where the domain is initialised with an exact solution, troubled-cells identified by the indicator may not be sufficient. In such cases, additional troubled cells corresponding to the optimal configuration may be required, or alternatively, the limiting-everywhere approach can be applied during the first time step.

\section{Conclusions}
\label{sec:Conclusion}
In this paper, we introduced a troubled-cell indicator for finite volume methods adapted from discontinuous Galerkin methods. We investigated the optimal number of troubled-cells required in the neighbourhood of a oblique shock to obtain a solution with minimal oscillations and improved convergence. We introduced a novel monotonicity parameter ($\mu$) to quantify the quality of the solution obtained using the limiting restricted region approach.

The results demonstrate that the limiting restricted region approach is better compared to the limiting everywhere approach for finite volume methods employing MUSCL reconstruction, similar to the observations made for discontinuous Galerkin methods. The limiting restricted approach produces the solution with quality similar to that of the limiting everywhere approach with improved convergence, provided there are sufficient number of troubled-cells in the neighbourhood of the shock. However, as the number of troubled-cells decreases, the solutions exhibits overshoots and undershoots.

For aligned shocks, at least two troubled-cells are required on each side of the shock to get a solution without any significant overshoots and undershoots. Increasing the number of troubled cells in the pre-shock region beyond two, while keeping the number of troubled cells in the post-shock region fixed, does not effect the solution, particularly for higher shock angles or higher pressure ratios. However, adding additional troubled-cells in post-shock region, while keeping the number of troubled cells in the pre-shock region fixed, improves the quality of the solution. Overall, three troubled-cells on each side of the shock (i.e., the configuration `33') are sufficient to obtain a solution comparable to that obtained by the limiting everywhere approach, but with improved convergence.

For non-aligned shocks, the troubled-cells identified by tracing the shock and at least two lines parallel to it, on each side of the shock separated by the grid spacing, are required to obtain a solution without any significant overshoots and undershoots. For higher shock angles, more troubled-cells are required in the pre-shock region, while for smaller shock angles, more troubled-cells are required in the post-shock region. Overall, the troubled-cells identified by tracing the shock and four parallel lines on each side of the shock (i.e., the configuration `44') are sufficient to get a solution comparable to that obtained by the limiting everywhere approach, but with improved convergence. The requirement of the number of troubled-cells depends on the shock angle with respect to the grid but not on the actual shock angle with respect to the flow.

The troubled-cell indicator using a threshold constant $K = 0.05$ identified all the troubled-cells of the optimal configuration for the respective test cases, with a few extra in the pre-shock region. Consequently, the resulting solutions closely match and, in turn, closely match with those obtained using the limiting everywhere approach.

For unsteady test cases, when the domain is initialised with the exact solution, only a few troubled-cells are identified by the troubled-cell indicator. In such cases, additional troubled-cells, based on the optimal configuration, need to be incorporated or the limiting everywhere approach must be applied at the initial time step to ensure a stable solution.

\bibliographystyle{ieeetr}
\bibliography{references}

\end{document}